\documentclass{amsart}

\pdfoutput = 1

\usepackage[utf8]{inputenc}
\usepackage[english]{babel}
\usepackage{amsfonts}
\usepackage{amsmath}
\usepackage{graphicx}
\usepackage{array}
\usepackage{subfigure}
\usepackage{geometry}
\title[Opt. penalty meth. for an hyperbolic syst. modeling the edge plasma transp. in a tokamak]{An optimal penalty method for an hyperbolic system modeling the edge plasma transport in a tokamak}

\date\today

\DeclareMathOperator{\minmod}{minmod}
\DeclareMathOperator{\sign}{sign}
\renewcommand{\d}{\operatorname{d}}

\begin{document}
\newtheorem{definition}{Definition}[section]
\newtheorem{theorem}{Theorem}[section]
\newtheorem{proposition}{Proposition}[section]
\newtheorem{lemme}{Lemme}[section]
\newtheorem{corollary}{Corollary}[section]
\newtheorem{assumption}{Assumption}[section]

\renewcommand{\b}[1]{\mathbf{#1}}
\renewcommand{\u}[1]{\underline{#1}}


\author{Philippe Angot, Thomas Auphan, Olivier Gu\`es}

\maketitle
\begin{center}
\small
[angot,tauphan,gues]@cmi.univ-mrs.fr\\
Aix Marseille Universit\'e, CNRS, LATP UMR 7353, 39 rue F. Joliot Curie, 13453 Marseille Cedex 13, France
\end{center}

\begin{abstract}
 The penalization method is used to take account of obstacles, such as the limiter, in a tokamak. Because of the magnetic confinement of the plasma in a tokamak, the transport occurs essentially in the direction parallel to the magnetic field lines. We study a 1D nonlinear hyperbolic system as a simplified model of the plasma transport in the area close to the wall. A penalization which cuts the flux term of the momentum is studied. We show numerically that this penalization creates a Dirac measure at the plasma-limiter interface which prevents us from defining the transport term in the usual distribution sense. Hence, a new penalty method is proposed for this hyperbolic system. For this penalty method, an asymptotic expansion and numerical tests give an optimal rate of convergence without spurious boundary layer. Another two-fields penalization has also been implemented and the numerical convergence analysis when the penalization parameter tends to $0$ reveals the presence of a boundary layer.
 
\end{abstract}

\underline{keywords:} penalization method, nonlinear hyperbolic system, boundary layer, finite volumes, plasma transport, tokamak 


\section{Introduction}

A tokamak is a machine to study plasmas and the fusion reaction induced by the magnetic confinement. The plasma at high temperature ($10^{8} K$, in the center) is confined in a toro\"idal chamber thanks to a magnetic field. One of the main goals is to perform controlled fusion with enough efficiency to be a reliable source of energy.
But, since the magnetic confinement is not perfect, the plasma is in contact with the wall. In order to preserve the integrity of the wall and to limit the pollution of the plasma, it is crucial to control these interactions. Obviously, if it was possible to simulate wall plasma interactions, it would be significantly easier to optimize the configuration.

Plasma models can be classified into three main classes. First, there are the single particle models where we compute the trajectory of each particle, but, as the number of particles in tokamak is of the order of $10^{20}$, the computational cost is prohibitive. The kinetic models study the distribution function $f(\b x, \b v, t)$ which represents the density of particles of speed $\b v$, located in $\b x$, at the time $t$. The kinetic models have seven variables (in a three dimensional spatial domain) so the computational cost is still heavy. The fluid model is the most approximate one, since it considers that the plasma has the same behavior as a fluid and uses equations similar to Navier-Stokes equations. The fluid approximation seems to be verified for the scrape-off layer (temperature of the order of $10^{4} K$) whereas, in the center of the tokamak, a kinetic model is necessary.

To take into account the boundary conditions in the complex geometry of a tokamak, we can use volume penalty methods. These methods consist in embedding the original domain into a fictitious larger and simple domain and to modify the model equations outside the original domain so that the boundary conditions are verified. One advantage of these methods is that we don't need to use a mesh fitted to the geometry of the domain. But, as we shall see later, in addition to well-posedness issue, the penalty method adds a modeling error, which needs to be controlled. Such approaches have already been implemented successfully for elliptic and parabolic systems \cite{Ang05}, for incompressible or compressible flows  \cite{Ang99, Liu07}.

In the sequel, we study, using a fluid approximation of the plasma, a simplified system of equations governing the plasma transport in the scrape-off layer parallel to the magnetic field lines.
A penalty method has been introduced by Isoardi \emph{et al.} \cite{Iso10}, which gives interesting results. But the numerical study was incomplete and the fact that the momentum flux is cut inside the limiter may provide a Dirac measure next to the interface.

In this paper, Section \ref{Presentation of the hyperbolic system} is devoted to a presentation of the toy model considered and Section \ref{Numerical Schemes} gives the finite volume scheme which is used for the numerical tests provided in Sections \ref{First penalty approaches} and \ref{New_penalty_method}. In Section \ref{First penalty approaches},  after a numerical study of the penalization of Isoardi \emph{et al.} \cite{Iso10}, we modify the boundary conditions to ensure the well-posedness of the hyperbolic system and we study numerically another penalization which generates a boundary layer. In Section \ref{New_penalty_method}, we propose an optimal penalty method which is free of boundary layer, a theoretical result is stated for a slightly different problem. At the end of Section \ref{New_penalty_method}, the results of numerical tests are presented, with an extension to a two-sides limiter.

This work completes the first results presented by the authors in \cite{Ang11_2}.

\section{The model hyperbolic problem}\label{Presentation of the hyperbolic system}
At the center of the reactor, the transport along the field lines is almost free of 
constraint and fast enough to consider that at our time scale, physical quantities are constant along a magnetic field line.
This is not the case in the scrape-off layer: magnetic field lines are intercepted
by wall components (such as the limiter in TORE SUPRA). When the ion bumps into the limiter, 
a recombination process occurs and transforms the ion into a neutral particle 
which may be trapped into the limiter or re-injected in the plasma (and re-ionized later).
In this paper, we consider a very simple model taking only into account the transport in the direction parallel 
to the magnetic field lines, (see for example  \cite{Iso10, Tam07}).
It is a one dimensional $2\times 2$ nonlinear hyperbolic system of conservation laws 
for the particle density $N$ and the particle flux $\Gamma$, which reads:
\begin{equation}\label{O1}
\left\{
 \begin{array}{l}
     \partial_t N + \partial_x \Gamma=S_N \\
     \partial_t \Gamma + \partial_x \left(\dfrac{\Gamma^2}{N} + N\right)=S_{\Gamma} \\
     \textsf{Initial conditions: }N(0,.)=N_0 \textsf{ and } \Gamma(0,.)=\Gamma_0,
 \end{array}
    \right.
   (t,x) \in \mathbb{R}^+_* \times ]-L,L[ 
\end{equation}

\begin{figure}
\begin{center}
\includegraphics[scale=0.60, trim = 20mm 91mm 20mm 118mm, clip=true]{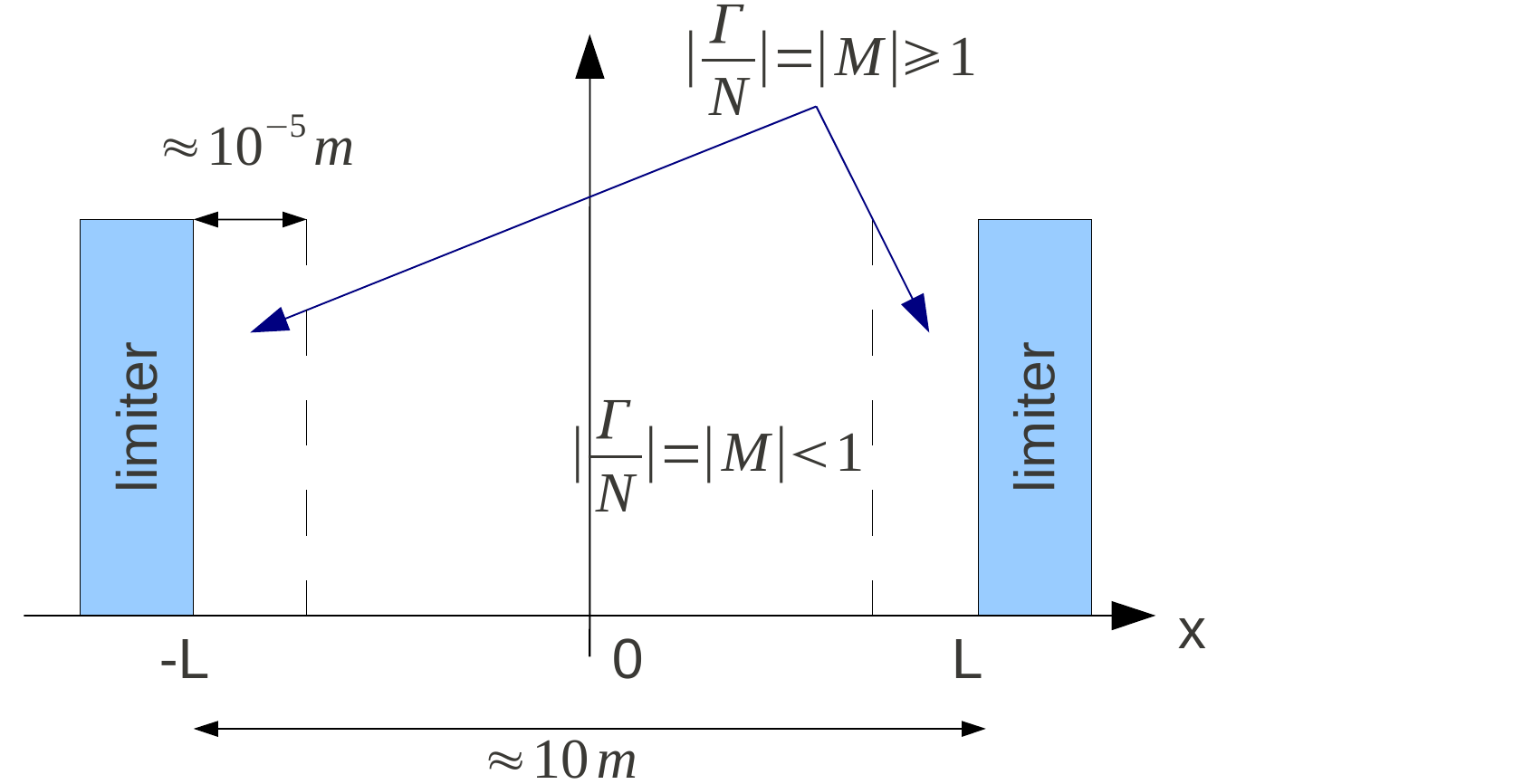}
\end{center}
\caption{
Schematic representation of the scrape-off layer close to the wall of the tokamak. The $x$-axis corresponds to the curvilinear coordinate along a magnetic line.
}
\label{Complete_domain}
\end{figure}

Here, the boundaries of the domain  $x=L$ and $x=-L$
correspond to the limiter ones, which are material obstacles for the fluid (see Fig. \ref{Complete_domain}). In the right-hand side,
$S_N$ and $S_{\Gamma}$ are given source terms. This hyperbolic system is similar to the 1-D isoentropic Euler equation with a linear pressure law
. For sufficiently regular solutions, it can be written in the following non-conservative quasilinear form:
\begin{align}
  &\partial_t \left( \begin{array}{c} N \\ \Gamma \end{array} \right) + 
  \left( \begin{array}{cc} 0 & 1\\ 1 - \dfrac{\Gamma^2}{N^2} & 2 \dfrac{\Gamma}{N} \end{array} \right)
  \partial_x \left(\begin{array}{c} N \\ \Gamma \end{array} \right)=
  \left(\begin{array}{c} S_N \\ S_{\Gamma} \end{array} \right) \qquad (t,x) \in \mathbb{R}^+_* \times ]-L,L[.
\end{align}
We note in the sequel $M = M(\Gamma, N):= \dfrac{\Gamma}{N}$ the Mach number. The eigenvalues of the matrix 
\begin{equation}\label{Matrix_flux}
 \left( \begin{array}{cc} 0 & 1\\ 1 - {M^2} & 2 M \end{array} \right)
\end{equation}
are $\Lambda_1=M-1$ and $\Lambda_2=M+1$, hence $\Lambda_1 <  \Lambda_2$ and the system is strictly hyperbolic.


\noindent{\bf The boundary conditions.}
there is a difficulty with the choice of the boundary conditions for the system (\ref{O1}) that we describe now.
From  physical arguments, it follows that the domain (namely the scrape-off layer) is basically divided
into two regions  \cite{Tam07}:
\begin{itemize}
 \item One region far from the limiter, the pre-sheath, where the plasma is neutral and the Mach number $M=\Gamma / N$ of the plasma satisfies $|M|\leq 1$.
 \item One region next to the limiter (in a thin layer called the sheath area, whose typical thickness is of the order of $10^{-5} m$), where the electroneutrality hypothesis does not hold and we have $|M|>1$. More precisely $M>1$ close to $x=L$ and
 $M<-1$ close to the boundary $x=-L$.
\end{itemize}
At first glance, it could seem natural to prescribe $M=1$ (resp. $M=-1$)
as a boundary condition at $x=L$ (resp. $x=-L$) for the system,
since the physical arguments imply that $M=\pm 1$ very close to the obstacle (Bohm criterion).  
These are exactly the boundary conditions which are chosen in  \cite{Iso10}.
However, in that case, since the eigenvalues are $\Lambda_1=M-1$ and $\Lambda_2=M+1$, it follows that,
at the plasma-limiter interface, one eigenvalue is $0$ (the boundary is characteristic) and the other one corresponds to an outgoing wave (it is also true at $x=-L$). Thus, the problem (\ref{O1}) does not satisfy the usual sufficient conditions for well-posedness, see \cite{Ben07, Gue90, Rau85}: the number of boundary conditions ($=1$) is not equal
to the number of incoming eigenvalues ($=0$). 

In order to test our penalty approach with a well-defined hyperbolic boundary value problem,
in Sections \ref{First penalty approaches} and \ref{New_penalty_method}, we slightly modify the boundary conditions
of the paper \cite{Iso10}, and impose $M= 1 -\eta$ on $x=L$ and $M= -1 + \eta$ on $x=-L$
with a fixed $\eta >0$, which leads to a well-posed hyperbolic problem.

\section{The numerical conservative scheme}\label{Numerical Schemes}
To study numerically the penalty method we need to construct approximate solutions with a sufficient accuracy to obtain a discretization error significantly smaller than the error due to the penalization. As finite volume methods are natural and
efficient for conservation laws, we present a second order finite volume scheme.
To solve the previous nonlinear hyperbolic problem, we use finite volume methods. We tested several schemes \cite{Lev02,God96}, such as Rusanov (modified Lax-Friedrichs), Roe, VFRoe ncv \cite{Buf00, Buf98, Gal03} (possibly with high-order extensions).  

\subsection{General considerations and notations}
We consider an hyperbolic system in the general conservative form in 1D:
\begin{align*}
 &\partial_t \b u + \partial_x \b f (\b u)=\b 0 \quad \textnormal{ in } \mathbb{R}^+_* \times \mathbb{R}\\
 &\b u(0,x)=\b u_0(x) \quad x \in \mathbb{R}.
\end{align*}
We provide a semi-discrete formulation of the problem:
with a uniform spatial mesh of step $\delta x>0$. For each $i \in \mathbb{Z}$, $x_i=i\, \delta x$ is the coordinate of the center of the cell $i$.
For each $t \geq 0$ and for each $i \in \mathbb{Z}$, $\b u_i(t)$ approximates the average of $\b u(t,.)$ over the cell $i$:
\begin{equation*}
 \b u_i(t) = \dfrac{1}{\delta x} \int_{\left( i-\frac12 \right) \delta x}^{\left( i + \frac12 \right)\delta x}{\b u(t,x) dx}.
\end{equation*} 
In this case, the semi-discrete finite volume scheme (of the first-order) has the form:
\begin{align}
&\d_t \b u_i + \dfrac{1}{\delta x} \left( \b F(\b u_i,\b u_{i+1} ) -  \b F (\b u_{i-1},\b u_i ) \right)= \b 0\\
&\b u_i(0)=\dfrac{1}{\delta x} \int_{\left( i-\frac12 \right) \delta x}^{\left( i+\frac12 \right) \delta x}{\b u(0,x) dx}.
\end{align}
The numerical flux $\b F$ is given by the choice of the finite volume scheme which must be consistent, \emph{i.e.}, for all $\b v$, $\b F(\b v, \b v)=\b f(\b v)$.


The VF Roe ncv is described in \cite{Buf00}. We consider  that, in each cell, the value of the unknown function is constant. The VF Roe ncv scheme is an approximate Godunov scheme. The linearized Riemann problem used for the evaluation of the numerical flux, is written with non conservative variables.
In order to have an entropy preserving scheme, a Rusanov entropy correction is implemented.

The VF Roe ncv numerical scheme is of first-order accuracy for smooth solutions. For our applications, the accuracy may not be sufficient, that is why we proposed a second-order extension based on the MUSCL method (Monotone Upwind Scheme for Conservation Laws, see \cite{Van79}) with minmod slope limiter. This method allows us to extend some first-order finite volume schemes up to the second-order. The minmod function is defined as:

\begin{equation*}
\forall (a,b) \in \mathbb{R}^2, \minmod(a,b)=\dfrac{1}{2}(\sign(a)+\sign(b)) \min(|a|,|b|).
\end{equation*}
When $a$ and $b$ are vectors, we use the previous definition component by component. Let us now define the slope terms used for the linear reconstruction of $\b u$:

$\forall t \in \mathbb{R}^*_+, (\b u_i)_x(t)= \minmod\left(\dfrac{\b u_{i+1}(t)-\b u_i(t)}{\delta x},\dfrac{\b u_i(t)-\b u_{i-1}(t)}{\delta x}\right)$.

At the left side of the interface between the cells $i$ and $i+1$, the MUSCL reconstruction of $\b u$ is
\begin{equation*}
\b u_{i+\frac12,l}(t)=\b u_{i}(t) + \dfrac{\delta x}{2}(\b u_i)_x(t)
\end{equation*}
and, at the right side of this interface
\begin{equation*}
\b u_{i+\frac12,r}(t)=\b u_{i+1}(t) - \dfrac{\delta x}{2}(\b u_{i+1})_x(t).
\end{equation*}
Finally, the spatial discretization using the MUSCL reconstruction reads:
\begin{equation*}
\d_t \b u_i + \dfrac{1}{\delta x} \left( \b F\left(\b u_{i+\frac12,l}, \b u_{i+\frac12,r} \right) -  \b F \left(\b u_{i-\frac12,l}, \b u_{i-\frac12,r}\right) \right).
\end{equation*}
Obviously, to have the benefits of a second-order spatial discretization, we need to use a high-order time scheme, such as Heun scheme, to solve the problem.

 \subsection{Application to our problem}\label{Appl_scheme}
 
In our problem, we consider the non-conservative variables $N$ and $M$. 
 $N_i^n$ and $\Gamma_i^n$ approximate the mean values of $N$ and $\Gamma$ over the cell $i$, at the time $t_n$.

The reconstructions of $N$ are defined by
\begin{equation}\label{N_i_n_lr} 
 N_{i+\frac12,l}^n=N_i^n+\dfrac{\delta x}{2}(N_i^n)_x \quad \text{ and } \quad
 N_{i+\frac12,r}^n=N_{i+1}^n-\dfrac{\delta x}{2}(N_{i+1}^n)_x,
\end{equation}
where $(N_i^n)_x= \minmod\left(\dfrac{N_{i+1}^n- N_i^n}{\delta x},\dfrac{N_i^n-N_{i-1}^n}{\delta x}\right)$.
$\Gamma_{i,l}^n, \Gamma_{i,r}^n$ are defined in the same way and $M_{i,l}^n, M_{i,r}^n$ read
\begin{equation}\label{G_i_n_lr} 
 M_{i+\frac12,l}^n=\dfrac{\Gamma_{i+\frac12,l}^n}{N_{i+\frac12,l}^n} \quad \text{ and } \quad M_{i+\frac12,r}^n=\dfrac{\Gamma_{i+\frac12,r}^n}{N_{i+\frac12,r}^n}.
\end{equation}

  The numerical fluxes $f_{N,i+\frac12}^n$ and $f_{\Gamma,i+\frac12}^n$ are evaluated thanks to the following expressions:
\begin{itemize}
 \item Where the entropy correction is not needed:
\begin{equation}\label{flux_without_entrop}
\begin{array}{ll}
\left(
\begin{array}{c}
f_{N,i+\frac12}^n\\
f_{\Gamma,i+\frac12}^n
\end{array}
\right) & = \b F\left(\left(\begin{array}{c} N_{i+\frac12,l}^n \\ \Gamma_{i+\frac12,l}^n \end{array} \right),\left(\begin{array}{c} N_{i+\frac12,r}^n\\ \Gamma_{i+\frac12,r}^n\end{array} \right) \right) \\ & =
\left(\begin{array}{c} 
\widetilde{\Gamma}_{i+\frac12}^{n}(t_n^+,(i+\frac12) \delta x)\\ 
\dfrac{\left(\widetilde{\Gamma}_{i+\frac12}^{n}(t_n^+,(i+\frac12) \delta x)\right)^2}{\widetilde{N}_{i+\frac12}^{n}(t_n^+,(i+\frac12) \delta x)}+\widetilde{N}_{i+\frac12}^{n}(t_n^+,(i+\frac12) \delta x)
\end{array}\right),
\end{array}
\end{equation}
where $\widetilde{N}_{i+\frac12}^{n}(t_n^+,(i+\frac12) \delta x)$ and $\widetilde{\Gamma}_{i+\frac12}^{n}(t_n^+,(i+\frac12) \delta x)=\widetilde{M}_{i+\frac12}^{n}(t_n^+,(i+\frac12) \delta x) \widetilde{N}_{i+\frac12}^{n}(t_n^+,(i+\frac12) \delta x)$ are computed solving the linear Riemann problem written below ($t_n^+$ is any value strictly greater than $t_n$):
\begin{equation}\label{flux_VFRncv_N_G}
\left\{
\begin{array}{l}
\partial_t \left( \begin{array}{c}
\widetilde{N}_{i+\frac12}^{n}\\
\widetilde{M}_{i+\frac12}^{n}
\end{array} \right) +
\left( \begin{array}{cc}
\frac12\left(M_{i+\frac12,l}^n+M_{i+\frac12,r}^n\right) & \frac12\left(N_{i+\frac12,l}^n+N_{i+\frac12,r}^n\right) \\
\dfrac{2}{N_{i+\frac12,l}^n+N_{+\frac12i,r}^n} & \frac12\left(M_{i+\frac12,l}^n+M_{i+\frac12,r}^n\right)
\end{array} \right)
\partial_x \left( \begin{array}{c}
\widetilde{N}_{i+\frac12}^{n}\\
\widetilde{M}_{i+\frac12}^{n}
\end{array} \right)
=
\left( \begin{array}{c}
0\\
0
\end{array} \right)\\
\left( \begin{array}{c}
\widetilde{N}_{i+\frac12}^{n}(t_n,x)\\
\widetilde{M}_{i+\frac12}^{n}(t_n,x)
\end{array} \right) = \left( \begin{array}{c}
N_{i+\frac12,l}^{n}\\
M_{i+\frac12,l}^{n}
\end{array} \right) \text{ if } x < (i+\frac12) \delta x \qquad \text{ and } \qquad \left( \begin{array}{c}
N_{i+\frac12,r}^{n}\\
M_{i+\frac12,r}^{n}
\end{array} \right) \text{ if } x \geq (i+\frac12) \delta x.
\end{array}
\right.
\end{equation}
 \item Where the entropy condition is needed, \emph{i.e.} if $M_{i+\frac12,l}^n - 1 \leq 0 \leq M_{i+\frac12,r}^n - 1$ (and $M_{i+\frac12,l}^n \neq M_{i+\frac12,r}^n$) or if $M_{i+\frac12,l}^n + 1\leq 0 \leq M_{i+\frac12,r}^n + 1$ (and $M_{i+\frac12,l}^n \neq M_{i+\frac12,r}^n$), the flux is replaced by a Rusanov flux:
\begin{equation}\label{flux_Rusanov}
\begin{array}{ll}
\left(
\begin{array}{c}
f_{N,i+\frac12}^n\\
f_{\Gamma,i+\frac12}^n
\end{array}
\right) & = \b F\left(\left(\begin{array}{c} N_{i+\frac12,l}^n \\ \Gamma_{i+\frac12,l}^n \end{array} \right),\left(\begin{array}{c} N_{i+\frac12,r}^n\\ \Gamma_{i+\frac12,r}^n\end{array} \right) \right) \\
 & =\dfrac12 \left(\begin{array}{c}
 \Gamma_{i+\frac12,l}^n\!+\!\Gamma_{i+\frac12,r}^n\\ \dfrac{\left( \Gamma_{i+\frac12,l}^n\right)^2}{N_{i+\frac12,l}^n} \! + \! \dfrac{\left( \Gamma_{i+\frac12,r}^n\right)^2}{N_{i+\frac12,r}^n} \! + \! N_{i+\frac12,l}^n \! + \! 
 N_{i+\frac12,r}^n \end{array} \right)\\ & \qquad +\dfrac12\left(\!\max \left\{\!|M_{i+\frac12,l}^n|,|M_{i+\frac12,r}^n|\!\right\}\!+\!1  \!\right)\left(\begin{array}{c} N_{i+\frac12,r}^n \! - \! N_{i+\frac12,l}^n \\ \Gamma_{i+\frac12,r}^n \!-\! \Gamma_{i+\frac12,l}^n \end{array} \right),
 \end{array}
 \end{equation}
 remembering that the spectral radius of the matrix given in the expression (\ref{Matrix_flux}) is $\max\{|\Lambda_1|,|\Lambda_2|\}=|M|+1$.
\end{itemize}

Finally, the full discretization reads:
\begin{align*}
 & N_i^{1,n}= N_i^{n}-\dfrac{\delta t}{\delta x}\left(f_{N,i+\frac12}^n-f_{N,i-\frac12}^n \right)+\delta t\, S_{N,i}^n \\
 & \Gamma_i^{1,n}= \Gamma_i^{n}-\dfrac{\delta t}{\delta x}\left(f_{\Gamma,i+\frac12}^n-f_{\Gamma,i-\frac12}^n \right) +\delta t\, S_{\Gamma,i}^n \\
 & N_i^{n+1}= \dfrac12 (N_i^{1,n}+N_i^{n})-\dfrac{\delta t}{2 \delta x}\left(f_{N,i+\frac12}^{1,n}-f_{N,i-\frac12}^{1,n} + f_{N,i+\frac12}^{n}-f_{N,i-\frac12}^{n} \right)+\dfrac{\delta t}{2}\, (S_{N,i}^n +S_{N,i}^{n+1} ) \\
 & \Gamma_i^{n+1}= \dfrac12 (\Gamma_i^{1,n}+\Gamma_i^{n})-\dfrac{\delta t}{2 \delta x}\left(f_{\Gamma,i+\frac12}^{1,n}-f_{\Gamma,i-\frac12}^{1,n}+ f_{\Gamma,i+\frac12}^n-f_{\Gamma,i-\frac12}^n \right) +\dfrac{\delta t}{2}\, (S_{\Gamma,i}^n + S_{\Gamma,i}^{n+1}),
\end{align*}
where the upper index $1,n$ corresponds to the intermediate step of the Heun scheme. $f_{N,i+\frac12}^{1,n},f_{\Gamma,i+\frac12}^{1,n}$ are evaluated using the formulas (\ref{N_i_n_lr})-(\ref{flux_Rusanov}), replacing the terms $(N_i^n)_{i \in \mathbb{Z}},(\Gamma_i^n)_{i \in \mathbb{Z}}$ by $(N_i^{1,n})_{i \in \mathbb{Z}},(\Gamma_i^{1,n})_{i \in \mathbb{Z}}$.
This Section only concerns the Initial Value Problem, and we did not consider boundary conditions nor implementation of penalization methods. We used an adaptive time step based on a CFL-like condition: for all $n$, the time step $\delta t$ satisfies $\max_{i}\{|M_i^n|+1\} \dfrac{\delta t}{\delta x}=0.8$. Besides, as we shall see later, the penalization, which adds discontinuous terms of large amplitude compared to the others, destabilizes the numerical scheme, if these terms are treated explicitly. So we shall need to modify this scheme 
to deal with these terms implicitly.

All the numerical tests presented in the following Sections of this paper use the finite volume scheme presented in this subsection (see formulas (\ref{N_i_n_lr})-(\ref{flux_Rusanov})).

\section{First penalty approaches}\label{First penalty approaches}
\subsection{A first penalty method}\label{Ssect_first_penalty}

The following penalty approach has been proposed by Isoardi 
\emph{et al.} \cite{Iso10} for the problem (\ref{O1}) with Bohm criterion as boundary conditions, which is recalled below:
\begin{equation}\label{Hyp_pb_Bohm}
\left\{
 \begin{array}{l}
     \partial_t N + \partial_x \Gamma=S_N \\
       \partial_t \Gamma + \partial_x \left(\dfrac{\Gamma^2}{N} + N\right) = S_{\Gamma} \\
       M(.,-L)=-1 \textsf{ and } M(.,L)=1\\
       N(0,.)=N_0 \textsf{ and } \Gamma(0,.)=\Gamma_0.
    \end{array}
    \right .
    \quad \qquad (t,x) \in \mathbb{R}^+_* \times ]-L,L[
\end{equation}
 Let $\chi$ be the characteristic function of the limiter 
, \emph{i.e.} $\chi(x)=1$ if $x$ is inside the limiter, and  $\chi(x)=0$ elsewhere, and
$\varepsilon>0$ the penalization parameter. The penalized system is given by:
\begin{equation}\label{Penal_pb}
\left\{\begin{array}{l} 
 \partial_t N + \partial_x \Gamma + \dfrac{\chi}{\varepsilon} N = (1-\chi) S_N \qquad \textnormal{ in } \mathbb{R}^+_* \times \mathbb{R}\\
 \partial_t \Gamma + (1-\chi) \partial_x \left(\dfrac{\Gamma^2}{N} + N \right) +\dfrac{\chi}{\varepsilon}(\Gamma-M_0 N)=(1-\chi) S_{\Gamma}\\
  N(0,.)=N_0 \textsf{ and } \Gamma(0,.)=\Gamma_0,
\end{array}
\right.
\end{equation}
where $M_0$ is a function such that, at the plasma-limiter interface we have $|M_0|=1$.
Here, the two components of the unknown are penalized although there is no incoming wave. At least formally, 
$N$ is enforced to converge to $0$ inside the limiter when $\varepsilon$ tends to $0$, whereas $\Gamma$ is enforced to $M_0 N$ inside the limiter satisfying the Bohm criterion.

The flux of the second equation is cut inside of the limiter, and this causes 
some troubles from the mathematical point of view. Indeed, the system (\ref{Penal_pb}) is
an hyperbolic system with discontinuous coefficients and the meaning of the term
$$
(1-\chi) \partial_x \left(\dfrac{\Gamma^2}{N} + N \right)
$$
is not clear because it can involve the product of a measure with a discontinuous function which has
no distributional sense. As a confirmation of this fact,
our numerical tests show the existence of a strong singularity at the interface for the numerical discrete solution. 
Concerning the interpretation of this numerical singularity, it could happen
(but we don't have any rigorous proof and this is just an open question) that this system admits 
generalized solutions in the spirit of Bouchut-James  \cite{Bou98} (see also Poupaud-Rascle \cite{Pou97}, or  Fornet-Gu\`es \cite{For08}) such as measure-valued solutions, which can for example exhibit a Dirac measure at the interface,
and this generalized solution could be selected by the numerical approximation process.  

\begin{figure}
\begin{center}
\includegraphics[scale=0.50, trim = 10mm 203mm 90mm 10mm, clip=true]{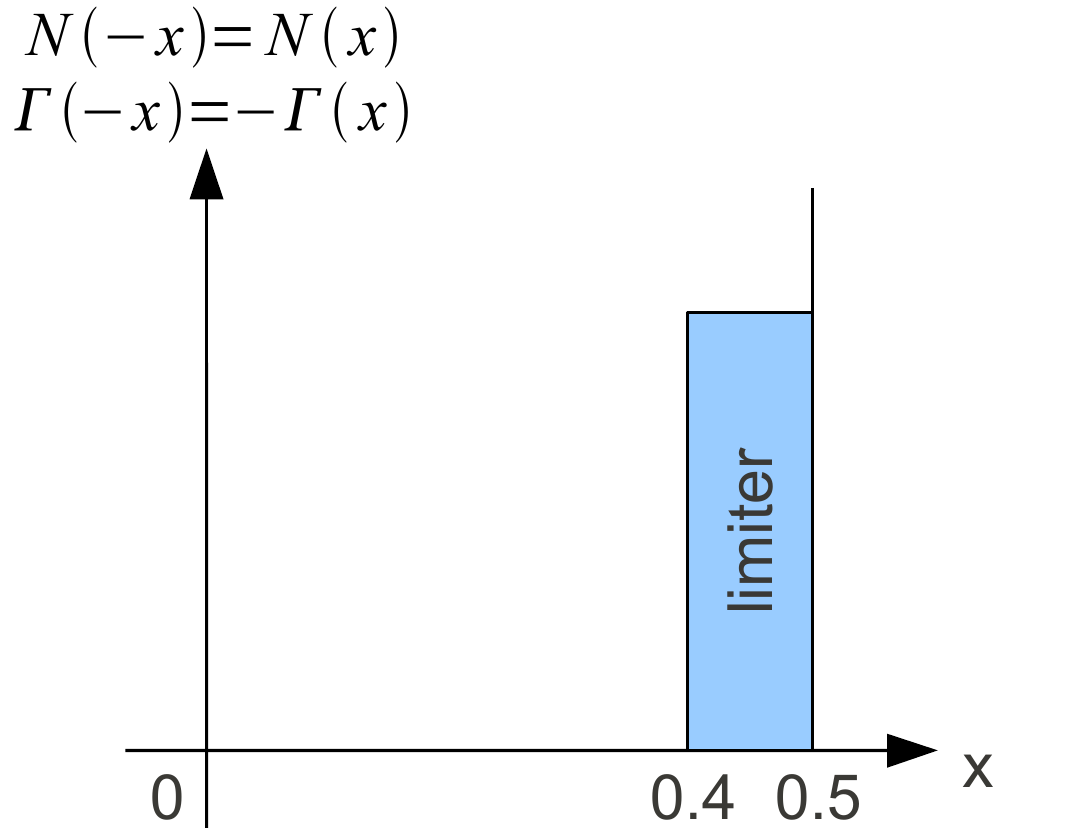}
\end{center}
\caption{
Representation of the computational domain, $x \in [0,0.5]$. The plasma area corresponds to $x \in [0,L]$, with $L=0.4$.
}
\label{Half_domain1}
\end{figure}

We choose $S_N$ and $S_{\Gamma}$ so that the following functions define a solution of the boundary value problem (\ref{O1}):
\begin{equation}\label{Sol_ex_pb_ini}
N(t,x)=\exp \left(\dfrac{-x^2}{0.16 (t+1)}\right) \qquad \Gamma(t,x)=\sin \left(\dfrac{\pi x}{0.8} \right) \exp \left(\dfrac{-x^2}{0.16 (t+1)}\right).
\end{equation}
This test solution is regular (at least inside the plasma area) and has no singularity at the plasma-limiter interface. The computational domain is $x \in [0,0.5]$ where the limiter set is $x \in [0.4,0.5]$ (see Fig. \ref{Half_domain1}).

\begin{figure}
\begin{center}
\includegraphics[scale=0.65, trim = 10mm 5mm 10mm 10mm, clip=true]{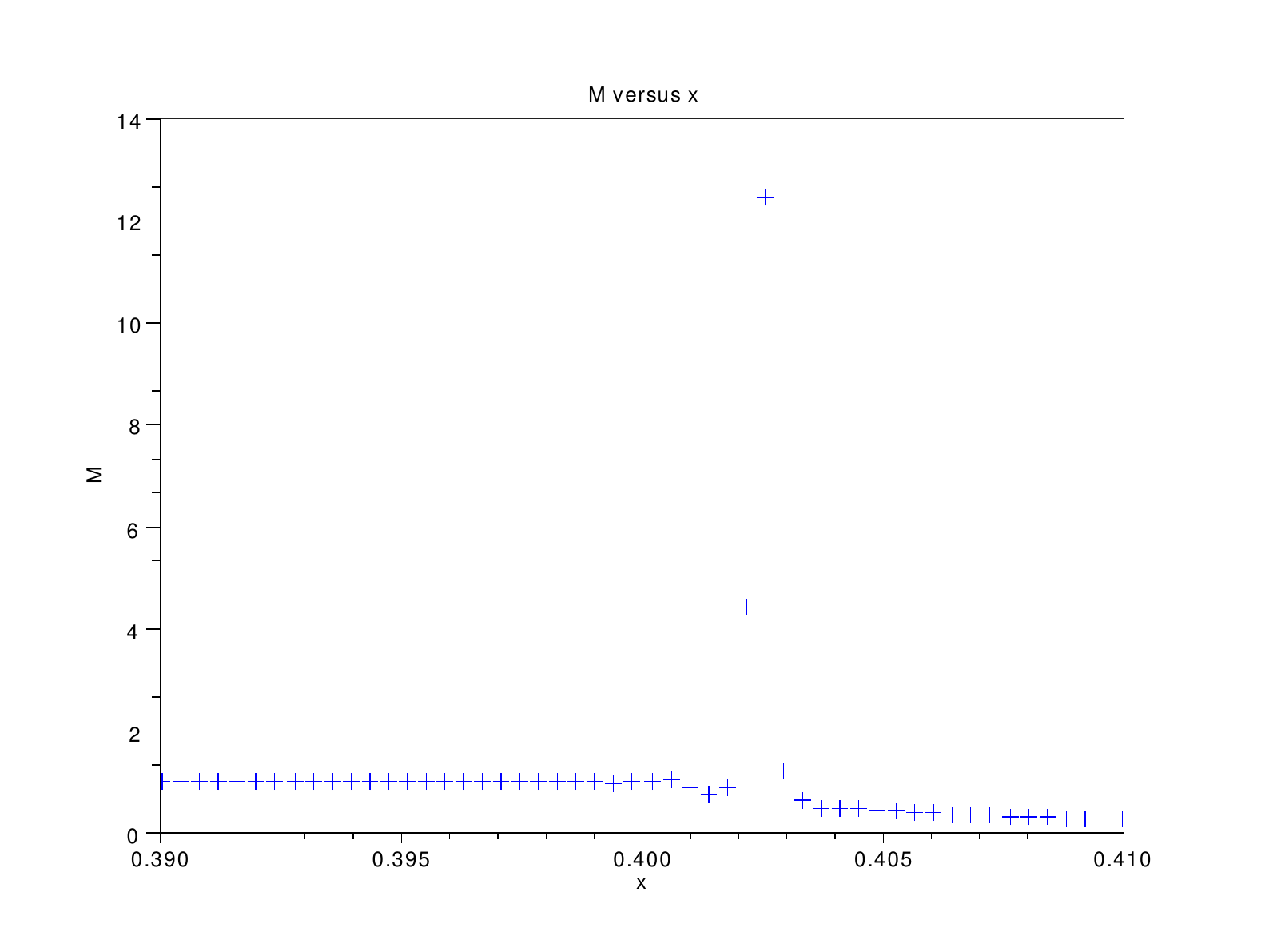}
\includegraphics[scale=0.65, trim = 10mm 5mm 10mm 10mm, clip=true]{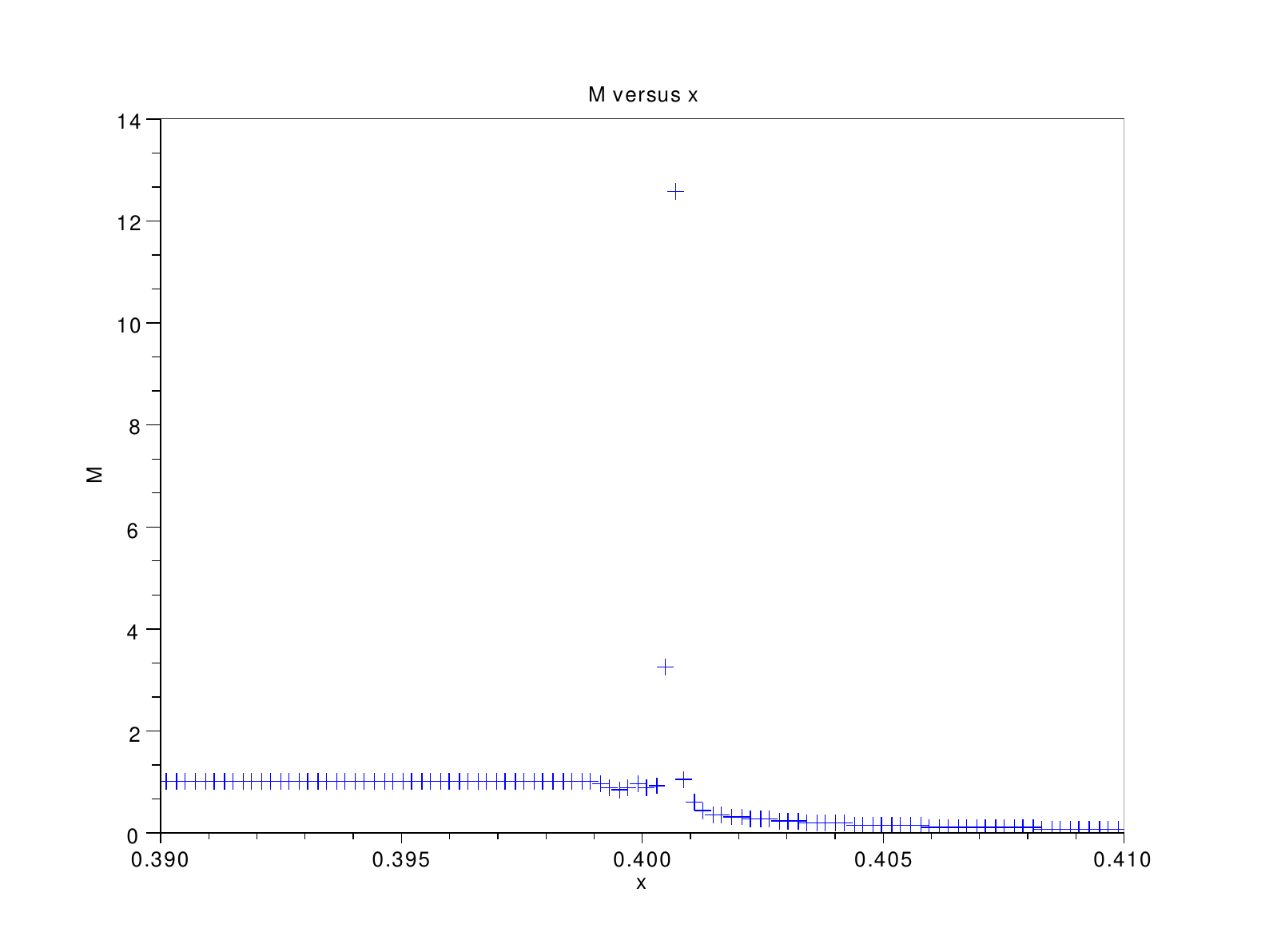}
\includegraphics[scale=0.65, trim = 10mm 5mm 10mm 10mm, clip=true]{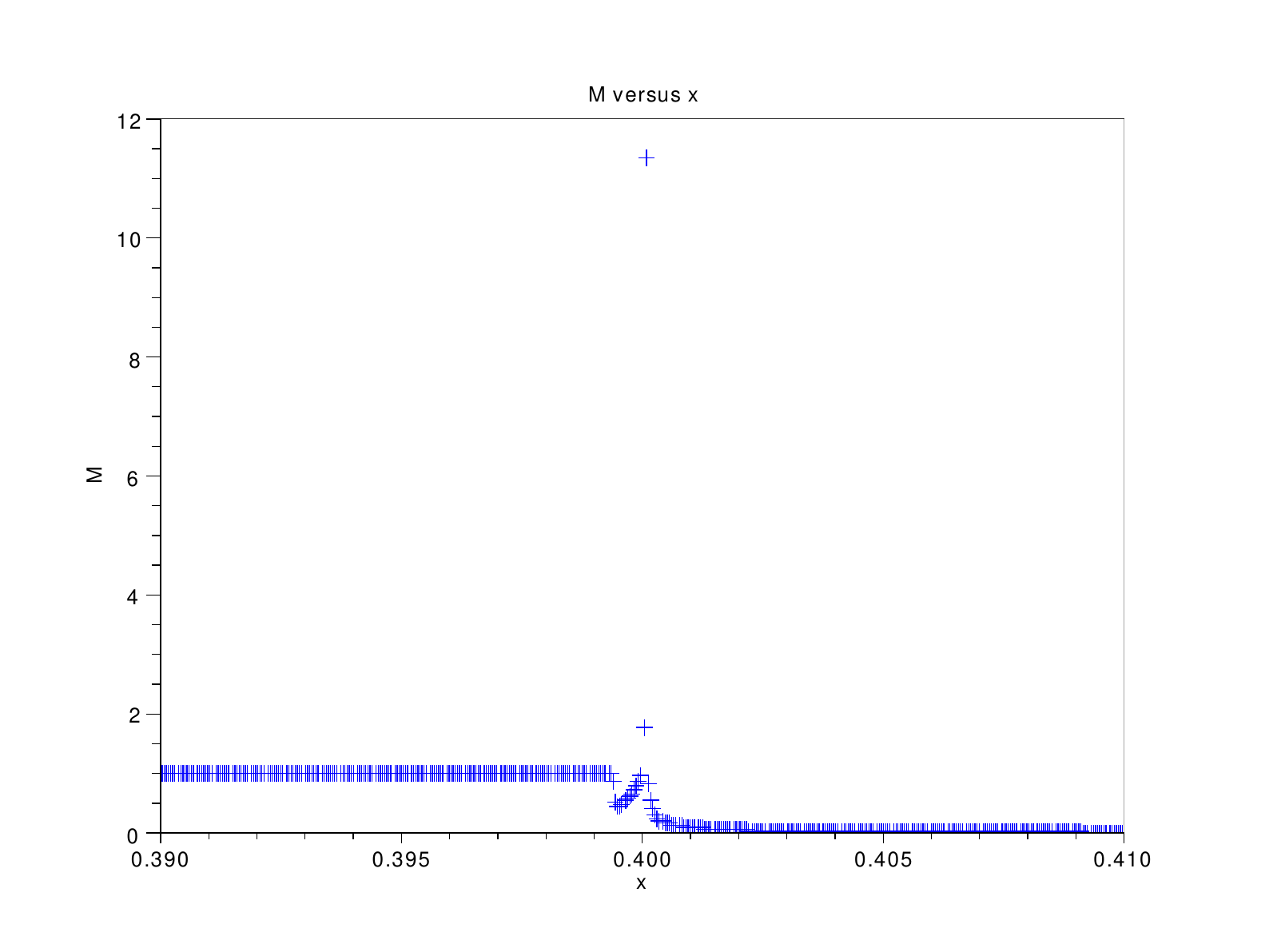}
\end{center}
\caption{
$M$ versus $x$ with $\varepsilon=10^{-3}$, for three different meshes (respectively $1280$, $2560$ and $10240$ cells), using the penalization of Isoardi \emph{et al.} \cite{Iso10}. The computations are stopped when $\max_{i \in \{1,\dots,J\}}(|M_i^n|)>10$, which corresponds to the following times: $t= 0.008822$, $t=0.004107$ and $t=0.0015834$. The computational domain was $[0,0.5]$ and $L=0.4$ (plasma-limiter interface). At $x=0$, we impose a symmetry boundary condition (see Fig. \ref{Half_domain1}).
}
\label{Dirac_measure}
\end{figure}

We perform a mesh convergence study with a fixed value for $\varepsilon=10^{-3}$, using the numerical scheme presented in Section \ref{Numerical Schemes}. In Fig. \ref{Dirac_measure}, we observe that a peak appears very quickly near the plasma-limiter interface. Then, $|M_i^n|$ becomes very large (about $10^8$) within a few points.
When the resolution increases, the peak is nearer and nearer to the plasma-limiter interface and appears earlier and earlier. We stop the computations when $\max_{i \in \{1,\dots,J\}}(|M_i^n|)>10$ but similar results have been obtained when the stop criterion is $\max_{i \in \{1,\dots,J\}}(|M_i^n|)>100$.
This leads us to believe that, if the solution converges to a generalized solution of the continuous problem, then 
this generalized solution must have a singularity supported by the interface (that could be a Dirac measure for example). 
We notice that the presence of a Dirac measure at the interface is not only a theoretical issue since it has been observed numerically and that the Dirac measure destabilizes the  numerical scheme. In the following section, we propose a modification of the boundary value problem to obtain a well-posed version.

A second motivation to modify the penalization approach, is the fact that the 
hyperbolic boundary value problem (\ref{Hyp_pb_Bohm})
does not satisfy the conditions to apply usual well-posedness results, as already explained in section 
\ref{Presentation of the hyperbolic system}.

\subsection{Penalization of the two-fields for the modified boundary conditions}\label{Subsec 2 field penal}

Obviously, the plasma density inside the limiter has to be negligible. So, we may expect that the penalization enforces $N$ to $0$ inside the limiter. The issue yields that the hyperbolic system is not valid when $N=0$. First, we tried to penalize only $N$ to $0$ inside the limiter without imposing any condition on $M$ or $\Gamma$, which implies that $M$ tends to infinity. To ensure the CFL stability condition, the time step must tend to $0$. So, if we want to penalize $N$ to $0$ inside the limiter, we also need to impose $M$ or $\Gamma$. 

The Bohm criterion ($|M|=1$ at the boundary) comes from the continuous connection of two different physical regimes between the pre-sheath and the sheath part. However, the hyperbolic system only models the plasma transport in the pre-sheath area where the electroneutrality hypothesis holds. So we consider that, at the boundary of the plasma domain, the Mach number is not exactly equal to $\pm 1$ but to $\pm (1-\eta)$ for $\eta>0$ sufficiently small. In the configuration of Fig. \ref{Complete_domain}, we impose $M(t,L)-1=-\eta<0$ and $M(t,-L)+1=\eta>0$, so that we have one incoming wave at the plasma-limiter interfaces $x=\pm L$. We choose these conditions as boundary conditions for the hyperbolic system, and thus model the plasma transport by the following initial boundary-value problem:
\begin{equation} \label{O2}
\left\{
 \begin{array}{l}
     \partial_t N + \partial_x \Gamma=S_N \\
       \partial_t \Gamma + \partial_x \left(\dfrac{\Gamma^2}{N} + N\right)=S_\Gamma \\
       M(.,-L)=-1+\eta \textsf{ and } M(.,L)=1-\eta   \\
       N(0,.)=N_0 \textsf{ and } \Gamma(0,.)=\Gamma_0.
    \end{array}
    \right.
    \quad \qquad (t,x) \in \mathbb{R}^+_* \times ]-L,L[
\end{equation}
For this problem, the boundary is not characteristic, and the boundary conditions are
maximally dissipative, as it is explained later in the paper (see the Definition \ref{Max_dissip}, at the end of Subsection \ref{asympt_exp}). Hence, the problem has a unique solution which is smooth up to a time $T$ for compatible initial data, see for example \cite{Gue90, Rau74} and Theorem 11.1 of \cite{Ben07}.

A natural penalized system could be:
\begin{equation}\label{Two_fields_penal}
 \left\{\begin{array}{l}
 \partial_t N + \partial_x \Gamma + \dfrac{\chi}{\varepsilon} N = S_N\\
 \partial_t \Gamma + \partial_x \left(\dfrac{\Gamma^2}{N} + N\right) + \dfrac{\chi}{\varepsilon}\left(\dfrac{\Gamma}{M_0} - N\right)=S_{\Gamma}
\end{array}\right. \qquad \textnormal{ in } \mathbb{R}^+_* \times \mathbb{R}.
\end{equation}
 In this situation, two-fields are penalized. The main goal of this subsection is to show how a boundary layer due to the penalization can be put in evidence thanks to a numerical study, as these tools are then used in Subsections \ref{ssect_1_side}, \ref{Subsec_2sides_penal} and \ref{ssect_M_O_to_1}.
 
For the numerical convergence analysis, we try to impose the following test solution:
\begin{equation}\label{Test_sol}
N(t,x)=\exp \left(\dfrac{-x^2}{0.16 (t+1)}\right) \qquad \Gamma(t,x)=M_0 \sin \left(\dfrac{\pi x}{0.8} \right) \exp \left(\dfrac{-x^2}{0.16 (t+1)}\right).
\end{equation}
 $S_N$ and $S_{\Gamma}$ are chosen so that (\ref{Test_sol}) is the solution of (\ref{O2}) in the plasma area and are null inside the limiter set. These source terms do not depend on $\varepsilon$. The formula (\ref{Test_sol}) differs from (\ref{Sol_ex_pb_ini}) because of the factor $M_0$ in $\Gamma(t,x)$. Besides its regularity, this solution has been chosen because it is not stationary. Obviously, the initial conditions are the traces of the imposed test solution at time $t=0$.
 
\begin{figure}
\begin{center}
\subfigure[Thick black: $N(1,x)$, Gray: $\Gamma(1,x)$, Narrow black: $M(1,x)$, $\varepsilon = 0.1$]{\includegraphics[scale=0.32, angle=-90, trim = 20mm 20mm 3mm 20mm, clip=true]{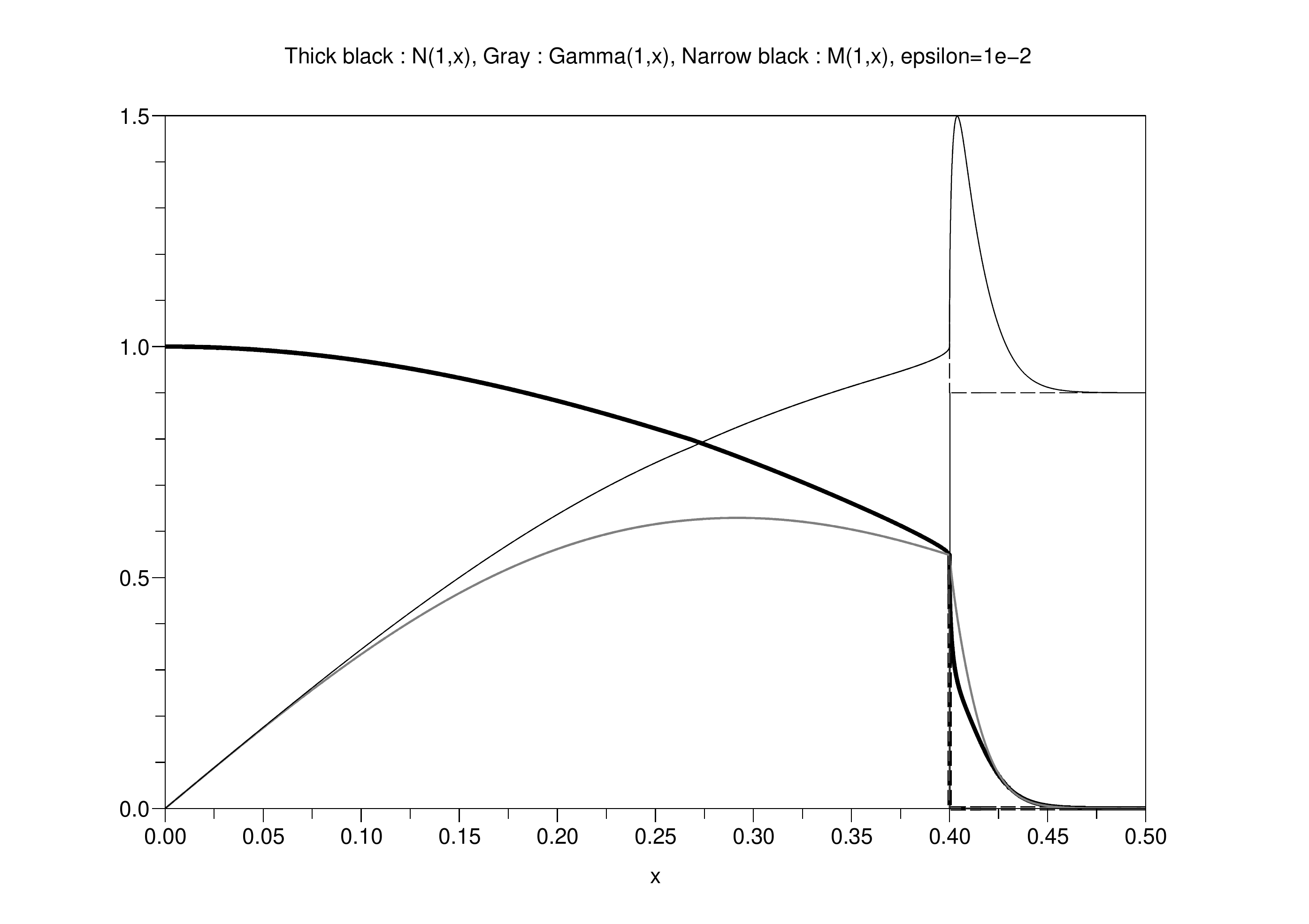}} $\quad$
\subfigure[Thick black: $N(1,x)$, Gray: $\Gamma(1,x)$, Narrow black: $M(1,x)$, $\varepsilon = 10^{-5}$]{\includegraphics[scale=0.32, angle=-90, trim = 20mm 20mm 3mm 20mm, clip=true]{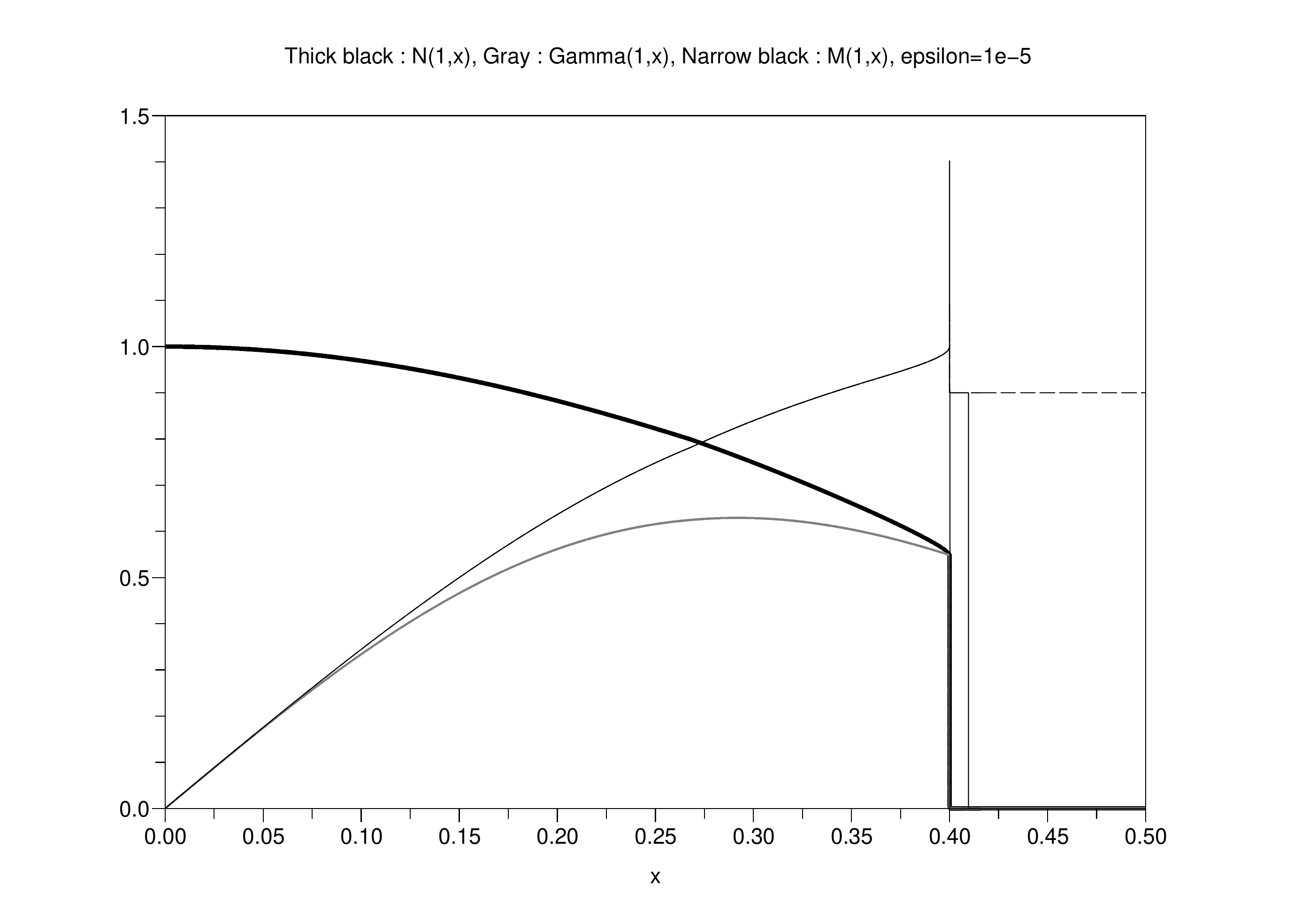}}
\end{center}
\caption{
Plot of $N$, $\Gamma$ and $M$ as functions of $x$ (at $t=1$) with the two-fields penalty method at the left for $\varepsilon=0.01$ and at the right for $\varepsilon=10^{-5}$. The continuous lines represent the numerical solutions whereas the dotted lines corresponds to the reference solution (when $\varepsilon= 10^{-20}$). The limiter corresponds to the area $x\in[0.4,0.5]$. The computational domain is the one presented in Fig. \ref{Half_domain1}. The mesh step is $\delta x=10^{-5}$. The penalization method considered is described in Subsection \ref{Subsec 2 field penal}.
}
\label{Plot_2f_regular}
\end{figure}

\begin{figure}
\begin{center}
\subfigure[$L^1$ error for $N$ in the plasma ($+$), $N$ in the limiter ($\times$), $\partial_x N$ in the plasma ($\circ$) and $\partial_x N$ in the limiter ($*$)]{\includegraphics[scale=0.56, trim = 7.5mm 3mm 19mm 13.3mm, clip=true]{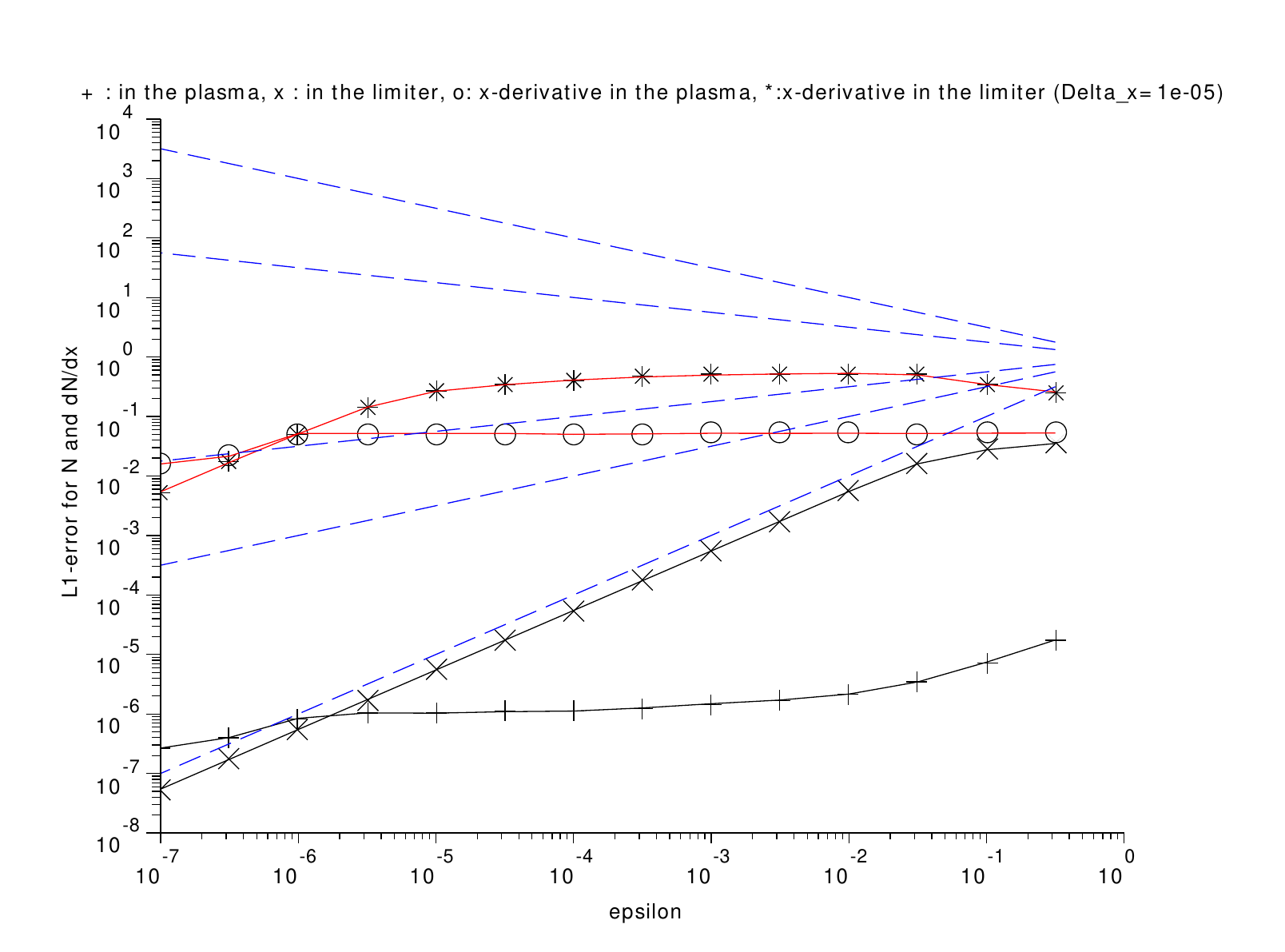}} $\quad$
\subfigure[$L^2$ error for $N$ in the plasma ($+$), $N$ in the limiter ($\times$), $\partial_x N$ in the plasma ($\circ$) and $\partial_x N$ in the limiter ($*$)]{\includegraphics[scale=0.56, trim = 7.5mm 3mm 19mm 13.3mm, clip=true]{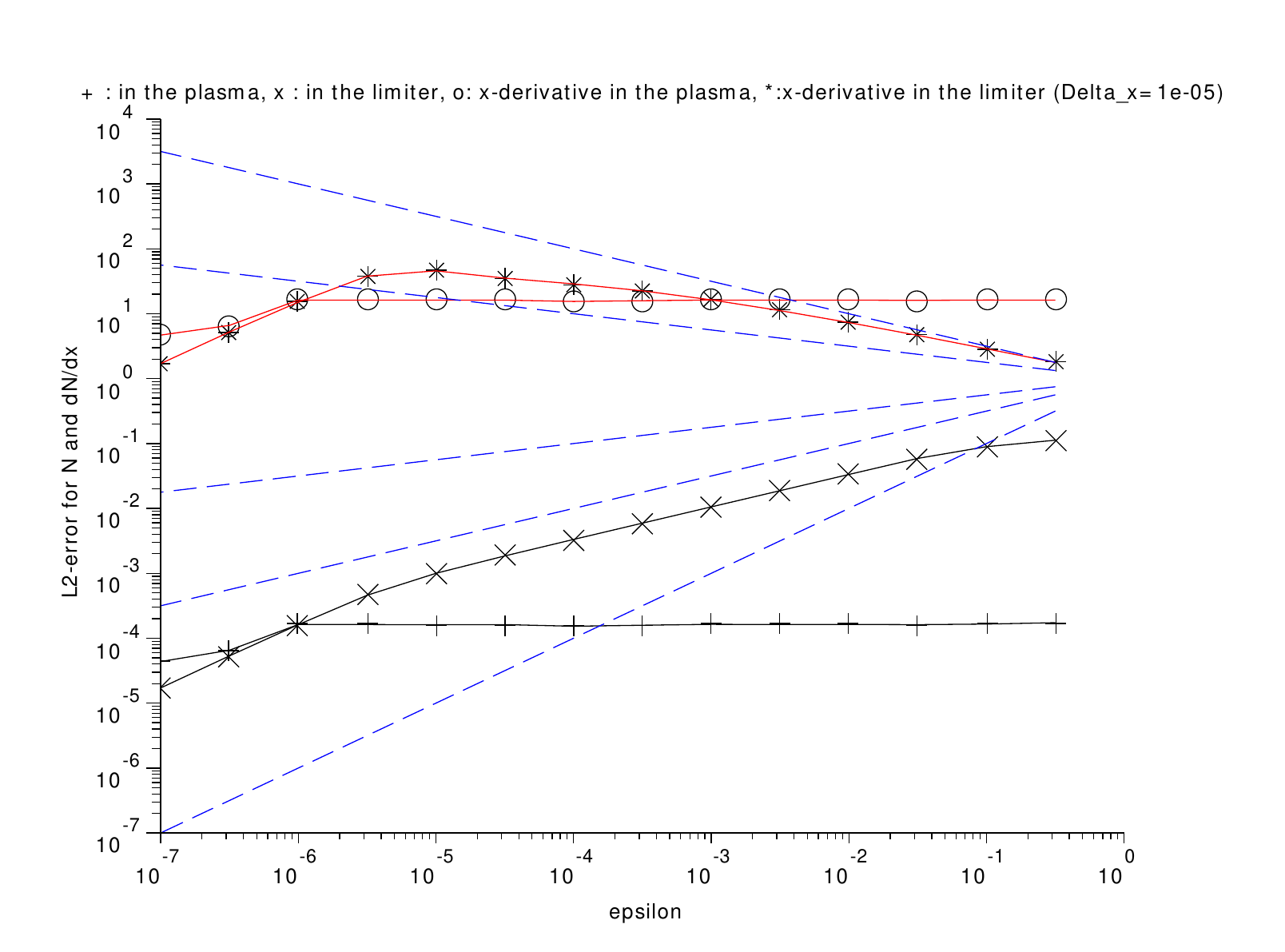}} 
\subfigure[$L^1$ error for $\Gamma$ in the plasma ($+$), $\Gamma$ in the limiter ($\times$), $\partial_x \Gamma$ in the plasma ($\circ$) and $\partial_x \Gamma$ in the limiter ($*$)]{\includegraphics[scale=0.56, trim = 7.5mm 3mm 19mm 13.3mm, clip=true]{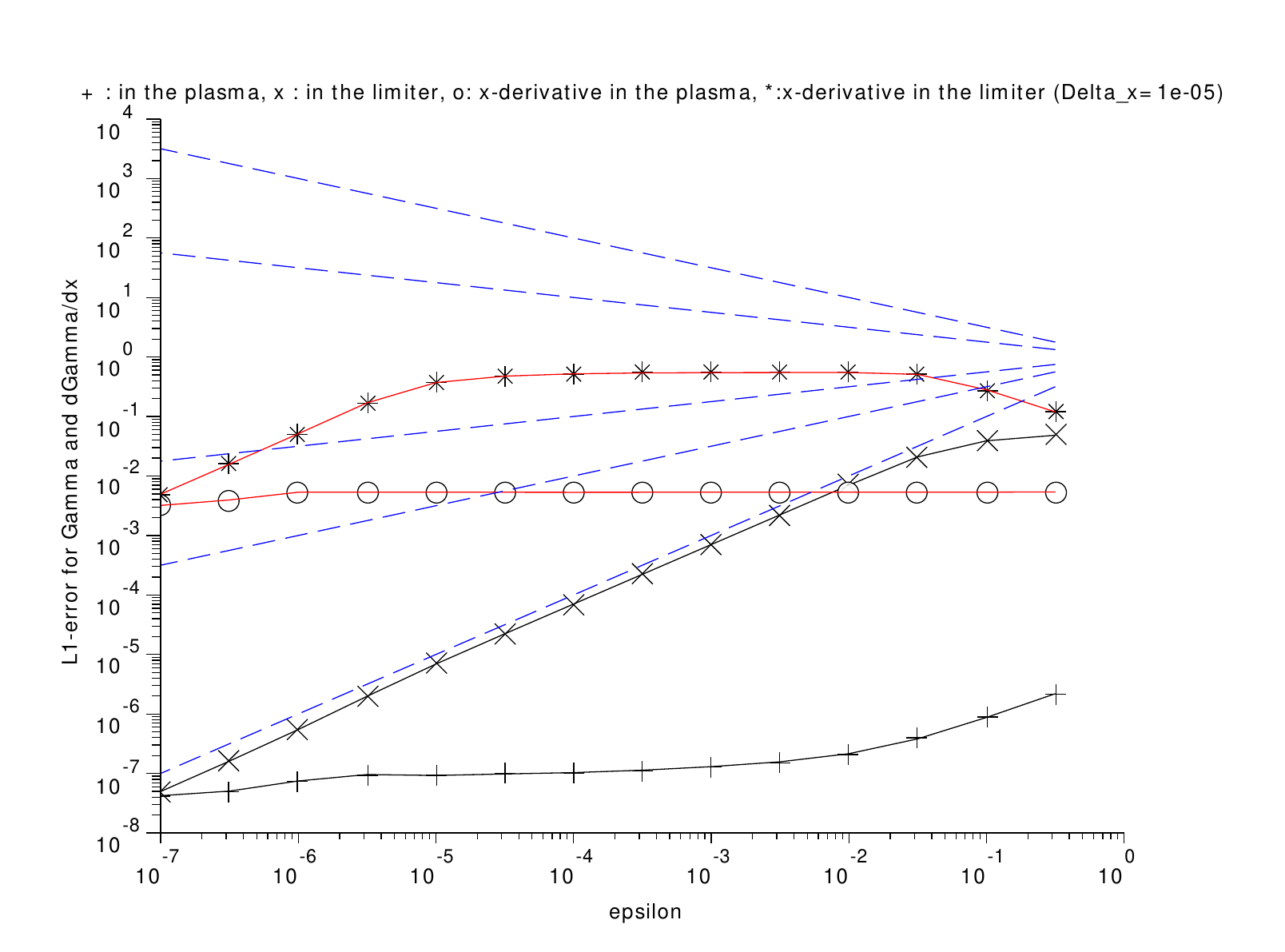}} $\quad$
\subfigure[$L^2$ error for $\Gamma$ in the plasma ($+$), $\Gamma$ in the limiter ($\times$), $\partial_x \Gamma$ in the plasma ($\circ$) and $\partial_x \Gamma$ in the limiter ($*$)]{\includegraphics[scale=0.56, trim = 7.5mm 3mm 19mm 13.3mm, clip=true]{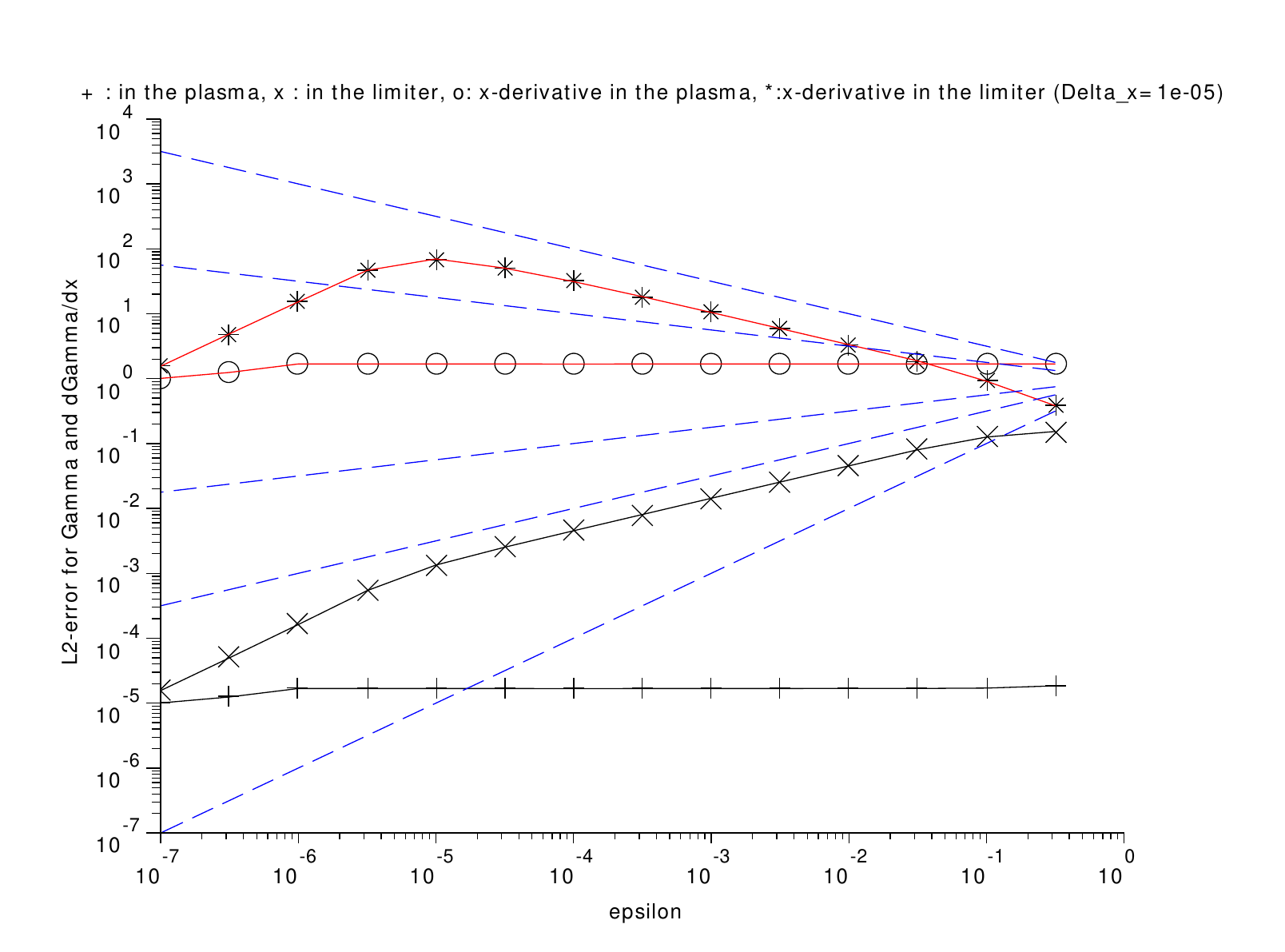}} 
\end{center}
\caption{
Errors for $N$, $\partial_x N$, $\Gamma$ and $\partial_x \Gamma$ in $L^1$ and $L^2$ norms with the two-fields penalty method, see (\ref{Two_fields_penal}). The dashed lines represent respectively the curves $\varepsilon^{-\frac12},\varepsilon^{-\frac14},\varepsilon^{1/4}, \varepsilon^{1/2}$ and $\varepsilon$. The mesh step is $\delta x=10^{-5}$. The penalization method considered is described in Subsection \ref{Subsec 2 field penal}.
}
\label{Error_2f_regular}
\end{figure}

 The numerical tests presented below (see Fig. \ref{Plot_2f_regular} and \ref{Error_2f_regular}) have been performed to show that the two-fields penalization generates a boundary layer which is captured by the numerical scheme when the mesh step is sufficiently small. The numerical approximation using the finite volume scheme described in Subsection \ref{Appl_scheme} is given by:
\begin{align*}
 & N_i^{1,n}= \dfrac{N_i^{n}-\dfrac{\delta t}{\delta x}\left(f_{N,i+\frac12}^n-f_{N,i-\frac12}^n \right)+\delta t\, S_{N,i}^n}{1+ \frac{\delta t}{\varepsilon} \chi} \\
 & \Gamma_i^{1,n}= \dfrac{\Gamma_i^{n}-\dfrac{\delta t}{\delta x}\left(f_{\Gamma,i+\frac12}^n-f_{\Gamma,i-\frac12}^n \right)+ \delta t \dfrac{\chi}{\varepsilon} M_0 +\delta t\, S_{\Gamma,i}^n}{1+\delta t \dfrac{\chi}{\varepsilon N_i^{1,n}}} \\
 & N_i^{n+1}= \dfrac{\dfrac12 (N_i^{1,n}+N_i^{n})-\dfrac{\delta t}{2 \delta x}\left(f_{N,i+\frac12}^{1,n}-f_{N,i-\frac12}^{1,n} + f_{N,i+\frac12}^n - f_{N,i-\frac12}^n \right)+\dfrac{\delta t}{2}\, (S_{N,i}^n +S_{N,i}^{n+1} )}{1+ \frac{\delta t}{\varepsilon} \chi} \\
 & \Gamma_i^{n+1}= \dfrac{\dfrac12 (\Gamma_i^{1,n}+\Gamma_i^{n})-\dfrac{\delta t}{2 \delta x}\left(f_{\Gamma,i+\frac12}^{1,n}-f_{\Gamma,i-\frac12}^{1,n}+ f_{\Gamma,i+\frac12}^n-f_{\Gamma,i-\frac12}^n \right)+ \delta t \dfrac{\chi}{\varepsilon} M_0 +\dfrac{\delta t}{2}\, (S_{\Gamma,i}^n + S_{\Gamma,i}^{n+1})}{1+\delta t \dfrac{\chi}{\varepsilon N_i^{n+1}} },
\end{align*}
 where the numerical fluxes $f_{N,i+\frac12}^n,f_{\Gamma,i-\frac12}^n$ are evaluated using the formulas (\ref{N_i_n_lr})-(\ref{flux_Rusanov}). The upper index $1,n$ corresponds to the intermediate step of the Heun scheme. To improve the stability of the scheme, the penalized terms have been treated implicitly.
 
 The Fig. \ref{Plot_2f_regular} shows that the limit solution, when $\varepsilon$ tends to $0$, is not the regular one imposed in (\ref{Test_sol}) but it appears that $M$ is close to $1$ at the plasma-limiter interface (though $M_0=1-\eta=0.9$). Similar results are observed for $\eta=0.01$. Ghendrih \emph{et al.} in \cite{Ghe11} explain this phenomenon as a consequence of the fact that $N$ is enforced to $0$ in the limiter. Finally, as we don't have the exact solution of (\ref{Two_fields_penal}), for the numerical tests, we consider that the reference solution is the one obtained by our numerical scheme with $\varepsilon= 10^{-20}$.
 
\begin{figure}
\begin{center}
\includegraphics[scale=0.56, trim = 7.5mm 3mm 19mm 13.3mm, clip=true]{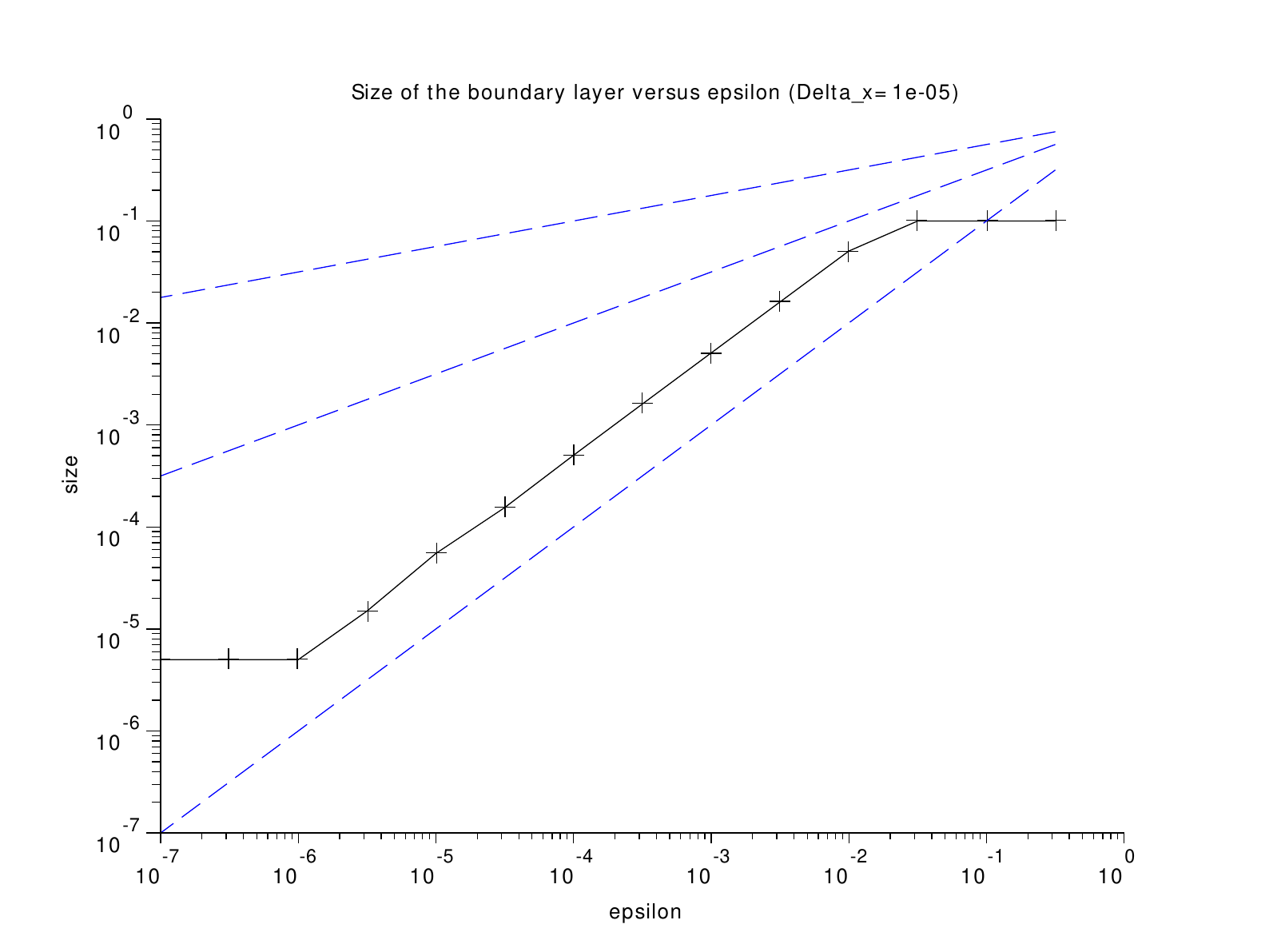}
\end{center}
\caption{
Numerical estimation of the size of the boundary layer versus the penalization parameter $\varepsilon$. The thickness is computed using
$x_{J_{bl}} - 0.4$ where $J_{bl} = \max\{|N_i^n - 0|<0.01 |N_I^n - 0|\}$ with $I$ the number of the cell in the plasma area in contact with the plasma-limiter interface. Indeed, $x_{J_{bl}}$ is the $x$-coordinate were $N$ reached $99 \%$ of the enforced value inside the penalization area ($0$ in the penalization method considered in Subsection \ref{Subsec 2 field penal}). This definition of the thickness of the boundary layer has been done by analogy with the case of a laminar flow around a flat plate (for instance, see page 30 of \cite{Sch00}).
The dashed lines represent respectively the curves $\varepsilon^{1/4}, \varepsilon^{1/2}$ and $\varepsilon$. The mesh step is $\delta x=10^{-5}$.
}
\label{Size_bl_2f_regular}
\end{figure}

 For the continuous problem, a boundary layer can be obtained by an explicit calculation of the solution (see, for instance, \cite{Pac05}, Chapter VII, Section 2 of \cite{Boy12}) or by the theoretical asymptotic expansion, as in \cite{Car03, For09}.
 The boundary layer is characterized by a non optimal convergence rate when $\varepsilon$ vanishes to $0$, which depends on the chosen norms. A boundary layer can also be observed through the plots of the solution as a quick variation of the solution between the value inside the plasma area and the enforced value inside the limiter ($N=0$ and $M=M_0$ in this case).
The numerical study presented below consists in researching evidences of boundary layers using its properties for the continuous problem (\ref{Two_fields_penal}). Hence, the need for a sufficiently accurate resolution of the problem (\ref{Two_fields_penal}).
 Usually a boundary layer has a size which decreases when the penalty parameter tends to $0$, see Fig. \ref{Size_bl_2f_regular} which indicates that the thickness decreases as $\mathcal{O}(\varepsilon)$ for $N$. Thus, if the boundary layer is too small compared to the mesh step, \emph{i.e.} when there are not enough cells in the boundary layer to resolve it, the numerical scheme does not capture the boundary layer and the rate of convergence looks sharp. So, in the presence of a boundary layer, when we study the convergence with the penalty parameter (with a fixed mesh step) we first observe a non optimal rate of convergence then the slope increases and we recover the optimal rate corresponding to the unresolved boundary layer.
   Besides, as $M \approx 1$ at the plasma-limiter interface, it appears that there is almost no wave going from the limiter to the plasma which explains why the errors in the plasma (see Fig. \ref{Error_2f_regular}) seem to be independent of the penalty parameter $\varepsilon$.
 
 In Fig. \ref{Error_2f_regular}, we notice that:
 \begin{itemize}
  \item For the $L^1$ norm in the limiter, the rate of convergence is in $\mathcal{O}(\varepsilon)$. 
  \item For the $L^1$ and the $L^2$ norms in the plasma (for $N$, $\Gamma$ and their $x$-derivatives), the errors remain almost constant. 
  \item For the $x$-derivatives in the $L^2$ norm inside the limiter, the error increases when $\varepsilon$ decreases until $\varepsilon \approx 10^{-5}$ (for $N$, the error is in $\mathcal{O}(\varepsilon^{-\frac14})$ and for $\Gamma$, the error is in $\mathcal{O}(\varepsilon^{-\frac12})$). When $\varepsilon$ is smaller than $10^{-5}$ the boundary layer is so small that there are not enough finite volume cells to resolve it, so the numerical scheme behaves as if there was no boundary layer.
  \item Inside the limiter, for $N$ and $\Gamma$ in the $L^2$ norm, we have a convergence in $\mathcal{O}(\sqrt{\varepsilon})$ for $\varepsilon \geq 10^{-5}$ which is also an evidence of the presence of a boundary layer. For $\varepsilon < 10^{-5}$, we recover a convergence in $\mathcal{O}(\varepsilon)$ due to the not sufficiently fine mesh.
 \end{itemize}

In Subsection \ref{Ssect_first_penalty}, the numerical study of the penalty presented by  Isoardi \emph{et al.} \cite{Iso10} reveals the presence of a Dirac measure next to the plasma limiter interface. The two-fields penalization presented in Subsection \ref{Subsec 2 field penal} generates a boundary layer which is not wanted.


\section[An optimal penalty method for the modified boundary conditions]{A new and optimal penalty method for the modified boundary conditions}\label{New_penalty_method}
In this paper, our goal is to provide a penalty method which is free of boundary layer, in order to ensure an optimal convergence rate when the penalization parameter $\varepsilon$ goes to $0$. Such a method is presented in this subsection. After a presentation of this method, a formal asymptotic expansion is done to provide a first evidence of the absence of boundary layer.

We are now going to describe a volume penalization method for the hyperbolic system (\ref{O1}),
that converges to the boundary problem (\ref{O2}). For the theoretical part, since we focus on a boundary value problem, we work in the domain $x<0$ as the plasma area and $x>0$ for the fictitious domain (\emph{i.e.}, the limiter set). In comparison with Section \ref{First penalty approaches}, this is just a translation in the computational domain, only considering the effects of one plasma-limiter interface.
We begin with a change of unknown to get an homogeneous Dirichlet boundary condition by defining:
\begin{align*}
 \widetilde{u}(t,\b x)=\ln \left( N(t,\b x) \right)\\
 \widetilde{v}(t,\b x)=\dfrac{\Gamma(t,\b x)}{N(t,\b x)} - M_0.
\end{align*}
The new system reads:
\begin{equation}\label{u_v_tilde}
 \left\{ \begin{array}{l}
  \partial_t \widetilde{u} + (M_0 + \widetilde{v}) \partial_x \widetilde{u} +  \partial_x \widetilde{v} =S_{\widetilde{u}}\\
  \partial_t \widetilde{v} + \partial_x \widetilde{u} + (M_0 + \widetilde{v}) \partial_x \widetilde{v} = S_{\widetilde{v}}\\
  \textsf{Boundary condition: } \widetilde{v}(.,0)=0\\
  \textsf{Initial conditions: }\widetilde{u}(0,.)  \textsf{ and } \widetilde{v}(0,.) \textsf{ are known,}
 \end{array}\right. \qquad \textnormal{ in } \mathbb{R}^+_* \times \mathbb{R}^-_*\\
\end{equation}
with the source terms $S_{\widetilde{u}}=\frac{1}{N} S_N$ and $S_{\widetilde{v}}= \frac{1}{N} S_\Gamma - \frac{M}{N} S_N$.

For Lipschitz solutions, the system (\ref{u_v_tilde}) is equivalent to (\ref{O2}) with a one side limiter.

We use a method developed in the semi-linear case by Fornet and Gu\`es \cite{For09}. Although the system (\ref{u_v_tilde}) is quasi-linear (and not semi-linear), the method can be extended to this case. An interesting feature of the method is that it yields a convergence
result without generation of a boundary layer inside the limiter. 

We assume that $M_0$ is a constant such that $0<M_0<1$.
The plasma corresponds to the region $x<0$ and the limiter is in the region $x>0$.
We denote by $\chi$ the characteristic function associated to the limiter, $\chi(x)=1$ if the point $x$ is in the limiter (i.e. $x>0$), otherwise, $\chi(x)=0$.
The penalized system is the following one:
\begin{equation}\label{penal_u_v_tilde}
 \left\{ \begin{array}{l}
  \partial_t \widetilde{u} + (M_0 + \widetilde{v}) \partial_x \widetilde{u} +  \partial_x \widetilde{v} =S_{\widetilde{u}}\\
  \partial_t \widetilde{v} + \partial_x \widetilde{u} + (M_0 + \widetilde{v}) \partial_x \widetilde{v} + \chi \dfrac{\widetilde{v}}{\varepsilon \, M_0} = S_{\widetilde{v}}\\
  \widetilde{u}(0,.)  \textsf{ and } \widetilde{v}(0,.) \textsf{ are given.}
 \end{array}\right. \qquad \textnormal{ in } \mathbb{R}^+_* \times \mathbb{R}
\end{equation}
\subsection{Asymptotic expansion of the solution}\label{asympt_exp}
We show that a formal asymptotic expansion of a regular solution can be built at any order and without boundary layer term. 
This is a first evidence of the absence of boundary layer.
At the end of this subsection, a theorem which asserts the absence of boundary layer in a slightly different problem is given.

The method consists in looking for solutions of (\ref{penal_u_v_tilde}) of the form:
\begin{align}
 \forall t\geq 0, \forall x \in \mathbb{R}, \quad \widetilde{u}_{\varepsilon}(t,x)  & \sim \left\{	\begin{array}{l}	
 																						\sum_{n=0}^{+\infty}{\varepsilon^n U^{n,-}(t,x)} \textsf{ if } x \leq 0\\
 																						\sum_{n=0}^{+\infty}{\varepsilon^n U^{n,+}(t,x)} \textsf{ if } x \geq 0 
 																					\end{array} \right. \label{asympt_u}\\
 \widetilde{v}_{\varepsilon}(t,x) & \sim \left\{	\begin{array}{l} 
 														\sum_{n=0}^{+\infty}{\varepsilon^n V^{n,-}(t,x)} \textsf{ if } x \leq 0\\
 														\sum_{n=0}^{+\infty}{\varepsilon^n V^{n,+}(t,x)} \textsf{ if } x \geq 0 ,
 													\end{array} \right. \label{asympt_v}
\end{align}
where the character $\sim$ must be read in the sense of asymptotic expansions.

We make the following assumptions:
\begin{assumption}\label{ass}
\begin{itemize}
 \item The initial condition is smooth and satisfies the compatibility conditions at the plasma limiter interface.
 \item $M_0 \in ]0,1[$ does not depend on $(t,x)$
 \item $\forall (t,x) \in \mathbb{R}^+ \times \mathbb{R}, (M_0+V^{0,\pm}(t,x))^2 < 1$
 \item $\forall n \in \mathbb{N}, U^{n,-}(.,0)=U^{n,+}(.,0)$ and $V^{n,-}(.,0)=V^{n,+}(.,0)$ . 
 \item The source terms $S_{\widetilde{u}}$ and $S_{\widetilde{v}}$ do not depend on $\widetilde{u}^{\pm}_{\varepsilon}$ and $\widetilde{v}^{\pm}_{\varepsilon}$.
\end{itemize}
\end{assumption}

The first and the fourth hypotheses are not essential: we could consider that $M_0$ varies with $(t,x)$ assuming that there exists some $c>0$ such that, for all $(t,x)$, $0<c<M(t,x)<1$.
The third hypothesis means that the continuity on $\widetilde{u}_{\varepsilon}$ and $\widetilde{v}_{\varepsilon}$ is also reported on each term of the asymptotic expansion.

\begin{proposition}\label{prop_asympt}
Under the assumption \ref{ass}, the terms of the asymptotic expansion (\ref{asympt_u}), (\ref{asympt_v}) can be constructed up to any order $n$.
\end{proposition}

Proof of Proposition \ref{prop_asympt}:

Plugging $\widetilde{u}_{\varepsilon}$ and $\widetilde{v}_{\varepsilon}$ in the penalized hyperbolic problem (\ref{penal_u_v_tilde}) gives:
\begin{align*}
 & \sum_{n=0}^{+\infty} \varepsilon^n \partial_t U^{n,\pm} + M_0 \sum_{n=0}^{+\infty} \varepsilon^n \partial_x U^{n,\pm} +\sum_{n=0}^{+\infty} \varepsilon^n V^n \sum_{k=0}^{+\infty} \varepsilon^k \partial_x U^{k,\pm} + \sum_{n=0}^{+\infty} \varepsilon^n \partial_x V^{n,\pm} =S_{\widetilde{u}}\\
 & \sum_{n=0}^{+\infty} \varepsilon^n \partial_t V^{n,\pm} + \sum_{n=0}^{+\infty} \varepsilon^n \partial_x U^{n,\pm} + M_0 \sum_{n=0}^{+\infty} \varepsilon^n \partial_x V^{n,\pm} + \sum_{n=0}^{+\infty} \varepsilon^n V^{n,\pm} \sum_{k=0}^{+\infty} \varepsilon^k \partial_x V^{k,\pm} + \dfrac{\chi}{\varepsilon} \dfrac{\sum_{n=0}^{+\infty} \varepsilon^n V^{n,\pm}}{M_0}=S_{\widetilde{v}}.
\end{align*}
Ordering the terms, we obtain:
\begin{align}
&\sum_{n=0}^{+\infty} \varepsilon^n \Big( \partial_t U^{n,\pm} + M_0 \partial_x U^{n,\pm} + \sum_{k=0}^{n}{V^{k,\pm} \partial_x U^{n-k,\pm} }+ \partial_x V^{n,\pm} \Big)=S_{\widetilde{u}} \label{Rearranged_equation_BKW1} \\
& \dfrac{1}{\varepsilon} \chi \dfrac{V^0}{M_0} + \sum_{n=0}^{+\infty} \varepsilon^n \Big( \partial_t V^{n,\pm} + \partial_x U^{n,\pm}  + M_0 \partial_x V^{n,\pm} + \sum_{k=0}^{n} V^{k,\pm} \partial_x V^{n-k,\pm} + \chi \dfrac{V^{n+1,\pm}}{M_0} \Big) = S_{\widetilde{v}}.
\label{Rearranged_equation_BKW2}
\end{align}

\emph{Term in $\varepsilon^{-1}$:}
 
 \underline{If $x>0$:}
 We have $V^{0,+}(.,x) = 0$ (for all $x>0$).

\emph{Now, we consider the induction hypothesis:} $(\mathcal{H}^n):  \forall k \leq n, (U^{k,\pm},V^{k,\pm})$ are well-defined on $]0,T[ \times \mathbb{R}$ and $V^{n+1,+}$ is well-defined on $]0,T[\times \mathbb{R}^+$ for some $T>0$ independent of $n$.

\emph{Proof of the initial assumption $(\mathcal{H}^0)$, studying the terms in $\varepsilon^0$:}

\underline{For $x<0$ ($\chi(x)=0$):}

From the equations (\ref{Rearranged_equation_BKW1}) and (\ref{Rearranged_equation_BKW2}), we have:
\begin{equation}\label{U_V_0_minus}
\left\{
\begin{array}{l}
 \partial_t U^{0,-}+ M_0 \partial_x U^{0,-} + V^{0,-} \partial_x U^{0,-} + \partial_x V^{0,-} =S_{\widetilde{u}} \qquad \textnormal{ in } \mathbb{R}^+_* \times \mathbb{R}^-_*\\
 \partial_t V^{0,-} + \partial_x U^{0,-} + M_0 \partial_x V^{0,-} + V^{0,-} \partial_x V^{0,-} = S_{\widetilde{v}}\\
 V^{0,-}(.,0)=V^{0,+}(.,0)=0 \textnormal{ (by continuity).}
\end{array}\right.
\end{equation}
Since the boundary is non characteristic and the boundary conditions are maximally strictly dissipative, the system (\ref{U_V_0_minus}) is well-posed and has a unique regular solution up to a time $T$ sufficiently small for compatible initial data, see \cite{Ben07, Rau74} and Theorem 11.1 of \cite{Ben07}. 

A definition of maximally strictly dissipative boundary condition can be found in the end of this subsection (Definition \ref{Max_dissip}).

\underline{Now, we consider the case $x>0$ ($\chi(x)=1$):}

Taking into account that $\forall (t,x) \in \mathbb{R}^+ \times \mathbb{R}^+_*,V^{0,+}(t,x)=0$, we have 
\begin{align}
 & \partial_t U^{0,+} + M_0 \partial_x U^{0,+}=S_{\widetilde{u}} \qquad \textnormal{ in } ]0,T[ \times \mathbb{R}^+_* \label{U_0_plus_eq}\\
 & U^{0,+}(.,0)=U^{0,-}(.,0) \label{U_0_plus_cl}\\
 & \partial_x U^{0,+} = S_{\widetilde{v}} - \dfrac {V^{1,+}}{M_0}. \label{U_0_plus_eq_V1}
\end{align}
The hyperbolic problem (\ref{U_0_plus_eq}), (\ref{U_0_plus_cl}) is well-posed (as $M_0>0$, we have one incoming field and one boundary condition), so $U^{0,+}$ and $V^{1,+}$ are well-defined on $ ]0,T[ \times \mathbb{R}^+_*$.

\emph{Proof of the induction step ($(\mathcal{H}^{n-1}) \Rightarrow (\mathcal{H}^n)$):}

Assuming $\mathcal{H}^{n-1}$, by using the terms in  $\varepsilon^n$, one gets:
\begin{align*}
 & \partial_t U^{n,\pm} + M_0 \partial_x U^{n,\pm} + \sum_{k=0}^{n}V^{k,\pm}\partial_x U^{n-k,\pm} + \partial_x V^{n,\pm} =0\\
 & \partial_t V^{n,\pm} + \partial_x U^{n,\pm} + M_0 \partial_x V^{n,\pm} +  \sum_{k=0}^{n}V^{k,\pm} \partial_x V^{n-k,\pm} + \chi \dfrac{V^{n+1,\pm}}{M_0} =0.
\end{align*}
\underline{If $x<0$ ($\chi(x)=0$):}
\begin{equation}\label{U_V_n_minus}
\left\{\begin{array}{l}
 \partial_t U^{n,-} + M_0 \partial_x U^{n,-} + \sum_{k=0}^{n}V^{k,-}\partial_x U^{n-k,-} + \partial_x V^{n,,-} =0\\
 \partial_t V^{n,-} + \partial_x U^{n,-} + M_0 \partial_x V^{n,-} +  \sum_{k=0}^{n} V^{k,-} \partial_x V^{n-k,-} =0.
\end{array}\right.
\end{equation}

Hence, the hyperbolic system (\ref{U_V_n_minus}) gives, by sorting according to the powers of $\varepsilon$:
\begin{equation}\label{U_V_n_minus2}
\left\{\begin{array}{l}
 \partial_t U^{n,-} + (M_0 + V^{0,-}) \partial_x U^{n,-} + \partial_x V^{n,-} = - \sum_{k=1}^{n}V^{k,-}\partial_x U^{n-k,-} \qquad \textnormal{ in } ]0,T[ \times \mathbb{R}^-_*\\
 \partial_t V^{n,-} + \partial_x U^{n,-} + (M_0 + V^{0,-}) \partial_x V^{n,-}= - \sum_{k=1}^{n}V^{k,-} \partial_x V^{n-k,-}\\
 U^{n,-}(0,.) \textsf{ and } V^{n,-}(0,.) \textsf{ are known}\\
 V^{n,-}(.,0)=V^{n,+}(.,0) \textsf{ (thanks to } (\mathcal{H}^{n-1}) \textsf{ and the continuity relation).}
 \end{array}\right.
\end{equation}
As the system (\ref{U_V_n_minus2}) is non characteristic and has maximally dissipative boundary conditions we can deduce that it is well-posed. Besides, since the system (\ref{U_V_n_minus2}) is linear, the solutions are defined on the whole interval of time $]0,T[$, for compatible initial data, see \cite{Gue90, Rau74} and Theorem 4.3 of \cite{Ben07}. Hence $U^{n,-}$ and $V^{n,-}$ are uniquely defined.

\underline{For $x>0$ ($\chi(x)=1$):}

In this area, $V^{0,+}=0$ and $V^{n,+}$ are known. From the equations (\ref{Rearranged_equation_BKW1}) and (\ref{Rearranged_equation_BKW2}):
\begin{align*}
 & \partial_t U^{n,+} + M_0 \partial_x U^{n,+} + \sum_{k=1}^{n}V^{k,+}\partial_x U^{n-k,+}  + \partial_x V^{n,+}=0\\
 & \partial_t V^{n,+} + \partial_x U^{n,+} + M_0 \partial_x V^{n,+} +  \sum_{k=1}^{n-1}V^{k,+} \partial_x V^{n-k,+} + \dfrac{V^{n+1,+}}{M_0} =0.
\end{align*}
So, we find a linear hyperbolic problem which has a unique solution:
\begin{align*}
 & \partial_t U^{n,+} + M_0 \partial_x U^{n,+} =- \sum_{k=1}^{n}V^{k,+}\partial_x U^{n-k,+} - \partial_x V^{n,+} \qquad \textnormal{ in } ]0,T[ \times \mathbb{R}^+_*\\
 & U^{n,+}(0,.) \textsf{ is known}\\
 & U^{n,+}(.,0)=U^{n,-}(.,0) \textsf{ by continuity.}\\
\end{align*}
Then, we can compute $V^{n+1,+}$ using the following relation:
\begin{equation*}
V^{n+1,+}=-M_0 \left( \partial_t V^{n,+} + \partial_x U^{n,+} + M_0 \partial_x V^{n,+} +  \sum_{k=1}^{n-1}V^{k,+} \partial_x V^{n-k,+} \right).
\end{equation*}
Hence, the property $(\mathcal{H}^n)$ is true. 
This finishes the proof of Proposition \ref {prop_asympt}.

At this stage, we have constructed an asymptotic expansion free of boundary layer. To have a complete result, we need to ensure that this asymptotic expansion converges to the solution of the limit problem at the rate $\mathcal{O}(\varepsilon)$.

To provide a rigorous result, we change the context of the problem in order to avoid a compatibility issue for the initial condition. We consider instead that the solution exists and is null in the past, \emph{i.e.} for $t \in ]-T_0,0[$ with $T_0>0$. To give a physical interpretation of this condition in the past ($t<0$), it is possible to assert that this represents the state of the scrape-off layer before the tokamak is turned on. Besides, the solution is assumed to be smooth. Indeed, the goal is to focus our study on the penalization problem, and not on the compatibility of the initial data nor on the regularity of the solution. This leads to a slightly different problem:
\begin{equation}\label{Penal_demo}
 \left\{ \begin{array}{l}
  \partial_t \left(\begin{array}{c} \widetilde{u} \\ \widetilde{v} \end{array} \right)  + \b A( \widetilde{u}, \widetilde{v} ) \, \partial_x \left(\begin{array}{c} \widetilde{u} \\ \widetilde{v} \end{array} \right) + \dfrac{\chi}{M_0 \varepsilon} \b P \left(\begin{array}{c} \widetilde{u} \\ \widetilde{v} \end{array} \right) = \left(\begin{array}{c} S_{\widetilde{u}} \\ S_{\widetilde{v}}\end{array} \right)\\
  \widetilde{u}_{|]-T_0,0[}=0 \textsf{ and } \widetilde{v}_{|]-T_0,0[}=0 
 \end{array}\right. \qquad \textnormal{ in } ]-T_0,+\infty[ \times \mathbb{R}.
\end{equation}
Where the sources terms $S_{\widetilde{u}}$ and $S_{\widetilde{v}}$ are assumed to be null for $t<0$. In our case, the matrix $\b A(\widetilde{u}, \widetilde{v} )$ and $\b P$ writes:
\begin{equation*}
\b A( \widetilde{u}, \widetilde{v} )=\left(
\begin{array}{cc} 
M_0+\widetilde{v}	&	 1\\
1				& M_0+\widetilde{v}
\end{array}\right)
\qquad \b P = 
\left(
\begin{array}{cc} 
0		& 0 		\\
0		& 1
\end{array}\right).
\end{equation*}

The introduction of the abstract matrix $\b A$ has not been done for the asymptotic expansion because this simplifies the presentation of the calculations, in the case of the system (\ref{penal_u_v_tilde}).

We consider the framework of maximally strictly dissipative boundary conditions, whose definition for a general system is recalled below:

\begin{definition}[Maximally strictly dissipative boundary conditions]\label{Max_dissip}
Consider the following hyperbolic problem of unknown $\b u : ]0,T[\times \mathbb{R}^d \to \mathbb{R}^D$:
\begin{equation}\label{Original_problem_def_msd}
\displaystyle
\left\{\begin{array}{ll} 
\partial_t \b u(t,\b x) + \sum_{j=1}^{d}{ \b A_j(\b u(t,\b x)) \partial_j \b u(t,\b x)} = \b S(t,\b x) &(t,\b x) \in ]-T_0,T[ \times \mathbb{R}^d_+ \\
\b C \b u(t,(\b x',0))=\b 0 & (t,\b x') \in ]-T_0,T[ \times \mathbb{R}^{d-1}\\
\b u_{|t<0} = \b 0
\end{array}
\right.
\end{equation}
where:
\begin{itemize}
 \item $\b S : ]-T_0,T[ \times \mathbb{R}^d_+ \to \mathbb{R}^D$
 \item For all $j$, $\b A_j$ is a symmetric matrix.
 \item $\b C$ is a constant matrix.
\end{itemize}

In the sequel, $\langle . , . \rangle$ represents the Euclidean scalar product in $\mathbb{R}^D$, and $\|.\|$ the associated norm.

The boundary conditions of (\ref{Original_problem_def_msd}) are maximally strictly dissipative if, for all $\b U \in \mathbb{R}^D$ such that  $\b C \b U=\b 0$, the quadratic form $\b V \in \mathbb{R}^D \mapsto \langle \b A_d(\b U) \b V, \b V \rangle$ has the following properties:
\begin{enumerate}
 \item $\exists \mu>0, \forall \b W \in \ker \b C, \langle \b A_d(\b U) \b W, \b W \rangle \leq -\mu \| \b W \|^2$.
 \item $\dim \ker \b C$ is maximal for the property above.
\end{enumerate}
\end{definition}

In Definition \ref{Max_dissip}, if we replace the spatial domain $\mathbb{R}^d_+$ by $\mathbb{R}^d_-$, the first property of the maximally strictly dissipative boundary condition becomes: $\exists \mu>0, \forall \b W \in \ker \b C, \langle \b A_d(\b U) \b W, \b W \rangle \geq \mu \| \b W \|^2$.

The coefficients of the matrix $\b A(.)$ are indefinitely differentiable. 
For all $\widetilde{u}, \widetilde{v}$, $\b A( \widetilde{u}, \widetilde{v} )$ is symmetric.
$\b P$ is a constant projection matrix in $\mathcal{M}_{2 \times 2}(\mathbb{R})$ satisfying: for all $( \widetilde{u}, \widetilde{v} )^t \in \ker \,\b P$, the quadratic form $\b U \mapsto  \langle \b A( \widetilde{u}, \widetilde{v} ) \b U, \b U \rangle$ is positive definite on $\ker \b P$ and $\dim \ker \b P$ is maximal for this property (in our case, $\dim \ker \b P=1$). This assumption means, by definition (see Definition \ref{Max_dissip}), that the condition $\b P \left(\begin{array}{c} \widetilde{u} \\ \widetilde{v} \end{array} \right) = \b 0$ is a maximally dissipative boundary condition at $x=0$ for the boundary value problem below:
\begin{equation}\label{BVP_demo}
 \left\{ \begin{array}{l}
  \partial_t \left(\begin{array}{c} \widetilde{u} \\ \widetilde{v} \end{array} \right)  + \b A( \widetilde{u}, \widetilde{v} ) \, \partial_x \left(\begin{array}{c} \widetilde{u} \\ \widetilde{v} \end{array} \right) = \left(\begin{array}{c} S_{\widetilde{u}} \\ S_{\widetilde{v}}\end{array} \right)\\
 \b P \left(\begin{array}{c} \widetilde{u}(.,0) \\ \widetilde{v}(.,0) \end{array} \right) = \b 0
 \end{array}\right. \qquad \textnormal{ in } ]-T_0,T[ \times \mathbb{R}^-.
\end{equation}

Then, using techniques similar to \cite{For09}, we can prove the following theorem \cite{Aup13}:

\begin{theorem}\label{The_theorem}
There exist $T>0$ sufficiently small and $\varepsilon_0>0$ such that both the penalized problem (\ref{Penal_demo}), for all $0 < \varepsilon < \varepsilon_0$, and the boundary value problem (\ref{BVP_demo}) admit a regular solution, respectively $\widetilde{u},\widetilde{v}$ on $]-T_0,T[ \times \mathbb{R}$ and $U^{0,-},V^{0,-}$ on $]-T_0,T[ \times \mathbb{R}^-$. 

Moreover, we have the following error estimates:
\begin{align*}
 \forall s \in \mathbb{N}, \quad & \| \widetilde{u} - U^{0,-}\|_{H^{s}(]-T_0,T[ \times \mathbb{R}^-_*)} = \mathcal{O}(\varepsilon)\\
 & \| \widetilde{v} - V^{0,-}\|_{H^{s}(]-T_0,T[ \times \mathbb{R}^-_*)} = \mathcal{O}(\varepsilon).
\end{align*}

\end{theorem}

Returning to the conservative variables, the penalized problem writes:
\begin{equation}\label{Penal_pb_OK}
\left\{\begin{array}{l}
 \partial_t N + \partial_x \Gamma = S_N\\
 \partial_t \Gamma + \partial_x \left(\dfrac{\Gamma^2}{N} + N\right) + \dfrac{\chi}{\varepsilon}\left(\dfrac{\Gamma}{M_0} - N\right)=S_{\Gamma}
\end{array} \right. \qquad \textnormal{ in } \mathbb{R}^+_* \times \mathbb{R}.\\
\end{equation}
As $N > 0$, the system (\ref{Penal_pb_OK}) is equivalent to (\ref{penal_u_v_tilde}) for Lipschitz solutions.
It follows from Theorem \ref{The_theorem} that (\ref{Penal_pb_OK}) admits a unique solution $N,\Gamma$ and we have:
\begin{corollary}
\begin{align*}
 \forall s \in \mathbb{N}, \quad & \| N - N^{0,-}\|_{H^{s}(]-T_0,T[ \times \mathbb{R}^-_*)} = \mathcal{O}(\varepsilon)\\
 & \| \Gamma - \Gamma^{0,-}\|_{H^{s}(]-T_0,T[ \times \mathbb{R}^-_*)} = \mathcal{O}(\varepsilon).
\end{align*}
where $N^{0,-}, \Gamma^{0,-}$ is the solution of
\begin{equation}\label{pb_OK}
\left\{\begin{array}{l}
 \partial_t N^{0,-} + \partial_x \Gamma^{0,-} = N^{0,-} S_{\widetilde{u}}\\
 \partial_t \Gamma^{0,-} + \partial_x \left(\dfrac{\Gamma^{0,- \,^2}}{N^{0,-}} + N^{0,-}\right) = N^{0,-} S_{\widetilde{v}} + \Gamma^{0,-} S_{\widetilde{u}}\\
 \dfrac{\Gamma^{0,-}(.,0)}{N^{0,-}(.,0)}=M_0 \quad \textnormal{ on } x=0.
\end{array} \right. \qquad \textnormal{ in } \mathbb{R}^+_* \times \mathbb{R}\\
\end{equation}

\end{corollary}

\subsection{Numerical tests for a one-side limiter}\label{ssect_1_side}

In the form (\ref{Penal_pb_OK}), it is easy to use the former finite volume code, as detailed in Section \ref{Numerical Schemes}, to compute the solutions of the hyperbolic problem. To avoid a stability issue, the penalized terms have to be treated implicitly.
So, for the numerical simulations, we use a semi-implicit time discretization based on the Heun scheme.
We recall from Section \ref{Numerical Schemes} that $N_i^n$ and $\Gamma_i^n$ approximate respectively the mean values of $N$ and $\Gamma$ over the cell $i$ (whose center is at $x=i \delta x$), at the time $t=t_n$. Then $f_{N,i+\frac12}^n$ and $f_{\Gamma,i+\frac12}^n$ are the numerical fluxes  for $N$ and $\Gamma$ at the time $t_n$, at the interface of the cells $i$ and $i+1$: they are computed with the VF Roe ncv scheme with second order extension as described in Subsection \ref{Appl_scheme}, see the formulas (\ref{N_i_n_lr})-(\ref{flux_Rusanov}). 
The upper index $1,n$ corresponds to the intermediate step of the Heun scheme.
\begin{align}
 & N_i^{1,n}= N_i^{n}-\dfrac{\delta t}{\delta x}\left(f_{N,i+\frac12}^n-f_{N,i-\frac12}^n \right)+\delta t\, S_{N,i}^n \label{sch_penal1}\\
 & \Gamma_i^{1,n}= \dfrac{\Gamma_i^{n}-\dfrac{\delta t}{\delta x}\left(f_{\Gamma,i+\frac12}^n-f_{\Gamma,i-\frac12}^n \right)+ \delta t \dfrac{\chi}{\varepsilon} N_i^{1,n} +\delta t\, S_{\Gamma,i}^n}{1+\delta t \dfrac{\chi}{M_0 \varepsilon}} \\
 & N_i^{n+1}= \dfrac12 (N_i^{1,n}+N_i^{n})-\dfrac{\delta t}{2 \delta x}\left(f_{N,i+\frac12}^{1,n}-f_{N,i-\frac12}^{1,n} + f_{N,i+\frac12}^n - f_{N,i-\frac12}^n \right)+\dfrac{\delta t}{2}\, (S_{N,i}^n +S_{N,i}^{n+1} ) \\
 & \Gamma_i^{n+1}= \dfrac{\dfrac12 (\Gamma_i^{1,n}+\Gamma_i^{n})-\dfrac{\delta t}{2 \delta x}\left(f_{\Gamma,i+\frac12}^{1,n}-f_{\Gamma,i-\frac12}^{1,n}+ f_{\Gamma,i+\frac12}^n-f_{\Gamma,i-\frac12}^n \right)+ \delta t \dfrac{\chi}{\varepsilon} N_i^{n+1} +\dfrac{\delta t}{2}\, (S_{\Gamma,i}^n + S_{\Gamma,i}^{n+1})}{1+\delta t \dfrac{\chi}{M_0 \varepsilon}}. \label{sch_penal4}
\end{align}
The computational domain is $[0,0.5]$ with a symmetry boundary condition at $x=0$ and the limiter set corresponds to $x \in [0.4,0.5]$ (see Fig. \ref{Half_domain1}). We study two test cases:
\begin{itemize}
 \item The first case with the regular solution (\ref{Test_sol}), that we recall here: 
\begin{equation*}
N(t,x)=\exp \left(\dfrac{-x^2}{0.16 (t+1)}\right) \qquad \Gamma(t,x)=M_0 \sin \left(\dfrac{\pi x}{0.8} \right) \exp \left(\dfrac{-x^2}{0.16 (t+1)}\right)
\end{equation*}
 and $S_N, S_{\Gamma}$ are well chosen such that $N, \Gamma$ is solution to (\ref{O2}) in the plasma area. As we have an expression of the exact solution (when $\delta x$  and $\delta t$ tend to $0$) of the limit problem (when $\varepsilon$ tends to $0$), the evaluation of the error and the convergence analysis would be easy.
 \item And with stationary solutions (as it has been studied in \cite{Iso10}).
\end{itemize}

We analyze the convergence when the penalization parameter $\varepsilon$ tends to $0$ using a uniform spatial mesh of step $\delta x=10^{-5}$. We calculate the error in $L^1$ and $L^2$ norms for $N$, $\partial_x N$, $\Gamma$ and $\partial_x \Gamma$. The goal is to confirm numerically the absence of boundary layer with an optimal rate of convergence as $\mathcal{O}(\varepsilon)$.

\begin{figure}
\begin{center}
\subfigure[Thick black: $N(1,x)$, Gray: $\Gamma(1,x)$, Narrow black: $M(1,x)$, $\varepsilon = 0.1$]{\includegraphics[scale=0.32, angle=-90, trim = 20mm 20mm 3mm 20mm, clip=true]{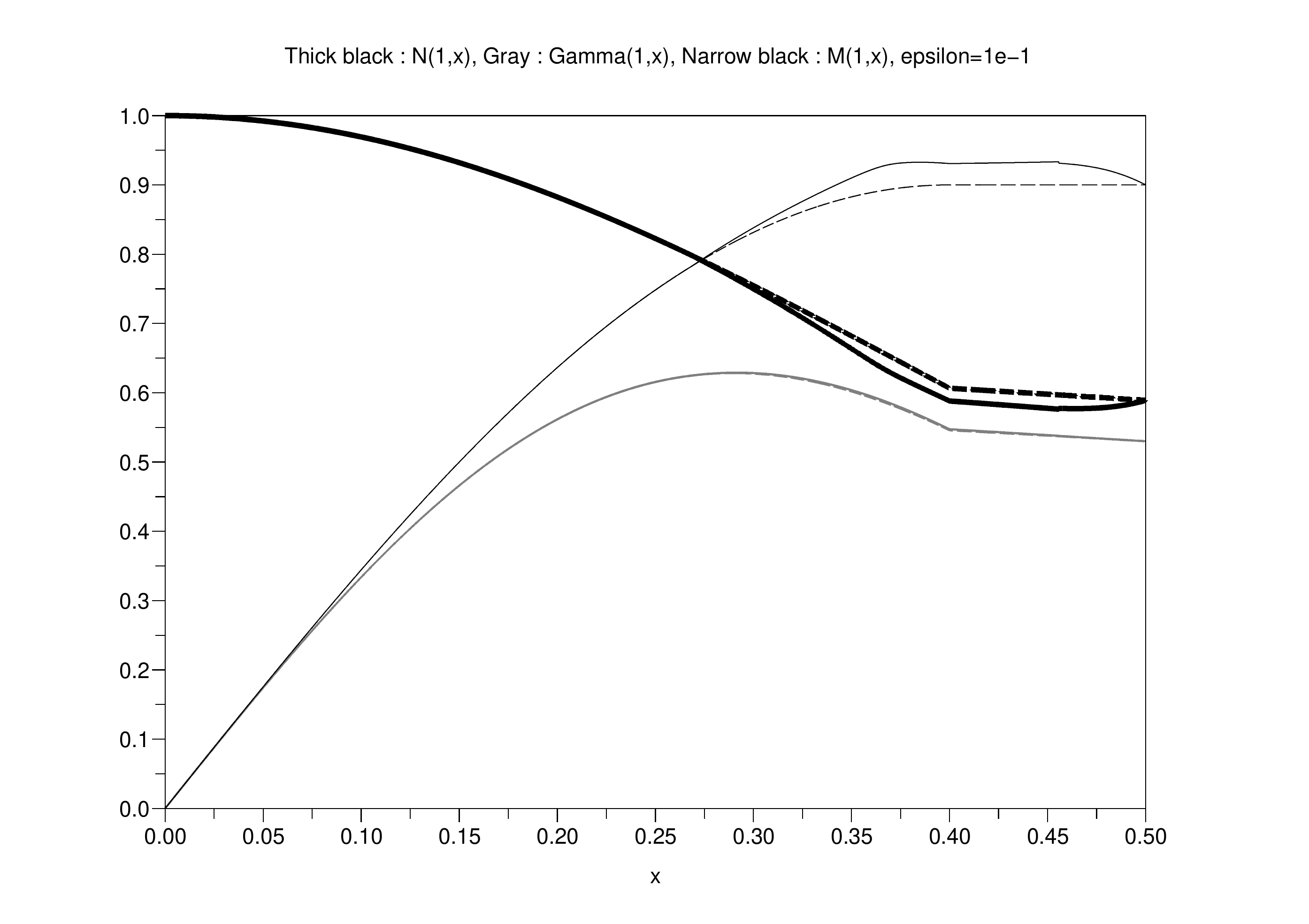}} $\quad$
\subfigure[Thick black: $N(1,x)$, Gray: $\Gamma(1,x)$, Narrow black: $M(1,x)$, $\varepsilon = 10^{-5}$]{\includegraphics[scale=0.32, angle=-90, trim = 20mm 20mm 3mm 20mm, clip=true]{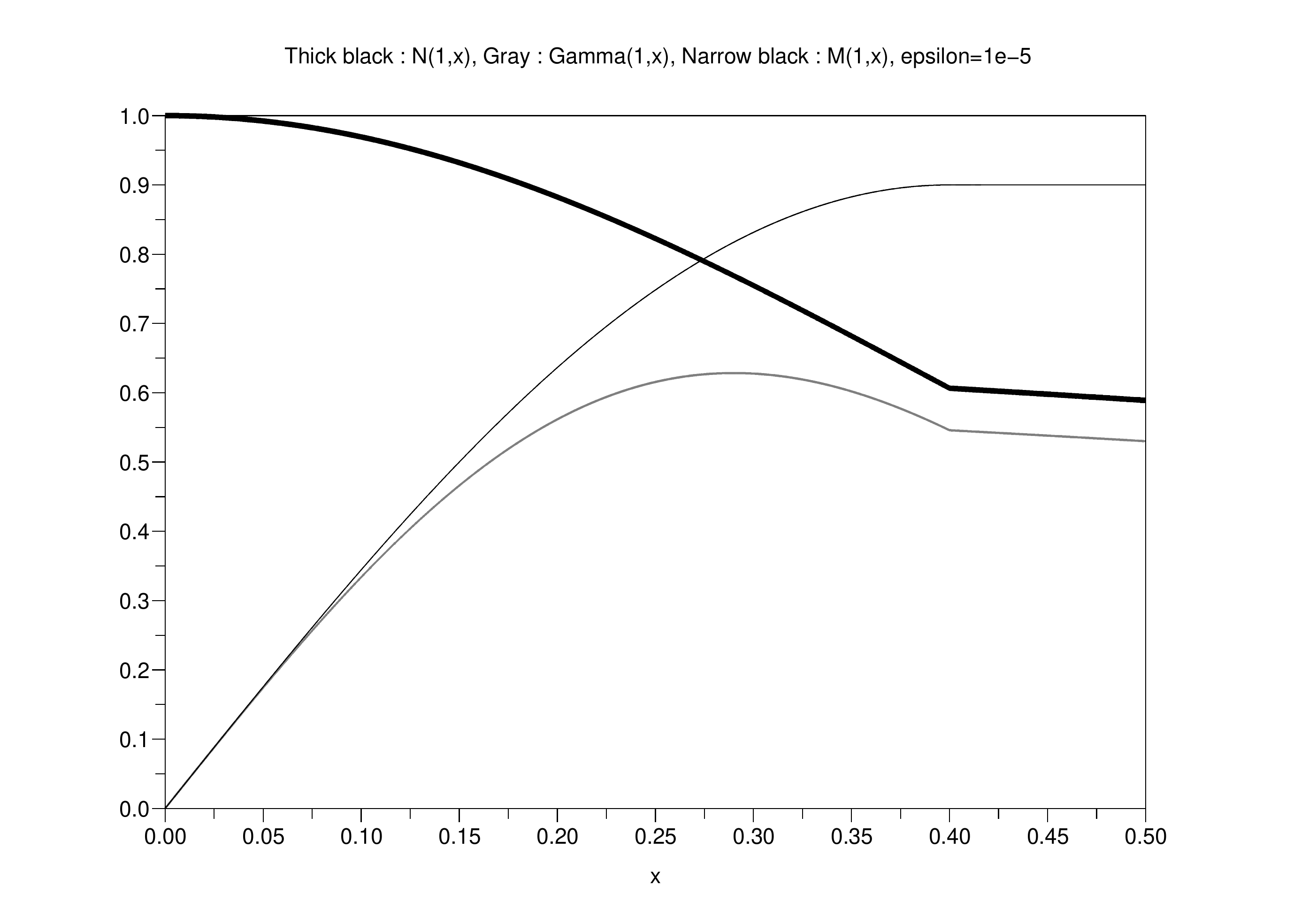}}
\end{center}
\caption{
Plot of $N$, $\Gamma$ and $M$ as functions of $x$ (at $t=1$) with the boundary layer free penalty method (at the left for $\varepsilon=0.1$ and at the right for $\varepsilon=10^{-5}$). The continuous lines represent the numerical solutions whereas the dotted lines correspond to the exact solution (when $\varepsilon$ tends to $0$). The limiter corresponds to the area $x\in[0.4,0.5]$. The mesh step is $\delta x = 10^{-5}$.
}
\label{Plot_newpen_regular}
\end{figure}

\begin{figure}
\begin{center}
\subfigure[$L^1$ error for $N$ in the plasma ($+$), $N$ in the limiter ($\times$), $\partial_x N$ in the plasma ($\circ$) and $\partial_x N$ in the limiter ($*$)]{\includegraphics[scale=0.56, trim = 7.5mm 3mm 19mm 13.3mm, clip=true]{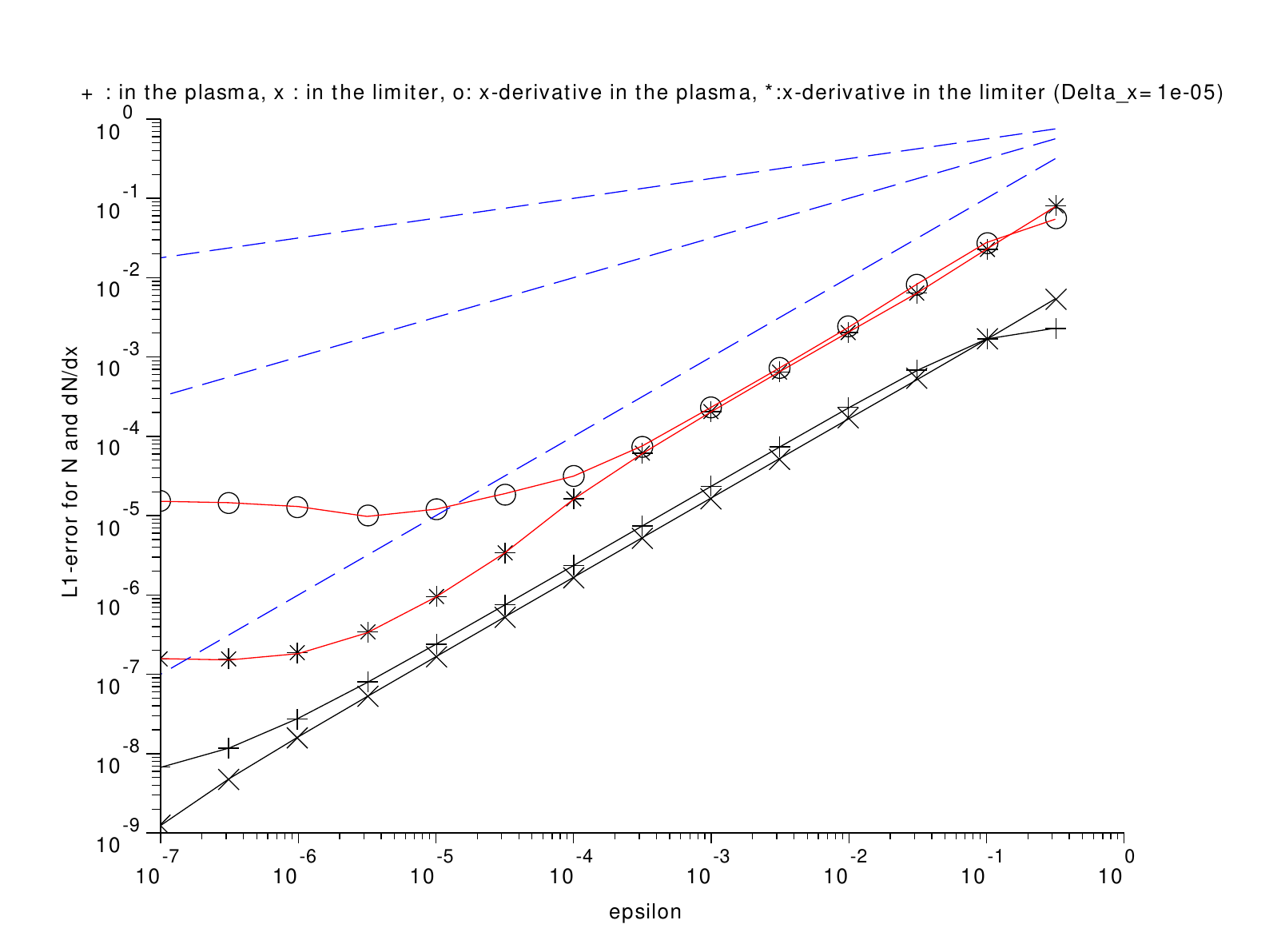}} $\quad$
\subfigure[$L^2$ error for $N$ in the plasma ($+$), $N$ in the limiter ($\times$), $\partial_x N$ in the plasma ($\circ$) and $\partial_x N$ in the limiter ($*$)]{\includegraphics[scale=0.56, trim = 7.5mm 3mm 19mm 13.3mm, clip=true]{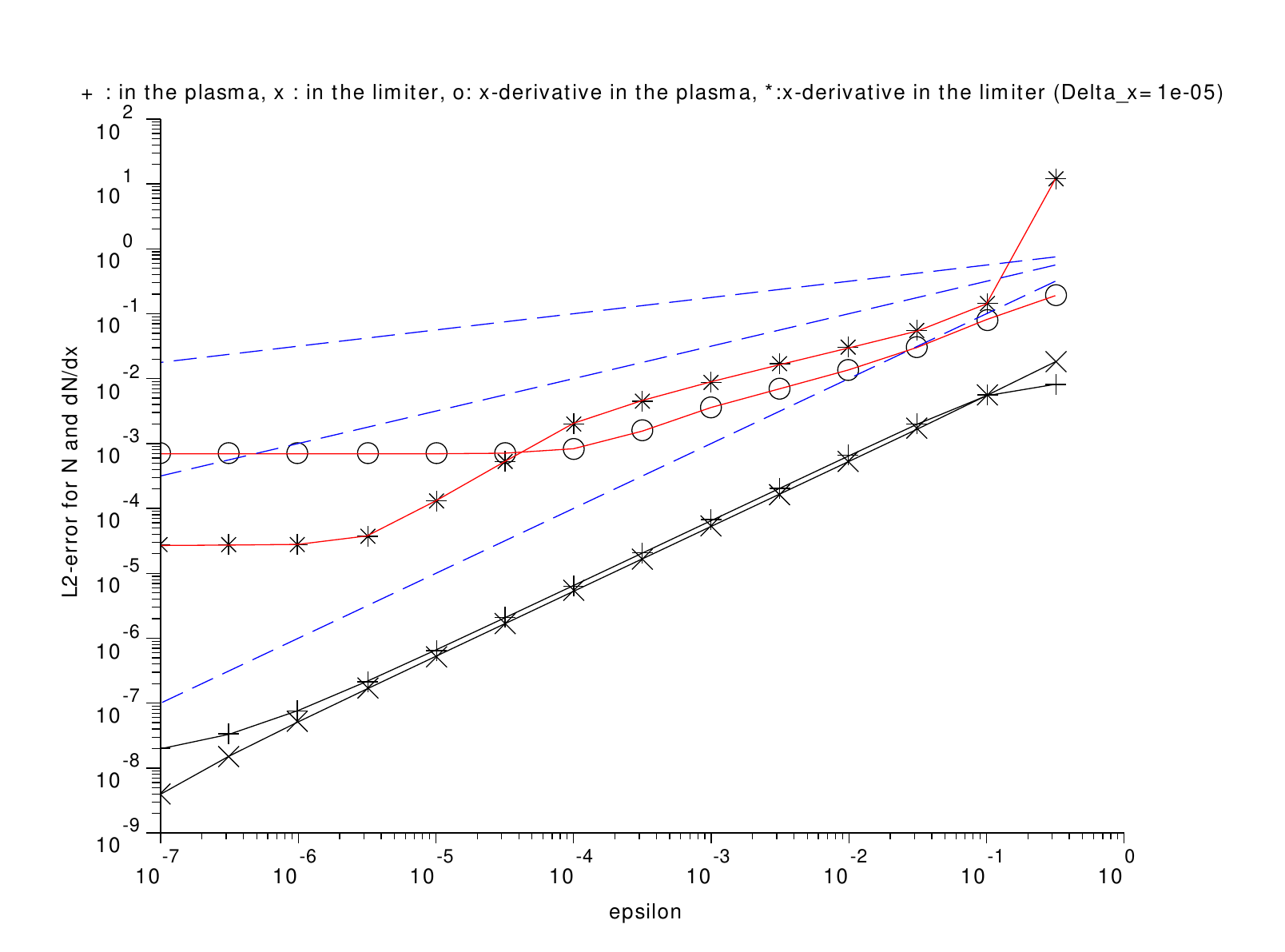}} 
\subfigure[$L^1$ error for $\Gamma$ in the plasma ($+$), $\Gamma$ in the limiter ($\times$), $\partial_x \Gamma$ in the plasma ($\circ$) and $\partial_x \Gamma$ in the limiter ($*$)]{\includegraphics[scale=0.56, trim = 7.5mm 3mm 19mm 13.3mm, clip=true]{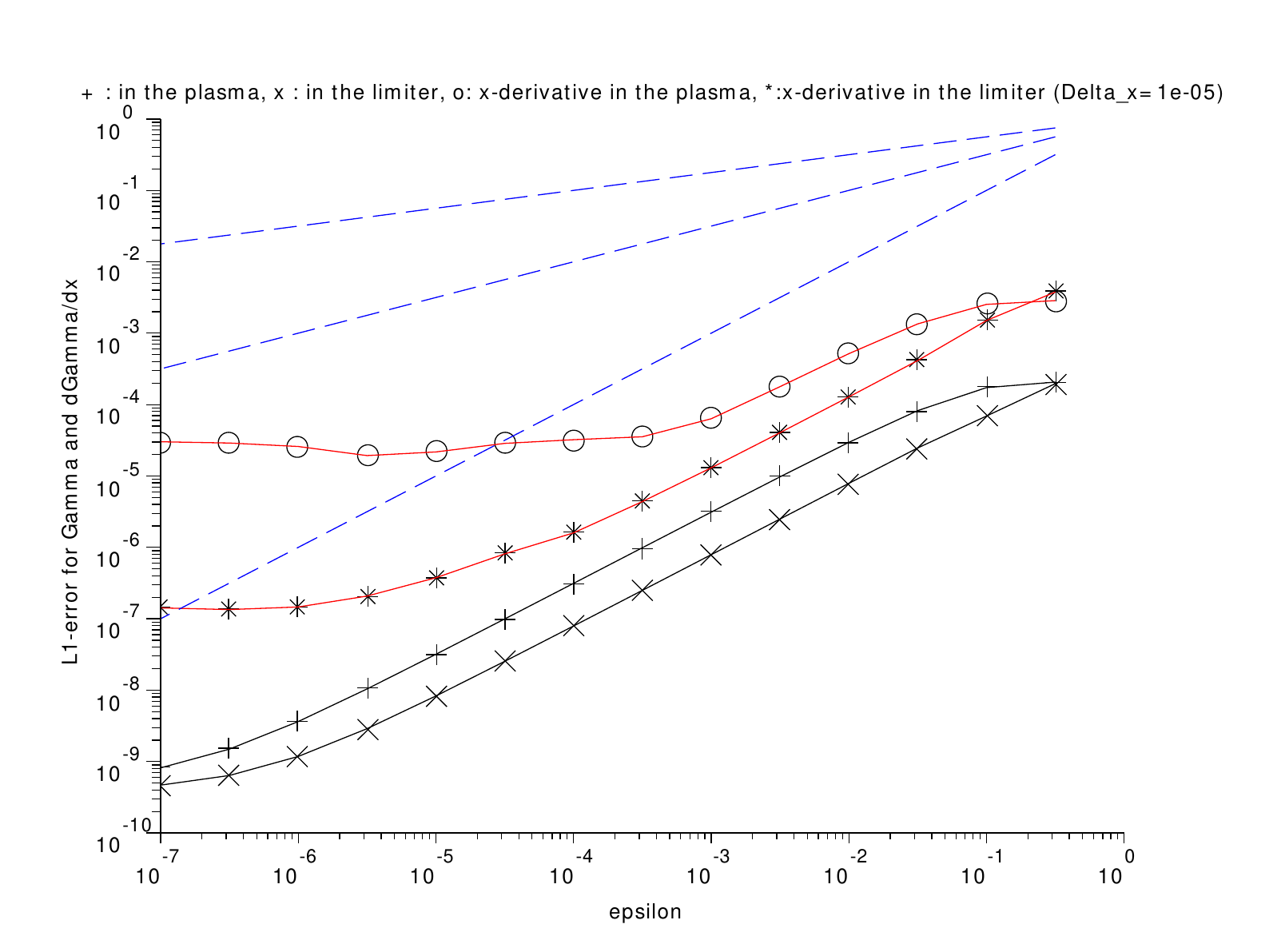}} $\quad$
\subfigure[$L^2$ error for $\Gamma$ in the plasma ($+$), $\Gamma$ in the limiter ($\times$), $\partial_x \Gamma$ in the plasma ($\circ$) and $\partial_x \Gamma$ in the limiter ($*$)]{\includegraphics[scale=0.56, trim = 7.5mm 3mm 19mm 13.3mm, clip=true]{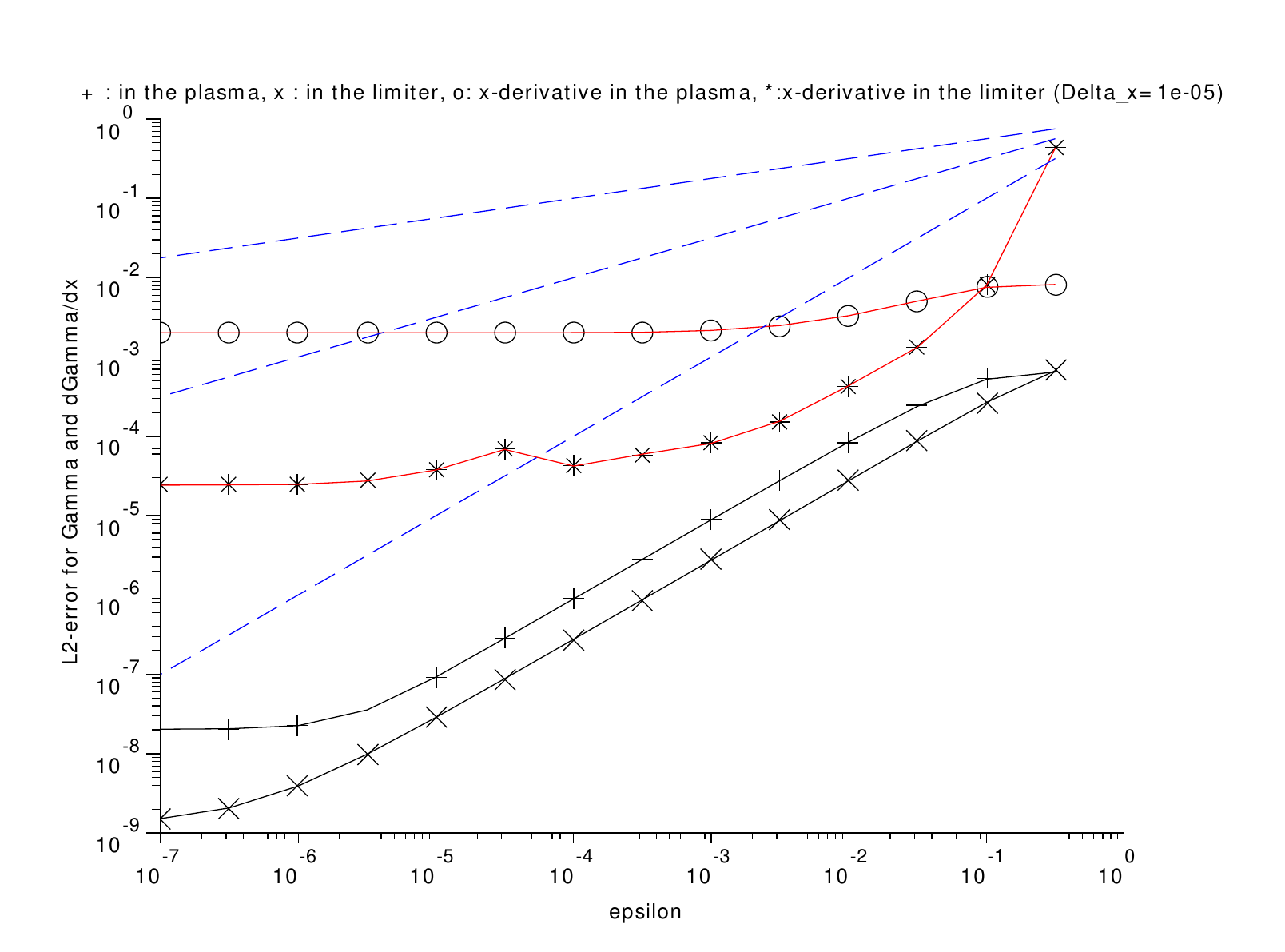}} 
\end{center}
\caption{
Errors for $N$, $\partial_x N$, $\Gamma$ (or $G$) and $\partial_x \Gamma$ in $L^1$ and $L^2$ norms with the boundary layer free penalization. The dashed lines represent the curves
$\varepsilon^{1/4}, \varepsilon^{1/2}$ and $\varepsilon$. The mesh step is $\delta x=10^{-5}$.
}
\label{Error_newpen_regular}
\end{figure}

One of the main difficulties for the implementation of the penalization, is the choice of boundary conditions at $x=0.5$ which is necessary for the numerical scheme. As only $\Gamma$ is penalized, we need a transparent boundary condition for $N$. 
For the numerical tests, the boundary condition comes from the zeroth order of the asymptotic expansion. 
For $x>0.5$, in the numerical scheme (\ref{sch_penal1})-(\ref{sch_penal4}), $N_i^{n}, \Gamma_i^{n}$, are replaced by respectively $N_{BC}(t_n,x_i), N_{BC}(t_n,x_i) M_{BC}(t_n,x_i)$ where $N_{BC}$ and $M_{BC}$ are given by the formulas (\ref{Asympt_t_inf})-(\ref{Asympt_t_sup}) given below:
\begin{itemize}
 \item If $t<\dfrac{x-0.4}{M_0}$: 
  \begin{equation}\label{Asympt_t_inf}
  \begin{array}{ll}
  N_{BC}(t,x)&=\exp(U^{0,+}(t,x))=\exp\left(-6.25 (x-t M_0)^2\right)\\
  M_{BC}(t,x)&=V^{0,+}(t,x)+M_0 = M_0.
  \end{array}
  \end{equation}
 \item Else:  
 \begin{equation}\label{Asympt_t_sup}
  \begin{array}{ll}
  N_{BC}(t,x) &= \exp(U^{0,+}(t,x))=\exp\left(-\dfrac{1}{t-\frac{x-0.4}{M_0}+1}\right)\\
  M_{BC}(t,x) &=V^{0,+}(t,x)+M_0 = M_0.
 \end{array}
  \end{equation}
\end{itemize}

\begin{figure}
\begin{center}
\includegraphics[scale=0.56, trim = 7.5mm 3mm 19mm 13.3mm, clip=true]{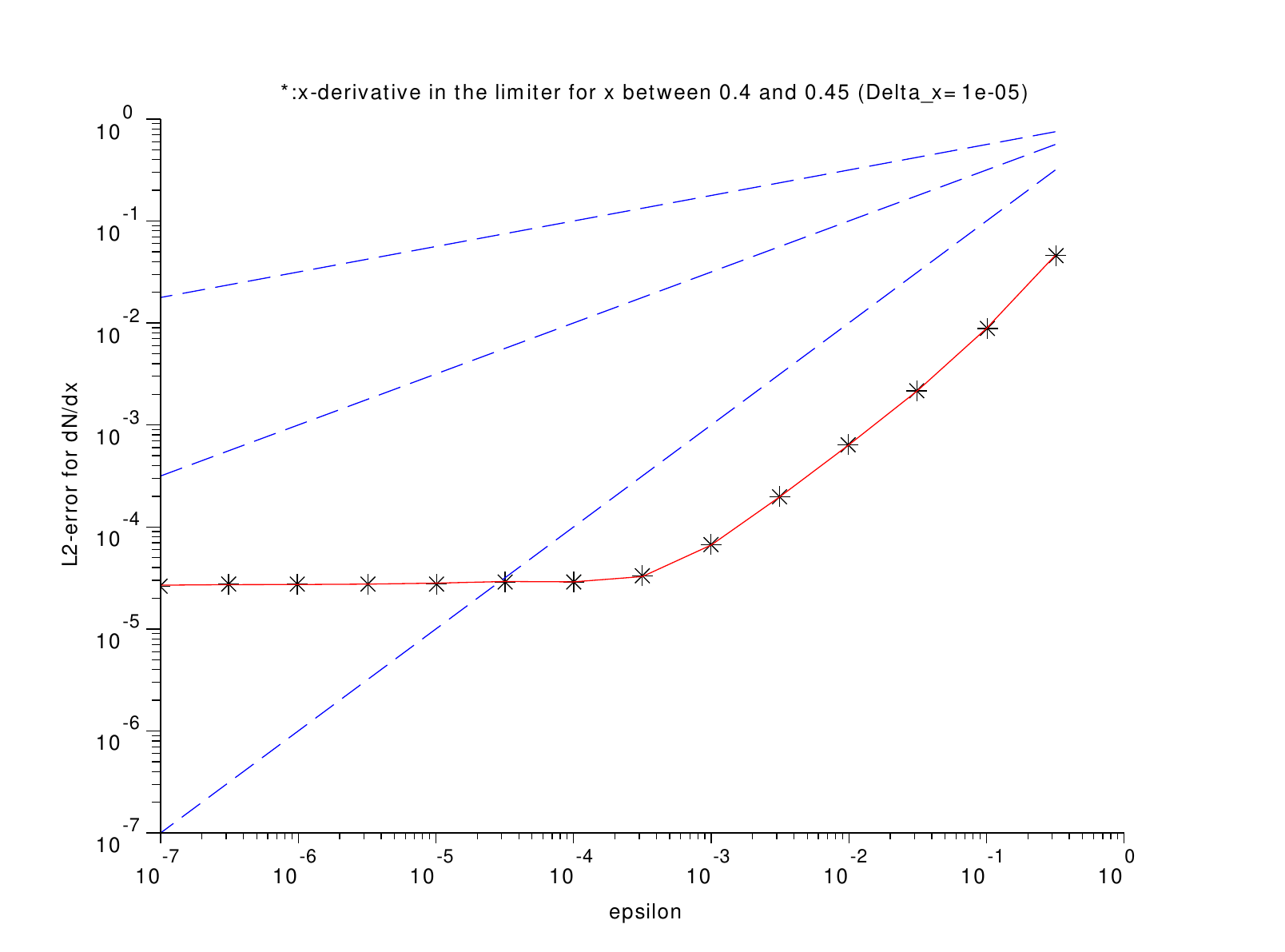}
\end{center}
\caption{
Error for $\partial_x N$ in $L^2$ norm with the boundary layer free penalization, for $0.4 \leq x \leq 0.45$. The dashed lines represent the curves
$\varepsilon^{1/4}, \varepsilon^{1/2}$ and $\varepsilon$. The mesh step is $\delta x=10^{-5}$. We observe that the optimal rate of convergence $\mathcal{O}(\varepsilon)$ is reached.
}\label{Error_newpen_regular_045}
\end{figure}

We performed the computations up to $t=1$ with an adaptive time step so that the CFL-like condition is always satisfied. The results are plotted in Fig. \ref{Plot_newpen_regular}. 
In Fig. \ref{Error_newpen_regular}, we observe that the optimal rate of convergence $\mathcal{O}(\varepsilon)$ is reached for the $L^1$ norm of error, even for the derivatives. In the $L^2$ norm, for the $x$-derivative of $N$ in the limiter, the rate of convergence seems non optimal but it can be partially explained by the difficulties to find a good artificial boundary condition at $x=0.5$ (hence the problem is not localized next to the limiter, but at $x\approx 0.5$, see Fig. \ref{Error_newpen_regular_045}). 
The anomaly inside the plasma area might also be caused by the non compatibility of the initial condition at the plasma-limiter interface.

In spite of these difficulties, the comparison of the error estimates between the two-fields penalization (Subsection \ref{Subsec 2 field penal}, Fig. \ref{Error_2f_regular}) and the penalization described in this Section (see Fig. \ref{Error_newpen_regular}) reveals that the last one has a better rate of convergence for all the norms considered.

The same numerical results in $\mathcal{O}(\varepsilon)$ are obtained if the penalty term in (\ref{Penal_pb_OK}) is replaced by
\begin{equation}\label{Penal_nat_OK}
\dfrac{\chi}{\varepsilon} \left(\dfrac{\Gamma}{N}-M_0 \right)
\end{equation}
 see \cite{Aup10}. This is due to the fact that (\ref{Penal_nat_OK}) is the penalization term in (\ref{Penal_pb_OK}) divided by $N$, which is not null in the limiter.

The stationary solution has been experimented for the problem considering that $S_N=(1-\chi) S$ and $S_{\Gamma}=0$:
$\Gamma(x)=S\,x$ and $N(x)=0.2 S\, \left(\dfrac{1}{M_0}+M_0\right)+\dfrac{S}{2} \sqrt{\left(0.4 \left(\dfrac{1}{M_0}+M_0\right) \right)^2- 4 x^2 }$.

\begin{figure}
\begin{center}
\includegraphics[scale=0.65, trim = 9mm 5mm 17mm 10mm, clip=true]{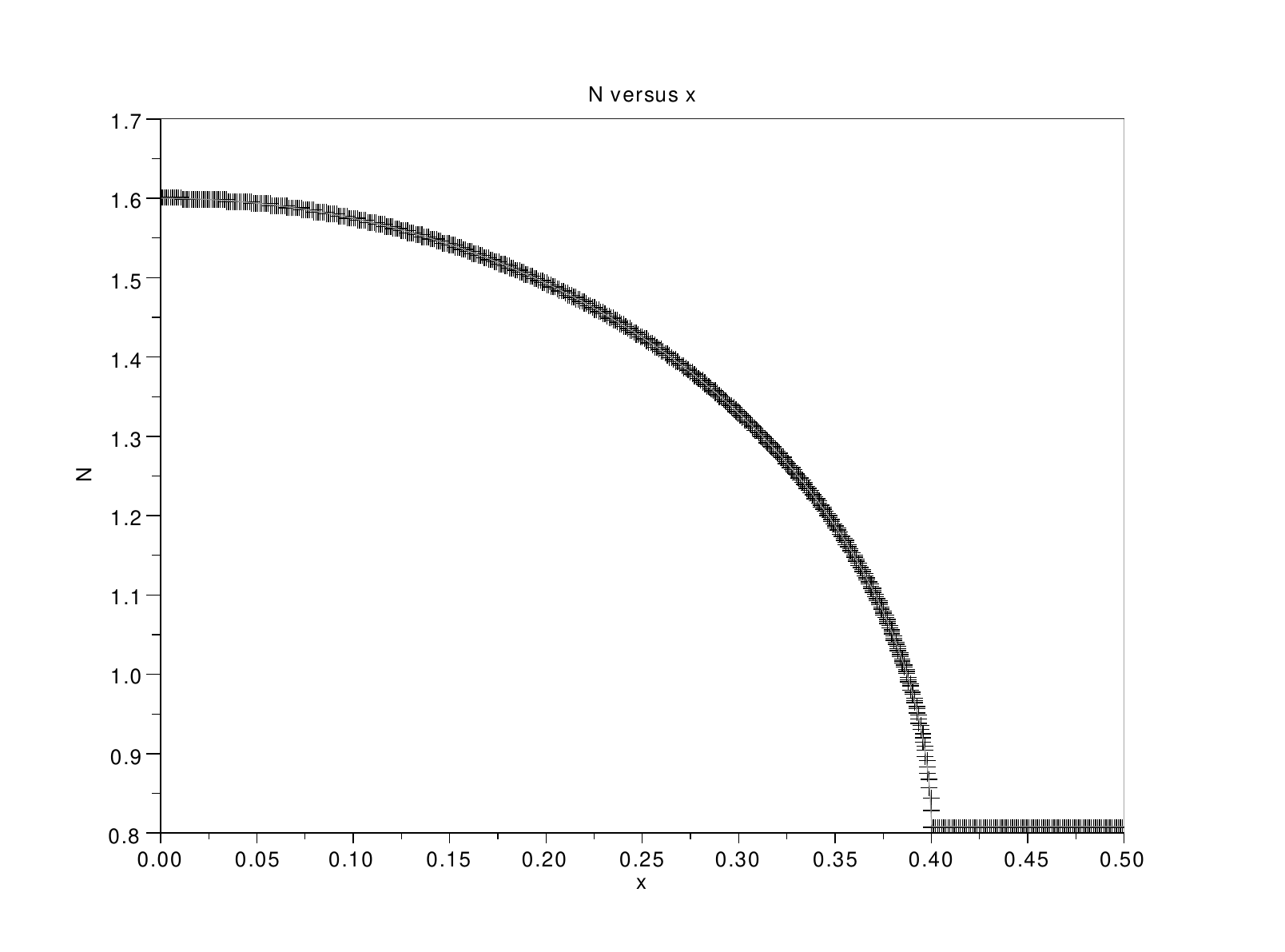}
\includegraphics[scale=0.65, trim = 9mm 5mm 17mm 10mm, clip=true]{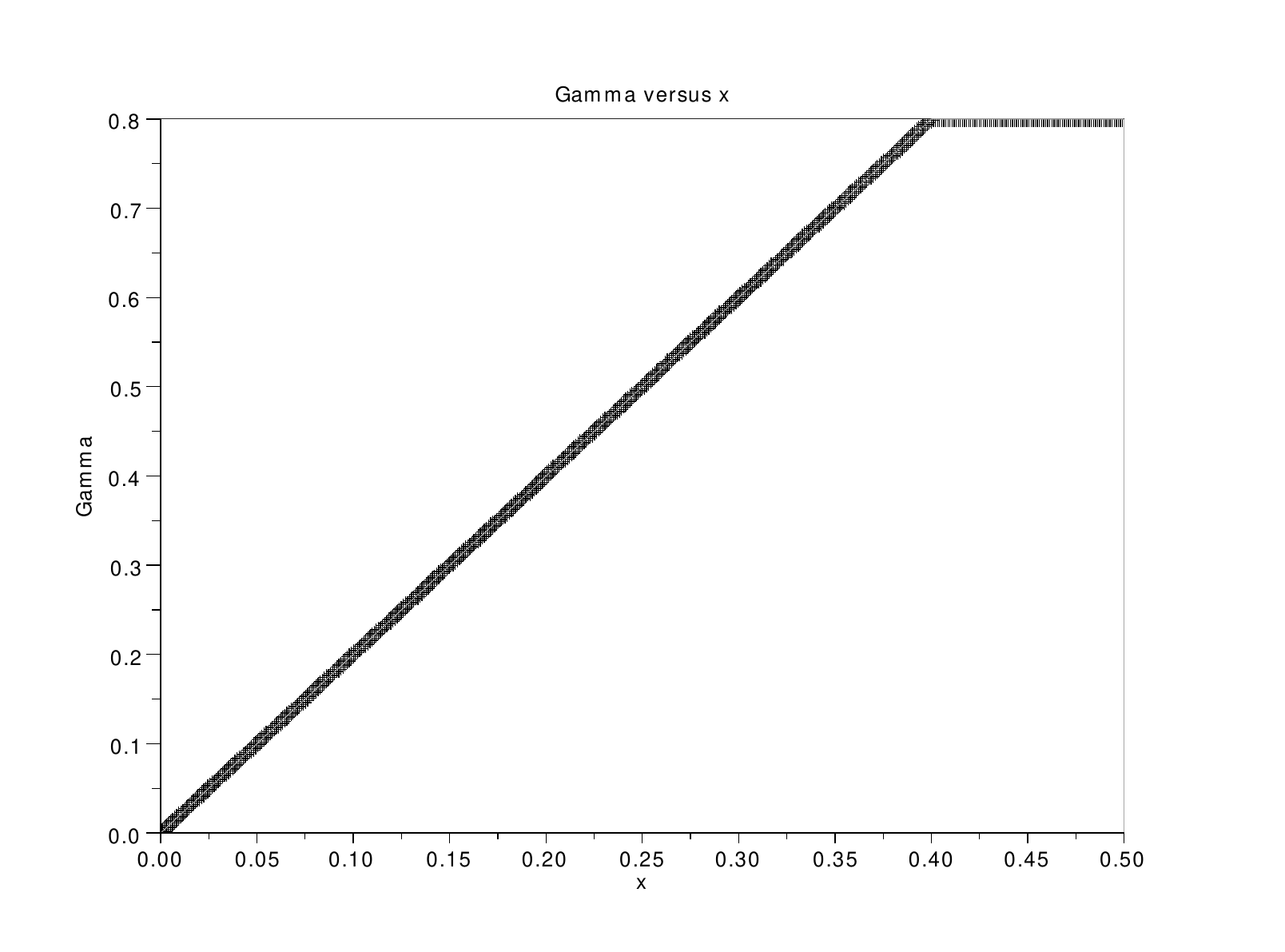}
\includegraphics[scale=0.65, trim = 9mm 5mm 17mm 10mm, clip=true]{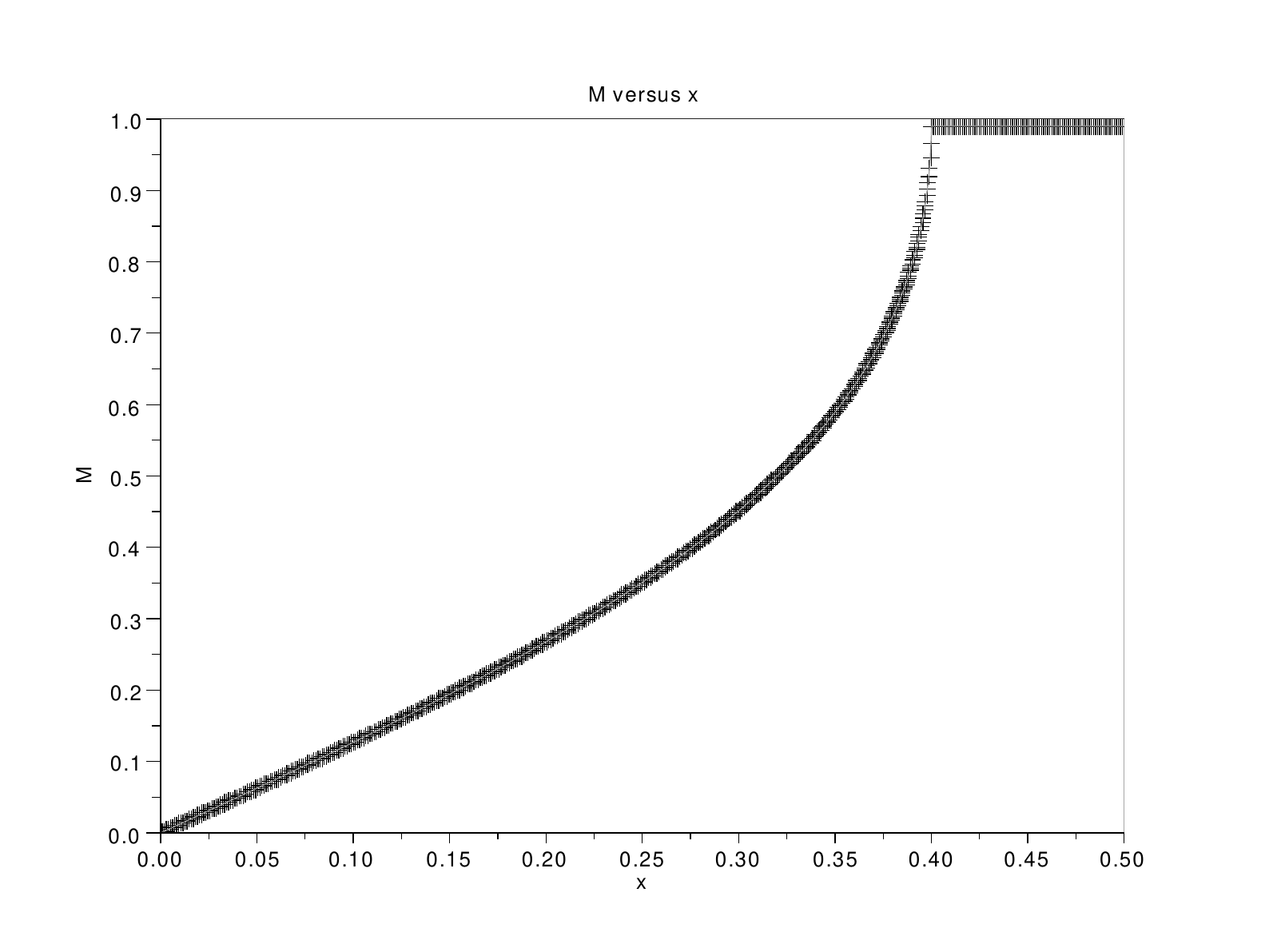}
\end{center}
\caption{
Plot of $N$, $\Gamma$ and $M$ as functions of $x$ (at $t=1$) with the boundary-layer free penalty method (the initial data is the stationary solution). The exact solution corresponds to the gray continuous line. The limiter corresponds to the area $x\in[0.4,0.5]$, see Fig. \ref{Half_domain1}. We consider $\varepsilon =10^{-3}$ and $M_0=0.99$.
}
\label{Plot_newpen_stationary}
\end{figure}

This case has been studied by Isoardi \emph{et al.} \cite{Iso10} with $M_0=1$ and with quite coarse meshes ($\delta x=0.01$) and a small penalization parameter $\varepsilon <10^{-2}$ (and for most tests $\varepsilon = 10^{-10}$), thus avoiding the issues presented in Subsection \ref{Ssect_first_penalty}. Though the boundary condition $M=1$ prevents us from using classical well-posedness theorems, the computations converge to the stationary solution. This might be due to the fact that the numerical scheme adds numerical diffusion. Tests have been conducted with $\delta x$ up to $5\cdot 10^{-4}$, $\varepsilon =10^{-3}$ or $10^{-7}$ and $M_0=0.9, 0.99$ or even $1$. The numerical solution inside the plasma domain converges towards the stationary solution and, inside the limiter, $N$ is constant but not null, as predicted by the asymptotic expansion. One test case has been represented in the Fig. \ref{Plot_newpen_stationary}.

This penalty method does not enforce $N=0$ inside the limiter, which implies that the variables $N$ and $\Gamma$ do not have any physical meaning in this area. Another point of view is to consider that this penalty method does not model the plasma-limiter interface but represents the boundary between the pre-sheath and the sheath.

\subsection{Penalization for a two-sides limiter}\label{Subsec_2sides_penal}

The penalty method presented in this Section assumes that only one side of the limiter is in interaction with the plasma. To provide a more realistic model as presented in \cite{Iso10}, we consider now that the limiter has two sides. As the $x$-axis follows a magnetic field line, which is a loop interrupted by the limiter, in this configuration, we can impose periodic boundary conditions.

\begin{figure}
\begin{center}
\includegraphics[scale=0.50, trim = 10mm 200mm 40mm 10mm, clip=true]{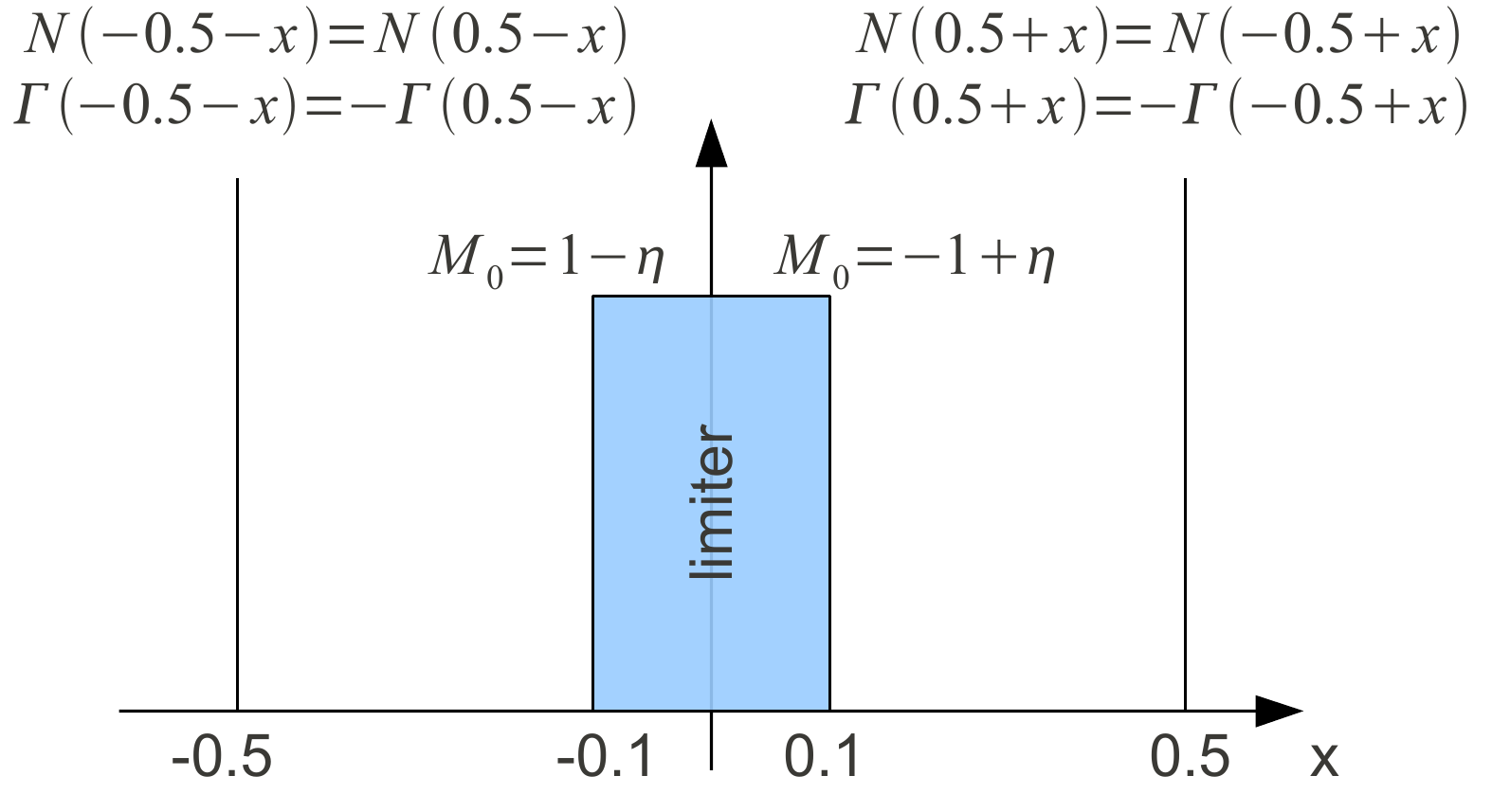}
\end{center}
\caption{
Representation of the simulated domain with a two-sides limiter at the center.
}
\label{Complete_domain_center}
\end{figure}

From the term of order $0$ in the asymptotic expansion, we deduce that information is propagating from the plasma-limiter interface to the interior of the limiter. But the limiter has now two faces, and no information must pass through it. To avoid this phenomenon, we multiply the flux by a smooth function $\alpha$, which is null in the central area inside the limiter and $\alpha=1$ elsewhere. The system obtained is still well-posed because of the smoothness of $\alpha$. For the numerical tests, the domain is $x \in ]-0.5,0.5[$ and the limiter set corresponds to $x \in [-0.1,0.1]$ (see Fig. \ref{Complete_domain_center}).

Now the penalized hyperbolic problem reads:
\begin{equation}\label{Penal_pb_2faces}
\left\{\begin{array}{l}
 \partial_t N + \partial_x \left(\alpha \Gamma\right) = S_N\\
 \partial_t \Gamma + \partial_x \left( \alpha\left(\dfrac{\Gamma^2}{N} + N\right)\right) + \sign(-x) \dfrac{\chi}{\varepsilon}\left(\dfrac{\Gamma}{M_0} - N\right)=S_{\Gamma}
\end{array} \right. \qquad \textnormal{ in } \mathbb{R}^+_* \times ]-0.5,0.5[.\\
\end{equation}

For $\alpha$, we use the following expression:
\begin{equation*}
 \alpha(x)=\left\{\begin{array}{l}
 1 \textsf{ if } x \in ]-0.5,-0.075] \\
 \dfrac12 \tanh\left(0.060 \left(-\dfrac{1}{x-0.015}-\dfrac{1}{x-0.075}\right)\right)+\dfrac12 \textsf{ if } \in ]-0.075,-0.015[\\
 0 \textsf{ if } x \in  ]-0.015,0.015[\\
 \dfrac12 \tanh\left(0.060 \left( \dfrac{1}{x+0.015}+\dfrac{1}{x+0.075} \right)\right)+\dfrac12 \textsf{ if } x \in ]0.015,0.075[\\
 1 \textsf{ if } x \in ]0.075,0.5]. \\
 \end{array} \right.
\end{equation*}

Following the idea developed by Greenberg and Le Roux \cite{Gre96}, we consider, for the implementation of the solver, that $\alpha$ is an unknown of the system. The new hyperbolic system in the non-conservative form reads:

\begin{equation}\label{Equa_Le_Roux}
\partial_t \left( \begin{array}{c} N \\ \Gamma \\ \alpha \end{array} \right) + 
 \left( \begin{array}{ccc} 0 & \alpha & \Gamma \\ \alpha \, \left( 1 - \dfrac{\Gamma^2}{N^2}\right) & 2 \alpha \dfrac{\Gamma}{N}  & \dfrac{\Gamma^2}{N}+N \\ 0 & 0 & 0\end{array} \right)
\partial_x \left( \begin{array}{c} N \\ \Gamma \\ \alpha \end{array} \right)+
\dfrac{\chi}{\varepsilon} \left( \begin{array}{c} 0 \\ \dfrac{\Gamma}{M_0}-N \\ 0 \end{array} \right) =
\left( \begin{array}{c} S_N \\ S_\Gamma \\ 0 \end{array} \right).
\end{equation}

The system (\ref{Equa_Le_Roux}) is solved using a VFRoe ncv scheme with MUSCL reconstruction, slope limiter and the modified Heun time discretization.

\begin{figure}
\begin{center}
\includegraphics[scale=0.32, angle=-90, trim = 0mm 20mm 0mm 20mm, clip=true]{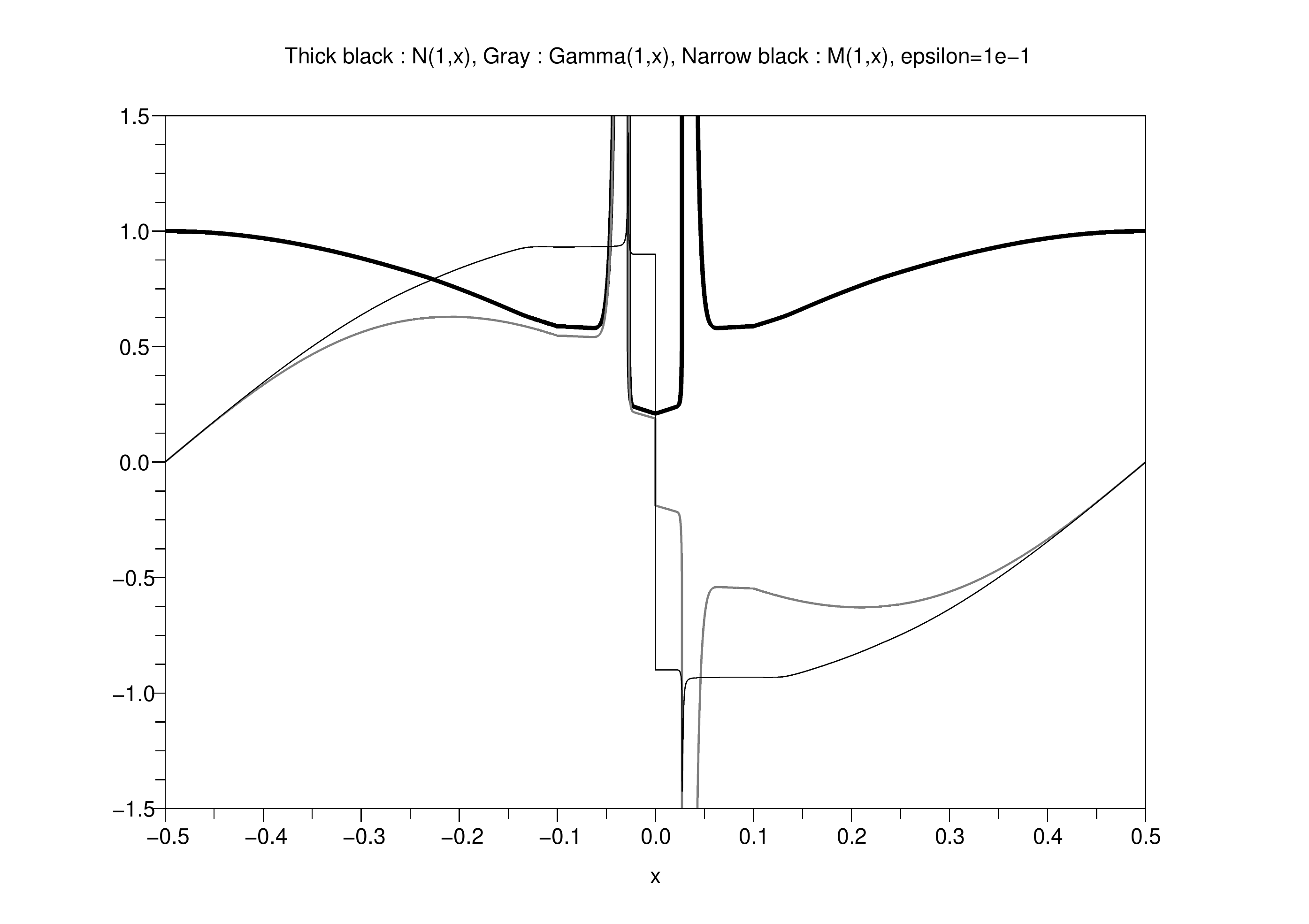}
\includegraphics[scale=0.32, angle=-90, trim = 0mm 20mm 0mm 20mm, clip=true]{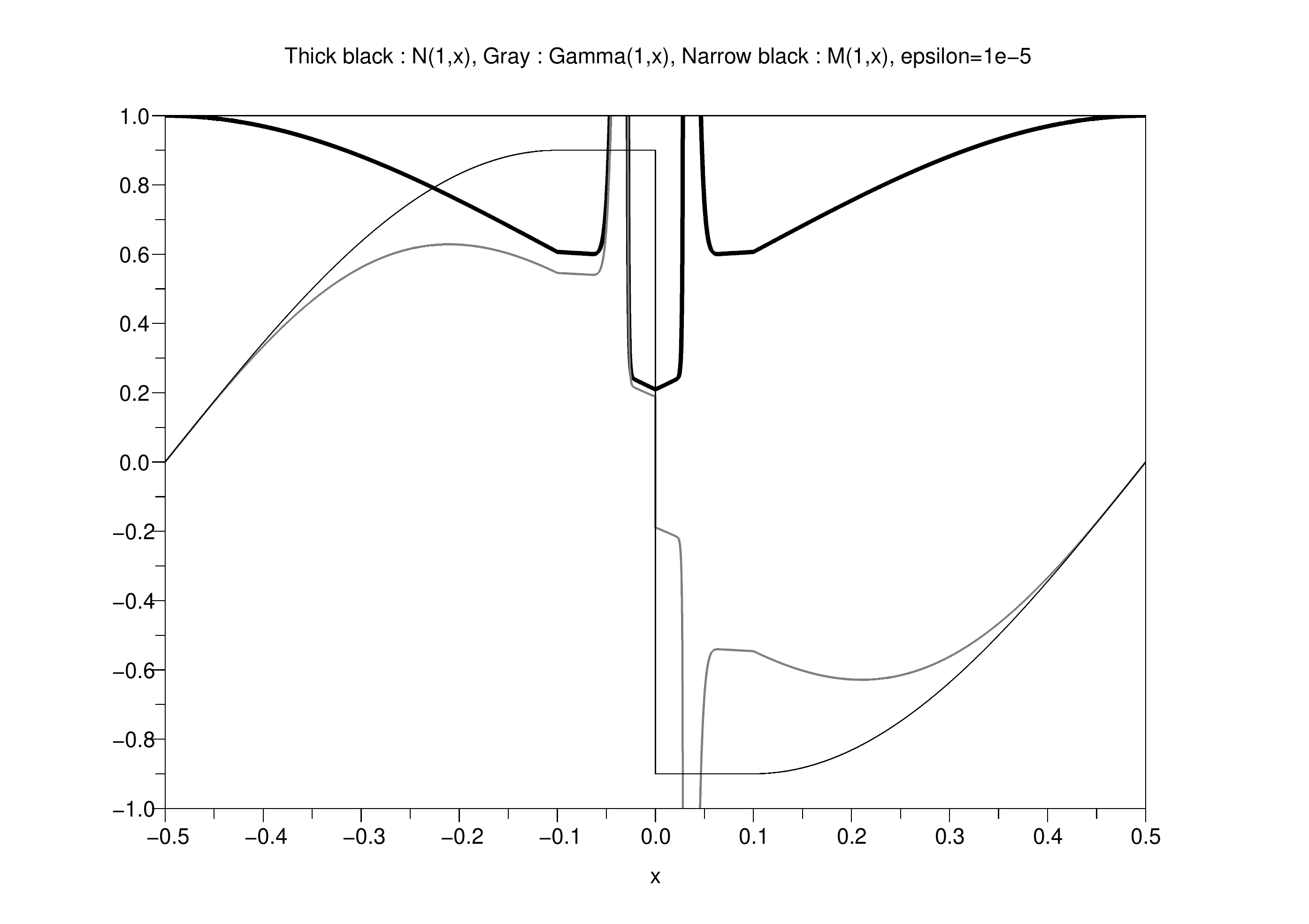}
\end{center}
\caption{
Plot of $N$, $\Gamma$ and $M$ as functions of $x$ (at $t=1$) with the boundary layer free penalty method for the two sides limiter (at the left for $\varepsilon=0.1$ and at the right for $\varepsilon=10^{-5}$). The limiter corresponds to the area $x\in[-0.1,0.1]$.
For $\varepsilon=0.1$, we have $\max(N)=115.65$ and $\max(|\Gamma |)=122.72$.
For $\varepsilon=10^{-5}$, we have $\max(N)=168.91$ and $\max(|\Gamma |)=152.04$.
}
\label{Plot_newpen_regular_2_bords}
\end{figure}

From the numerical tests (see Fig. \ref{Plot_newpen_regular_2_bords}), in the areas where $\alpha$ is close to $0$ (but not equal to $0$), we observe peaks for the variable $N$. This is not a Dirac measure though: even if the values are large, they remain bounded as the resolution increases. This can be intuitively explained by the order $0$ of the asymptotic expansion inside the limiter set: $N$ is mainly governed by a simple transport equation of speed $\alpha M_0$ towards the center of the limiter. So the variable $N$ is transported at the speed $\alpha M_0$ from the boundary and is stopped when $\alpha$ decreases to $0$. Hence, this generates the two accumulation areas observed.

\begin{figure}
\begin{center}
\subfigure[$L^1$ error for $N$ in the plasma ($+$), $N$ in the limiter ($\times$), $\partial_x N$ in the plasma ($\circ$) and $\partial_x N$ in the limiter ($*$)]{\includegraphics[scale=0.56, trim = 7.5mm 3mm 19mm 13.3mm, clip=true]{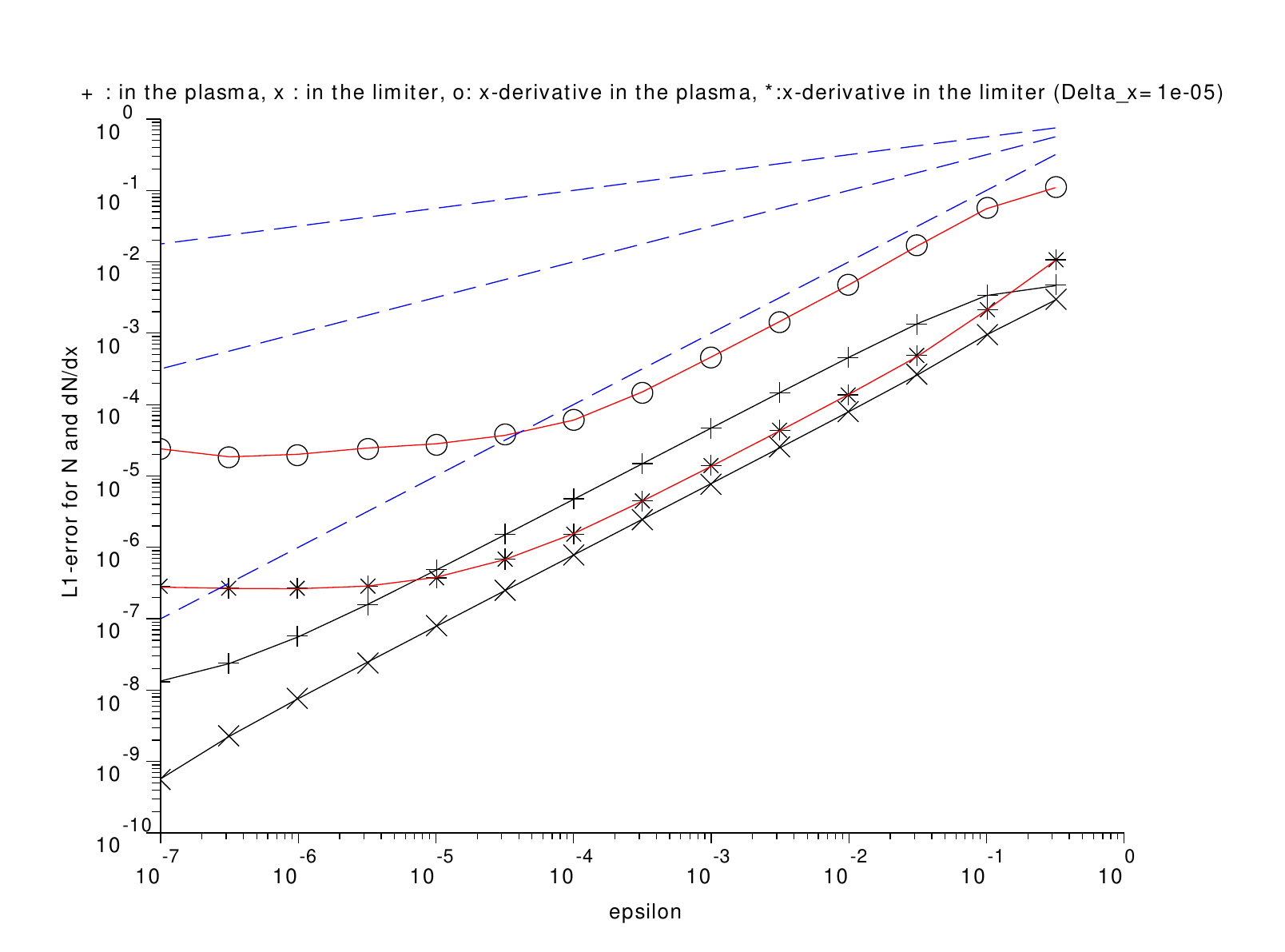}} $\quad$
\subfigure[$L^2$ error for $N$ in the plasma ($+$), $N$ in the limiter ($\times$), $\partial_x N$ in the plasma ($\circ$) and $\partial_x N$ in the limiter ($*$)]{\includegraphics[scale=0.56, trim = 7.5mm 3mm 19mm 13.3mm, clip=true]{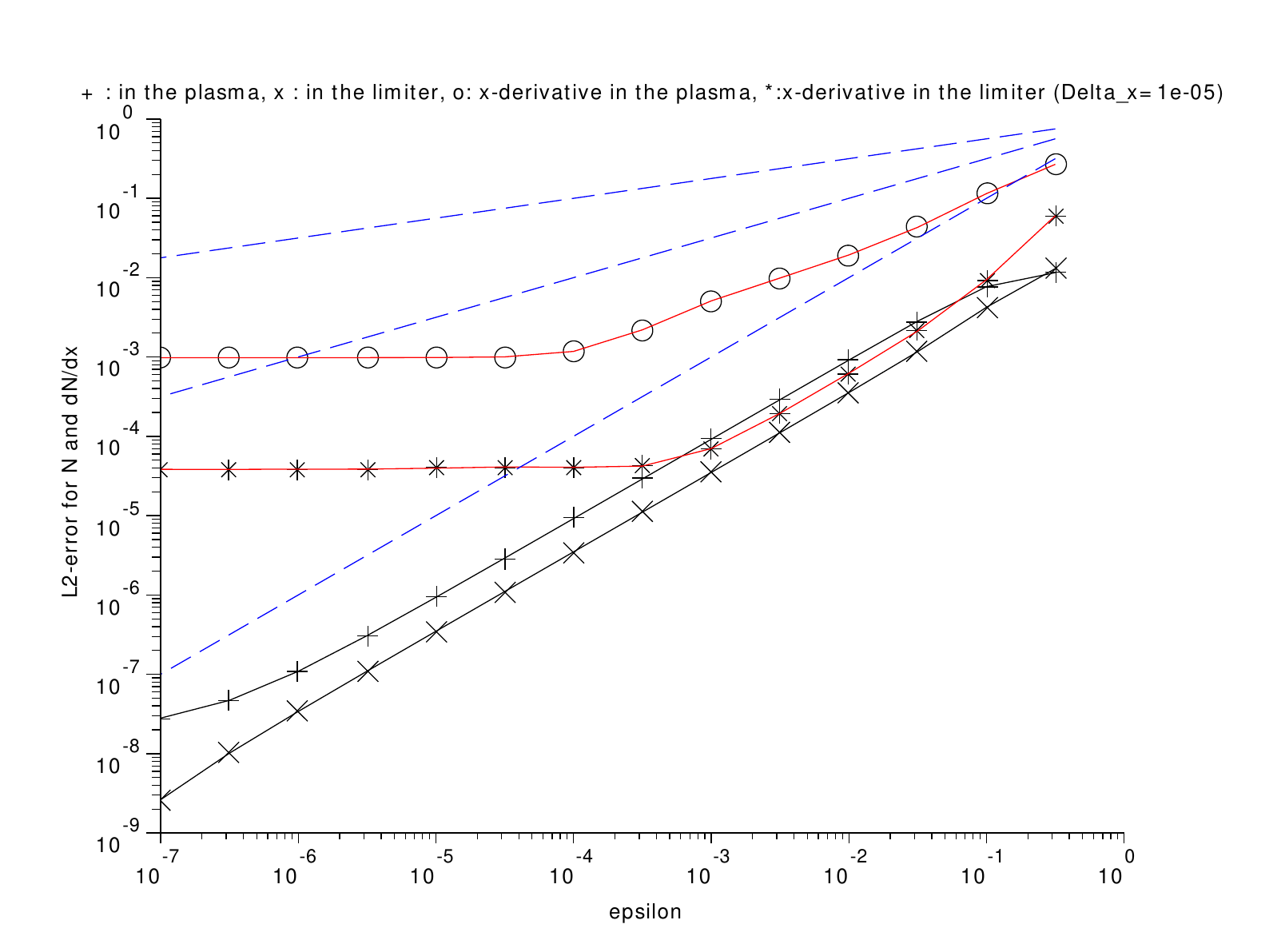}} 
\subfigure[$L^1$ error for $\Gamma$ in the plasma ($+$), $\Gamma$ in the limiter ($\times$), $\partial_x \Gamma$ in the plasma ($\circ$) and $\partial_x \Gamma$ in the limiter ($*$)]{\includegraphics[scale=0.56, trim = 7.5mm 3mm 19mm 13.3mm, clip=true]{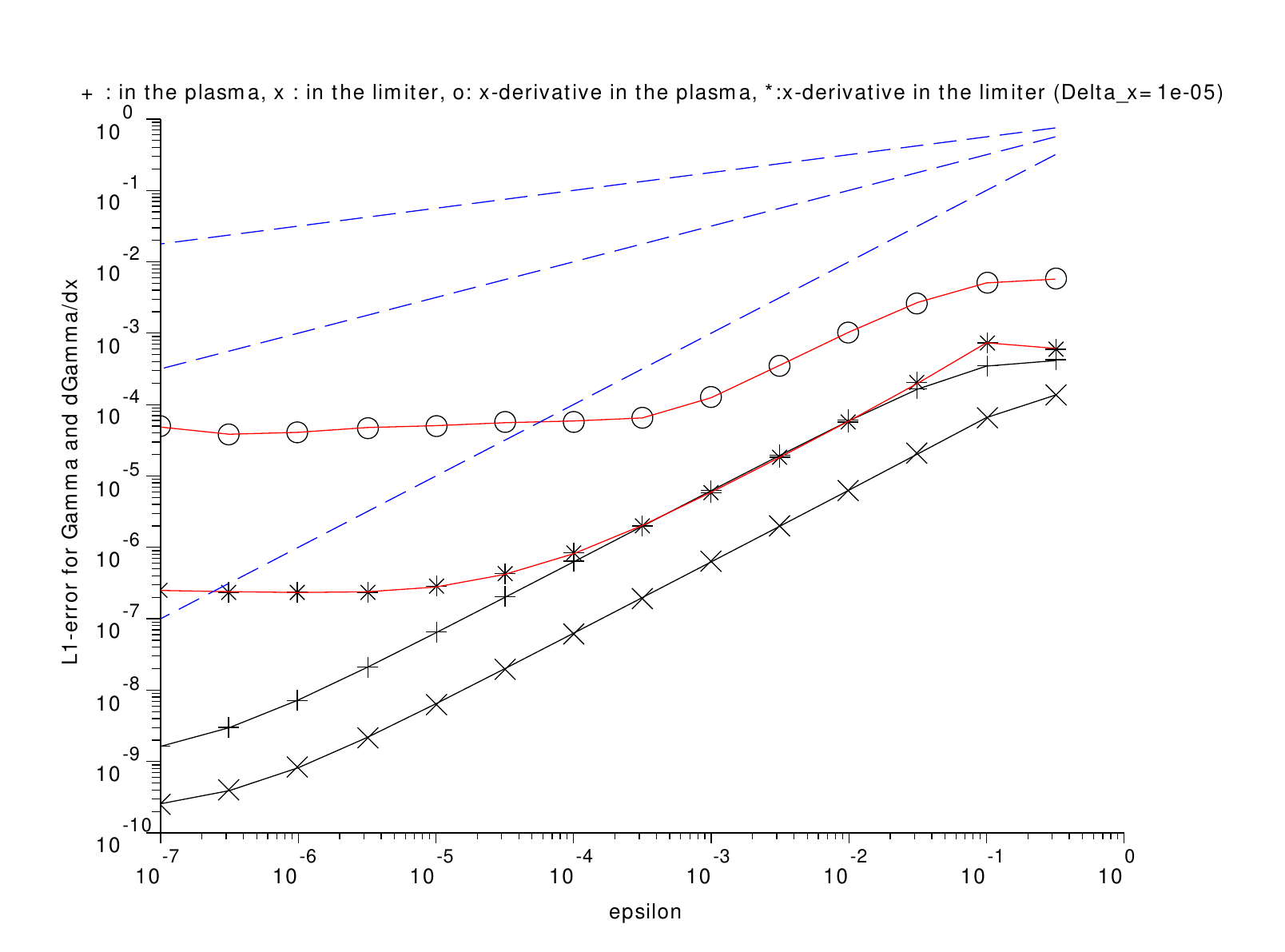}} $\quad$
\subfigure[$L^2$ error for $\Gamma$ in the plasma ($+$), $\Gamma$ in the limiter ($\times$), $\partial_x \Gamma$ in the plasma ($\circ$) and $\partial_x \Gamma$ in the limiter ($*$)]{\includegraphics[scale=0.56, trim = 7.5mm 3mm 19mm 13.3mm, clip=true]{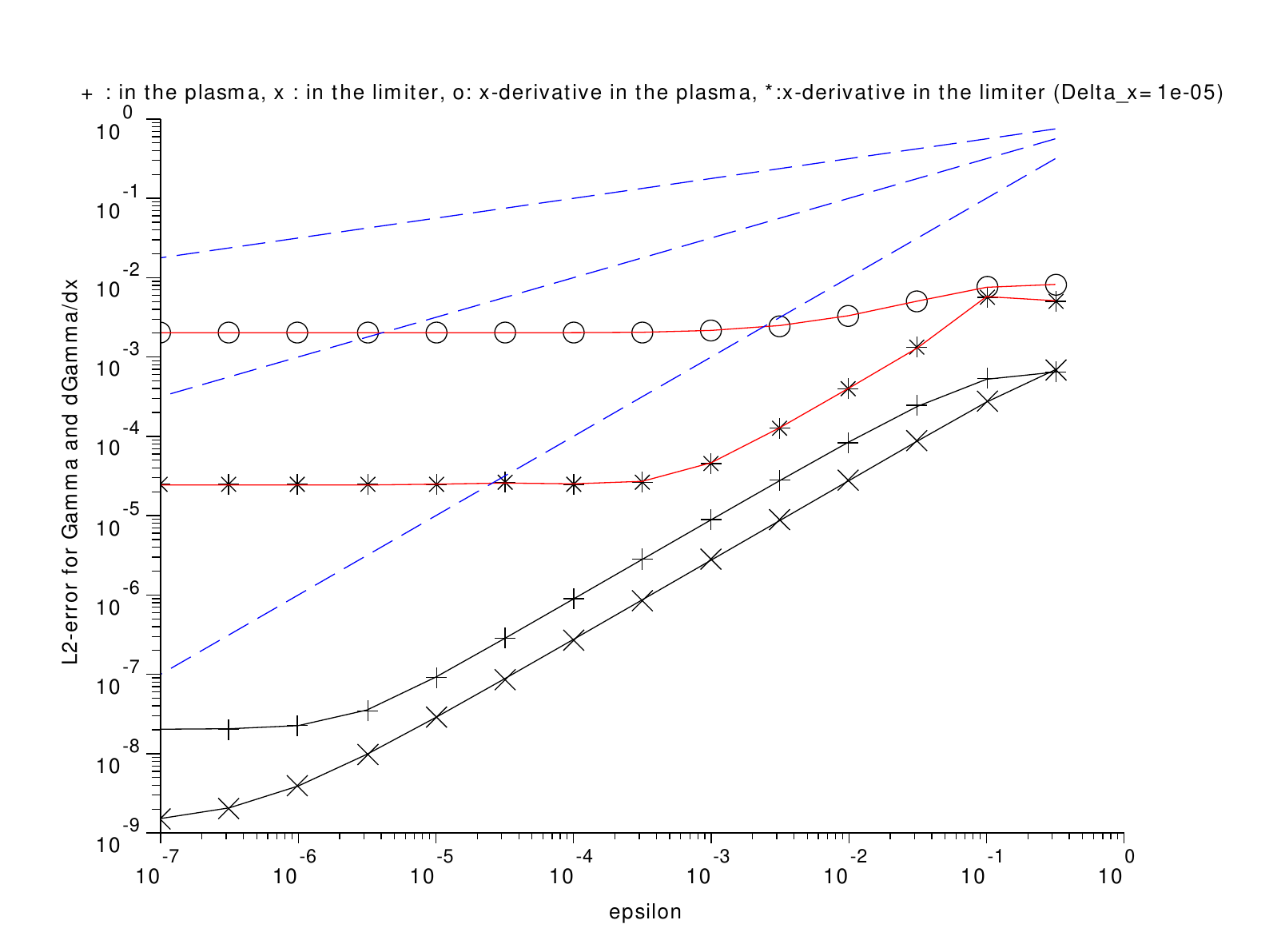}} 
\end{center}
\caption{
Errors for $N$, $\partial_x N$, $\Gamma$ (or $G$) and $\partial_x \Gamma$ in $L^1$ and $L^2$ norms with the boundary layer free penalization and the two sides limiter configuration (see Fig. \ref{Complete_domain_center}). The dashed lines represent the curves
$\varepsilon^{1/4}, \varepsilon^{1/2}$ and $\varepsilon$.
The error in the plasma area is estimated in the set $x \in ]-0.5,-0.1[\cup]0.1,0.5[$. In the limiter, the place where the error is computed is $x \in [-0.1,-0.075]\cup[0.075,0.1]$. So, the part where $\alpha(x)\neq 1$ is excluded.
Here, $M_0=0.9$ and the mesh step is $\delta x=10^{-5}$.
}
\label{Error_newpen_regular_2bords}
\end{figure}

The numerical convergence analysis when the penalization parameter tends to $0$ (see Fig. \ref{Error_newpen_regular_2bords}) leads to the same conclusions as for the one-side limiter model.

\subsection{Analysis when $|M_0|$ tends to $1$}\label{ssect_M_O_to_1}
In the Section \ref{Presentation of the hyperbolic system}, we modified the value of $|M_0|$ from $1$ to $1-\eta$ in order to ensure the well-posedness of the system and most of the numerical tests have been performed for $M_0=0.9$. The behavior of our optimal penalty model when we approach the characteristic boundary case (\emph{i.e.} $|M_0|=1$) is an interesting point and a natural question.

So, using the code for the two sides penalization, we tested the values $M_0=0.9$ (see Subsection \ref{Subsec_2sides_penal}, Fig. \ref{Error_newpen_regular_2bords}), $M_0=0.99$, $M_0=0.999$ and $M_0=0.9999$ (\emph{i.e.} $\eta=0.1, 10^{-2}, 10^{-3}, 10^{-4}$).

\begin{figure}
\begin{center}

\subfigure[$L^1$ error for $N$ in the plasma ($+$), $N$ in the limiter ($\times$), $\partial_x N$ in the plasma ($\circ$) and $\partial_x N$ in the limiter ($*$). $M_0=0.99$ ($\eta=10^{-2}$).]{\includegraphics[scale=0.56, trim = 7.5mm 3mm 19mm 13.3mm, clip=true]{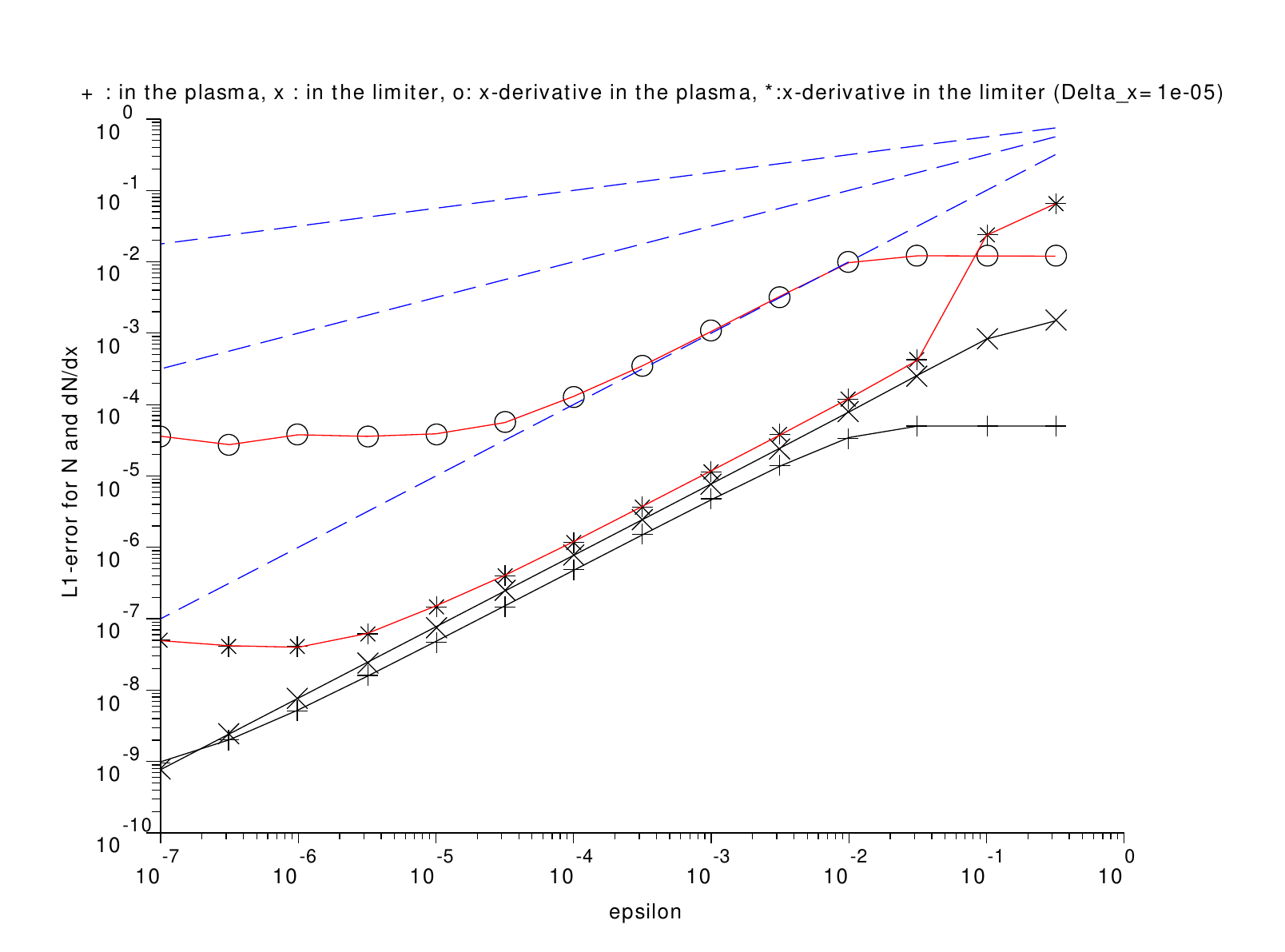}} $\quad$
\subfigure[$L^1$ error for $N$ in the plasma ($+$), $N$ in the limiter ($\times$), $\partial_x N$ in the plasma ($\circ$) and $\partial_x N$ in the limiter ($*$). $M_0=0.999$ ($\eta=10^{-3}$).]{\includegraphics[scale=0.56, trim = 7.5mm 3mm 19mm 13.3mm, clip=true]{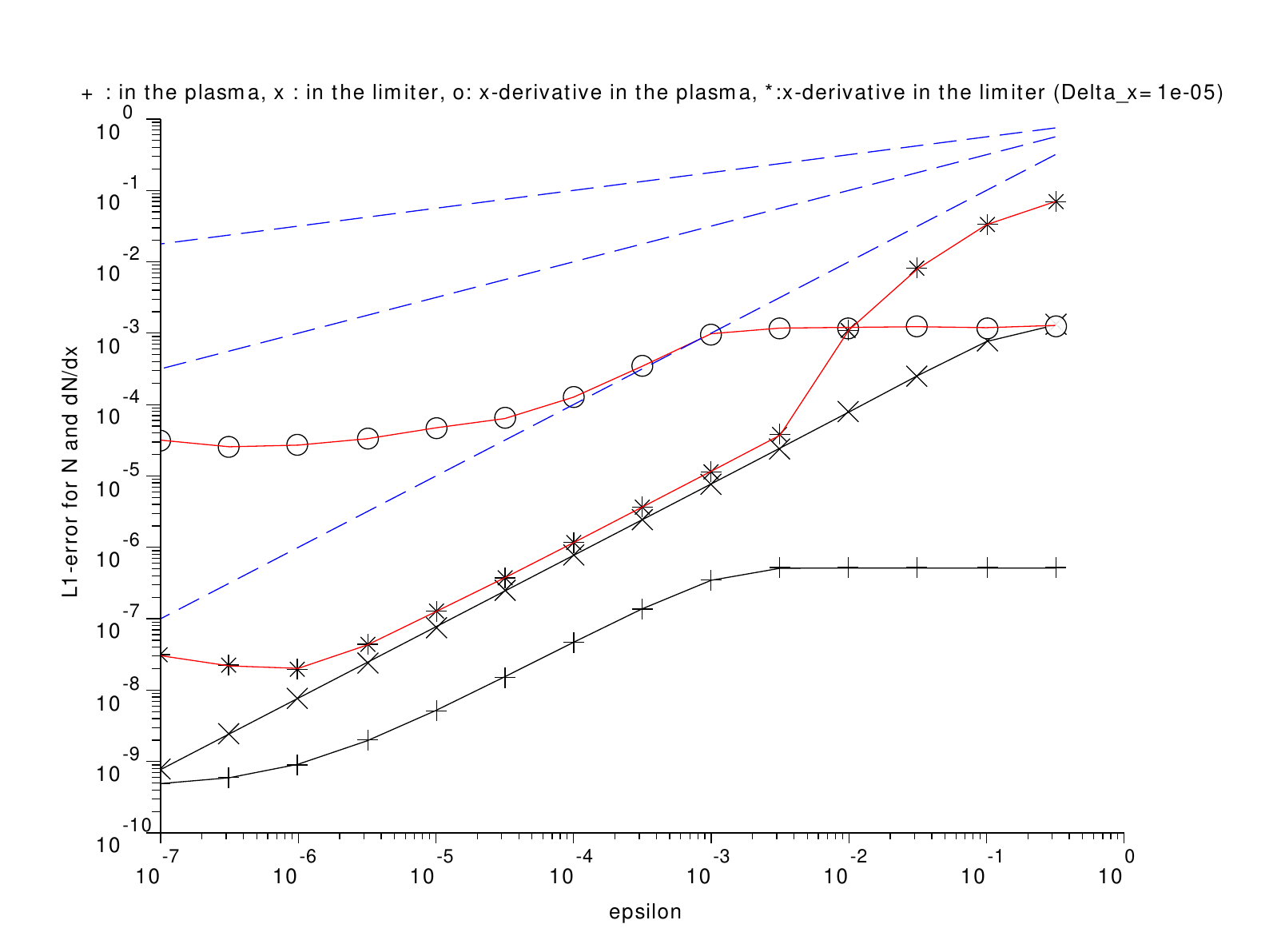}}\\
\subfigure[$L^1$ error for $N$ in the plasma ($+$), $N$ in the limiter ($\times$), $\partial_x N$ in the plasma ($\circ$) and $\partial_x N$ in the limiter ($*$). $M_0=0.9999$ ($\eta=10^{-4}$)]{\includegraphics[scale=0.56, trim = 7.5mm 3mm 19mm 13.3mm, clip=true]{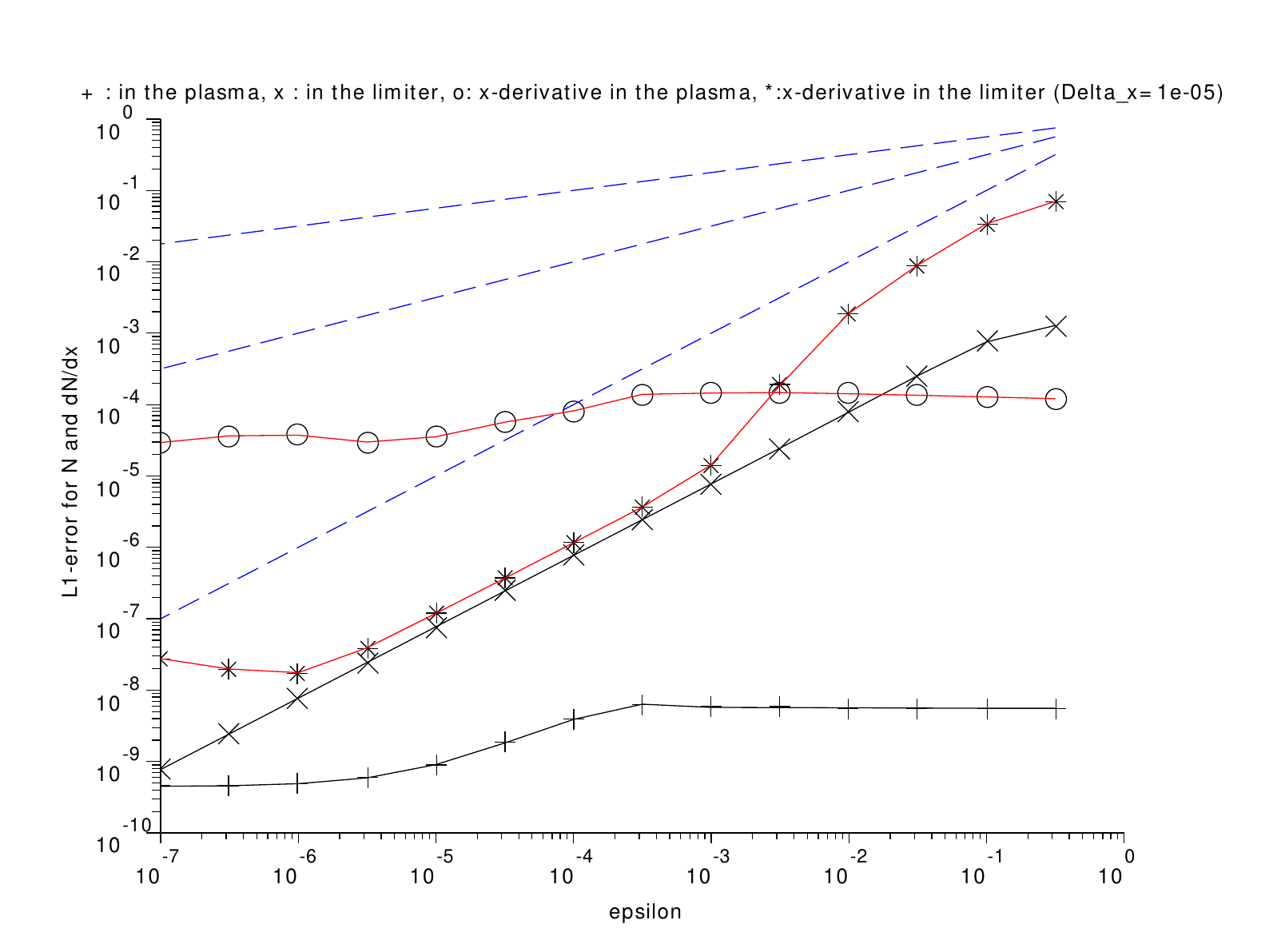}}

\end{center}
\caption{
Errors for $N$, $\partial_x N$, in $L^1$ norm versus $\varepsilon$ with the boundary layer free penalization and the two sides limiter configuration (see Fig. \ref{Complete_domain_center}). The dashed lines represent the curves $\varepsilon^{1/4}, \varepsilon^{1/2}$ and $\varepsilon$.
The error in the plasma area is estimated in the set $x \in ]-0.5,-0.1[\cup]0.1,0.5[$. In the limiter, the place where the error is computed is $x \in [-0.1,-0.075]\cup[0.075,0.1]$. So, the part where $\alpha(x)\neq 1$ is excluded.
}
\label{Error_newpen_regular_2bords_M_0_var}
\end{figure}

The computations show that (see Fig. \ref{Error_newpen_regular_2bords_M_0_var}), for $\varepsilon$ sufficiently small, such that $\varepsilon \leq \mathcal{O}(\eta)$, the convergence results are similar. This condition may come from the fact that $|M|$ must be less than $1$ and that the penalization error on $M$ is of the order of $\varepsilon$, see the first order of the asymptotic expansion in the Subsection \ref{asympt_exp}, equation (\ref{U_0_plus_eq_V1}).

From the theoretical point of view, the limit to the characteristic boundary case of the system (\ref{O2}) remains an open question. The simulations let us think that the numerical solution of the system (\ref{Penal_pb_2faces}) would converge when $\eta$ tends to $0$, but this might be due to the diffusivity of the scheme.

\section{Conclusion}
A nonlinear hyperbolic initial boundary problem has been studied in this paper. The set of equations is a simplified representation of the parallel plasma transport in the scrape-off layer of a tokamak. An interesting way to take into account the presence of obstacles such as limiter in the tokamak consists in using penalty method.

We first remark that the well-posedness of the hyperbolic problem (\ref{O1}) is not guaranteed. In order to ensure the well-posedness, we slightly modify the boundary condition on the Mach number $M$ and we do not impose $N=0$ at the boundary of the plasma.

In order to approximate the hyperbolic boundary value problem, we propose a carefully chosen penalty method which does not generate any artificial boundary layer: the convergence to the wished boundary value problem is sharp. This is in contrast with the previous approaches already known about this system. This is confirmed by our numerical tests which show an optimal rate of convergence in $\mathcal{O}(\varepsilon)$,
and also by an asymptotic analysis at any order of regular solutions.

This work has to be extended to a more complete model dealing with the plasma density, the momentum, the energy or the temperature and the electrical current. A penalization of the heat equation has been proposed by Paredes \emph{et al.} \cite{Par13} where the equations for $N$ and $\Gamma$ use the two-fields penalty methods described in Section \ref{Subsec 2 field penal}. The results from the future simulation codes are expected to provide a better understanding of the wall-plasma interactions in a tokamak and, perhaps, enable to fit the shape of the tokamak.

\subsection{Acknowledgements}
This work has been funded by the ANR ESPOIR (Edge Simulation of the Physics Of ITER Relevant turbulent transport) and the \emph{F\'ed\'eration nationale de Recherche sur la Fusion par Confinement Magn\'etique} (FR-FCM). We thank Guillaume Chiavassa, Guido Ciraolo and Philippe Ghendrih for fruitful discussions.

\bibliographystyle{plain}
\bibliography{biblio_article}

\begin{thebibliography}{10}

\bibitem{Ang05}
Ph. Angot.
\newblock A unified fictitious domain model for general embedded boundary
  conditions.
\newblock {\em Comptes Rendus Math\'ematique Acad. Sci. Paris}, 341(11):683 --
  688, 2005.

\bibitem{Ang11_2}
Ph. Angot, Ph. Auphan, and O.~Guès.
\newblock Penalty methods for the hyperbolic system modelling the wall-plasma
  interaction in a tokamak.
\newblock In {\em Finite Volumes for Complex Applications VI - Problems \&
  Perspectives}, volume~1, pages 31--38. Springer, June 2011.

\bibitem{Ang99}
Ph. Angot, Ch.-H. Bruneau, and P.~Fabrie.
\newblock A penalization method to take into account obstacles in an
  incompressible flow.
\newblock {\em Numerische Mathematik}, 81(4):497--520, 1999.

\bibitem{Aup10}
T.~Auphan.
\newblock Méthodes de pénalisation pour des systèmes hyperboliques et
  application au transport de plasma en bord de tokamak.
\newblock Master's thesis, Ecole Centrale Marseille, 2010.
\newblock Internship report.

\bibitem{Aup13}
T.~Auphan.
\newblock Penalization for non-linear hyperbolic system.
\newblock {\em Advances in Differential Equations}, 19(1/2):1--29, 2014.

\bibitem{Ben07}
S.~Benzoni-Gavage and D.~Serre.
\newblock {\em Multidimensional hyperbolic partial differential equations.
  First-order systems and applications}.
\newblock Oxford Mathematical Monographs. Oxford University Press, 2007.

\bibitem{Bou98}
F.~Bouchut and F.~James.
\newblock One-dimensional transport equations with discontinuous coefficients.
\newblock {\em Nonlinear Anal.}, 32:891--933, June 1998.

\bibitem{Boy12}
F.~Boyer and P.~Fabrie.
\newblock {\em Mathematical Tools for the Study of the Incompressible
  Navier-Stokes Equations and Related Models}.
\newblock Applied mathematical sciences. Springer, 2012.

\bibitem{Buf00}
T.~Buffard, T.~Gallou{\"e}t, and J.-M. H{\'e}rard.
\newblock A sequel to a rough {G}odunov scheme: application to real gases.
\newblock {\em Computers and Fluids}, 29(7):813 -- 847, 2000.

\bibitem{Buf98}
T.~Buffard, T.~Gallou{\"e}t, and J-M. Hérard.
\newblock Un schéma simple pour les équations de {S}aint-{V}enant.
\newblock {\em Comptes Rendus de l'Académie des Sciences - Series I -
  Mathematics}, 326(3):385 -- 390, 1998.

\bibitem{Car03}
G.~Carbou and P.~Fabrie.
\newblock Boundary layer for a penalization method for viscous incompressible
  flow.
\newblock {\em Differential Equations}, 8(12):1453--1480, 2003.

\bibitem{For08}
B.~Fornet.
\newblock Small viscosity solution of linear scalar 1-d conservation laws with
  one discontinuity of the coefficient.
\newblock {\em Comptes Rendus Mathematique}, 346(11-12):681 -- 686, 2008.

\bibitem{For09}
B.~Fornet and 0.~Guès.
\newblock Penalization approach of semi-linear symmetric hyperbolic problems
  with dissipative boundary conditions.
\newblock {\em Discrete and Continuous Dynamical Systems}, 23(3):827 -- 845,
  2009.

\bibitem{Gal03}
T.~Gallouët, J-M. Hérard, and N.~Seguin.
\newblock Some approximate {G}odunov schemes to compute shallow-water equations
  with topography.
\newblock {\em Computers and Fluids}, 32(4):479 -- 513, 2003.

\bibitem{Ghe11}
Ph. Ghendrih, K.~Bodi, H.~Bufferand, G.~Chiavassa, G.~Ciraolo, N~Fedorczak,
  L.~Isoardi, A.~Paredes, Y.~Sarazin, E.~Serre, F.~Schwander, and P.~Tamain.
\newblock Transition to supersonic flows in the edge plasma.
\newblock {\em Plasma Physics and Controlled Fusion}, 53(5):054019, 2011.

\bibitem{God96}
E.~Godlewski and P.-A. Raviart.
\newblock {\em Numerical approximation of hyperbolic systems of conservation
  laws}.
\newblock Springer, 1996.

\bibitem{Gre96}
J.~M. Greenberg and A.~Y. Le~Roux.
\newblock A well balanced scheme for the numerical processing of source terms
  in hyperbolic equation.
\newblock {\em J. Numer. Anal.}, 33(1):1--16, 1996.

\bibitem{Gue90}
0.~Guès.
\newblock Problème mixte hyperbolique quasi-linéaire caractéristique.
\newblock {\em Communications in Partial Differential Equations}, 15:595--654,
  1990.

\bibitem{Iso10}
L.~Isoardi, G.~Chiavassa, G.~Ciraolo, P.~Haldenwang, E.~Serre, Ph. Ghendrih,
  Y.~Sarazin, F.~Schwander, and P.~Tamain.
\newblock Penalization modeling of a limiter in the tokamak edge plasma.
\newblock {\em Journal of Computational Physics}, 229(6):2220 -- 2235, 2010.

\bibitem{Lev02}
R.~Leveque.
\newblock {\em Finite Volume Methods for Hyperbolic Problems}.
\newblock Cambridge University Press, 2002.

\bibitem{Liu07}
Qianlong Liu and Oleg~V. Vasilyev.
\newblock A {B}rinkman penalization method for compressible flows in complex
  geometries.
\newblock {\em Journal of Computational Physics}, 227(2):946 -- 966, 2007.

\bibitem{Pac05}
A.~Paccou, G.~Chiavassa, J.~Liandrat, and K.~Schneider.
\newblock A penalization method applied to the wave equation.
\newblock {\em Comptes Rendus Mécanique}, 333(1):79 -- 85, 2005.

\bibitem{Par13}
A.~Paredes, H.~Bufferand, F.~Schwander, G.~Ciraolo, E.~Serre, Ph. Ghendrih, and
  P.~Tamain.
\newblock Penalization technique to model wall-component impact on heat and
  mass transport in the tokamak edge.
\newblock {\em Journal of Nuclear Materials}, 438, Supplement(0):--, 2013.

\bibitem{Pou97}
F.~Poupaud and M.~Rascle.
\newblock Measure solutions to the linear multi-dimensional transport equation
  with non-smooth coefficients.
\newblock {\em Communications in Partial Differential Equations}, 22:225--267,
  1997.

\bibitem{Rau85}
J.~B. Rauch.
\newblock Symmetric positive systems with boundary characteristic of constant
  multiplicity.
\newblock {\em Trans. Amer. Math. Soc.}, 291(1):167--187, 1985.

\bibitem{Rau74}
J.~B. Rauch and F.~J.~III Massey.
\newblock Differentiability of solutions to hyperbolic initial-boundary value
  problems.
\newblock {\em Trans. Amer. Math. Soc.}, 189:303--318, 1974.

\bibitem{Sch00}
H.~Schlichting and K.~Gersten.
\newblock {\em Boundary Layer Theory}.
\newblock Physic and astronomy. MacGraw-Hill, 2000.

\bibitem{Tam07}
P.~Tamain.
\newblock {\em Etude des flux de matière dans le plasma de bord des tokamaks,
  alimentation, transport et turbulence}.
\newblock PhD thesis, Université de Provence, 2007.

\bibitem{Van79}
B~van Leer.
\newblock Towards the ultimate conservative difference scheme. v. a
  second-order sequel to {G}odunov's method.
\newblock {\em Journal of Computational Physics}, 32(1):101 -- 136, 1979.

\end{thebibliography}

\end{document}